\documentclass[10pt]{article}
\usepackage{amsfonts,amscd,amssymb}
\usepackage{theorem}
\input{diagrams}
\newtheorem{definition}{Definition}[section]
\newtheorem{lemma}[definition]{Lemma}
\newtheorem{proposition}[definition]{Proposition}

{\theorembodyfont{\rmfamily}\newtheorem{remark}[definition]{Remark}}
{\theorembodyfont{\rmfamily}}
\newtheorem{theorem}[definition]{Theorem}
{\theorembodyfont{\rmfamily}}
{\theorembodyfont{\rmfamily}}
{\theorembodyfont{\rmfamily}}
\def\va{\varepsilon}
\def\v{\varphi}

\def\tr{\triangleright}
\def\rh{\rightharpoonup}
\def\lh{\leftharpoonup}

\def\ra{\rightarrow}
\def\a{\alpha}
\def\b{\beta}

\def\l{\lambda}
\def\r{\rho}
\def\cd{\cdot}
\def\d{\delta}
\def\ov{\overline}
\def\un{\underline}
\newcommand{\tpra}{\mbox{$\tilde {p}^1$}}
\newcommand{\tprb}{\mbox{$\tilde {p}^2$}}

\newcommand{\tqra}{\mbox{$\tilde {q}^1$}}
\newcommand{\tqrb}{\mbox{$\tilde {q}^2$}}
\newcommand{\tQra}{\mbox{$\widetilde {Q}^1$}}
\newcommand{\tQrb}{\mbox{$\widetilde {Q}^2$}}

\newcommand{\smi}{\mbox{$S^{-1}$}}

\def\rawo\lonra{\longrightarrow}

\def\ot{\otimes}

\newcommand{\selabel}[1]{\label{se:#1}}

\newcommand{\eqref}[1]{(\ref{eq:#1})}

\newenvironment{proof}{{\it Proof.}}{\hfill $ \square $ \vskip 4mm}
\begin{document}
\title{Factorizable quasi-Hopf algebras. Applications}
\author{Daniel Bulacu
\thanks{During the time of the preparation of this work the first author was 
financially supported by INDAM and SB 2002-0286. He would like to thank 
University of Ferrara (Italy) and University of Almeria (Spain) for their warm 
hospitality.}\\
Faculty of Mathematics and Informatics\\
University of Bucharest\\ 
RO-010014 Bucharest 1, Romania
\and 
Blas Torrecillas\\
Department of Algebra and Analysis\\
University of Almeria\\
04071 Almeria, Spain}
\date{}
\maketitle

\begin{abstract}
We define the notion of factorizable quasi-Hopf algebra by using
a categorical point of view. We show that
the Drinfeld double $D(H)$ of any finite dimensional quasi-Hopf algebra
$H$ is factorizable, and we characterize $D(H)$ when $H$ itself is
factorizable. Finally, we prove that any finite dimensional factorizable 
quasi-Hopf algebra is unimodular. In particular, we obtain 
that the Drinfeld double $D(H)$ is a unimodular quasi-Hopf algebra. 
\end{abstract}
\section{Introduction}\selabel{0}
The concept of a quasi-triangular (or braided) bialgebra
is due to Drinfeld \cite{d2}. Roughly speaking, a bialgebra $H$
is quasi-triangular if the monoidal category of left
$H$-modules is braided in the sense of Joyal and Street \cite{js}.
In other words, $H$ is quasi-triangular if there exists an
invertible element $R\in H\ot H$ satisfying some
additional relations (see the complete definition below).
In the Hopf algebra case a reformulation of this definition
was given by Radford \cite{rad}. Ribbon and
factorizable Hopf algebras are special
classes of quasi-triangular Hopf algebras, and in this theory a
particular interest is produced by the Drinfeld double $D(H)$.
By the Drinfeld double construction \cite{d2},
every finite dimensional Hopf algebra $H$ can be embedded into a finite
dimensional quasi-triangular Hopf algebra $D(H)$.

As we pointed out, factorizable Hopf algebras belong to the class of
quasi-triangular Hopf algebras.
Suppose that $(H, R)$ is a quasi-triangular Hopf algebra
and denote by $R^1\ot R^2$ and $r^1\ot r^2$
two copies of the $R$-matrix $R$ of $H$. Then $(H, R)$ is
called factorizable if
$$
{\cal Q} : H^*\ra H,~~
{\cal Q}(\chi )=\chi (R^2r^1)R^1r^2~~
\forall~~\chi \in H^{*}
$$
is a linear isomorphism or, equivalently, if the map
$$
\ov{\cal Q} : H^*\ra H,~~
\ov{\cal Q}(\chi )=\chi (R^1r^2)R^2r^1~~
\forall~~\chi \in H^{*}
$$
is a linear isomorphism. Factorizable Hopf algebras were introduced and
studied by Reshetikhin and M. A. Semenov-Tian-Shansky \cite{rs}.
They are important in the Hennings investigation of $3$-manifold
invariants \cite{he}. Hennings shows how we can
construct $3$-manifold invariants using some finite
dimensional ribbon Hopf algebras which are, in particular,
factorizable. Afterwards, Kauffman reworked the Hennings construction,
see \cite{kauf} or \cite{rad1} for more details. We note that factorizable
Hopf algebras are also important in the representation theory
\cite{sch}, notably with applications to the classification of a
certain classes of Hopf algebras, see \cite{eg}.

Now, quasi-bialgebras and quasi-Hopf algebras were introduced by
Drinfeld \cite{d1}. They come out from categorical considerations: putting
some additional structure on the category of modules over an algebra $H$,
the definition of a quasi-bialgebra $H$ ensures that
the category of left $H$-modules ${}_H{\cal M}$ is a
monoidal category.
So $H$ is an unital associative algebra
together with a comultiplication $\Delta : H\ra H\ot H$ and a usual
counit $\va : H\ra k$ such that $\Delta $ and $\va $ are algebra maps,
and $\Delta $ is quasi-coassociative, in the sense that it is coassociative
up to conjugation by an invertible element $\Phi \in H\ot H\ot H$.
Consequently the definition of a quasi-bialgebra is not self dual.
For a quasi-Hopf algebra $H$
the definition ensures that the
category of finite dimensional left $H$-modules is a
monoidal category with left duality. In a similar manner
one can define quasi-triangular (ribbon, at least in the finite
dimensional case) quasi-bialgebras: a quasi-bialgebra
is called quasi-triangular (ribbon) if the monoidal  category
${}_H{\cal M}$ is braided (ribbon, respectively).
In the quasi-triangular case, this means that
there exists an invertible element $R\in H\ot H$ satisfying
some additional conditions. If $H$ is
a quasi-Hopf algebra then the definition of a
quasi-triangular quasi-bialgebra can be reformulated, see \cite{bn3}.
As we have already explained the study of quasi-Hopf algebras,
or quasi-triangular (ribbon) quasi-Hopf algebras
is strictly connected to the study
of monoidal, or braided (ribbon) categories. Thus, in general,
when we want to define some classes of quasi-Hopf algebras we should
look at the classical Hopf case in the sense that we should
try to reformulate their basic properties at a
categorical level, and then we must come back to the quasi-Hopf case;
if this is not possible then we have to be lucky in order to
define (and then study) them.

As far as we are concerned, in the classical case,
the definition of the map ${\cal Q}$ given
above has a categorical
interpretation due to Majid \cite{maj}. Hence,
if in the quasi-Hopf case
a suitable map satisfies the same
categorical interpretation then it makes sense to
define the factorizable notion.
This is way we propose in Section 2 the following
definition for the map ${\cal Q}$:
$$
{\cal Q}(\chi)=<\chi , S(X^2_2\tprb )f^1R^2r^1U^1X^3>
X^1S(X^2_1\tpra )f^2R^1r^2U^2,
$$
for all $\chi \in H^*$, where $r^1\ot r^2$ is another
copy of $R$, $\Phi =X^1\ot X^2\ot X^3=Y^1\ot Y^2\ot Y^3$, and
$f^1\ot f^2, U^1\ot U^2, \tpra \ot \tprb \in H\ot H$
are some special elements which we will define below. Moreover,
we will see that in the quasi-Hopf case
the analogues of the map $\ov{\cal Q}$ is
$$
\ov{\cal Q}(\chi )=<\chi , \smi (X^3)q^2R^1r^2X^2_2\tprb >
q^1R^2r^1X^2_1\tpra \smi (X^1),
$$
for all $\chi \in H^*$, where $q^1\ot q^2\in H\ot H$ is
another special element which will be defined.

Following \cite{maj}, in Section 4 we will give the categorical
interpretation of the map ${\cal Q}$ in the quasi-Hopf case.
For this, we develop first in Section 3 the transmutation
theory for dual quasi-Hopf algebras. Using the dual
reconstruction theorem (also due to Majid) we will show that
to any co-quasi-triangular dual quasi-Hopf algebra $A$ we can
associate a braided commutative Hopf algebra $\un{A}$ in the category
of right $A$-comodules. Keeping the same terminology as in the
Hopf case we will call $\un{A}$ the function algebra braided group
associated to $A$. This procedure is the formal dual of the one
performed in \cite{bn2} where to any quasi-triangular
quasi-Hopf algebra $H$ is associated a braided cocommutative
group $\un{H}$ in the braided category of left $H$-modules.
We call $\un{H}$ the associated enveloping algebra braided group of $H$.
We notice that, in the finite
dimensional case, $\un{A}$ cannot be obtained from $\un{H}$ by (usual)
dualisation. In fact, if $H$ is finite dimensional
then the map ${\cal Q}$ provides a braided Hopf algebra morphism between
the function algebra braided group $\un{H^*}$
associated to $H^*$ and $\un{H}$ (Proposition \ref{pr4.1}).
Moreover, $\un{H^*}$ is always
isomorphic to the categorical dual of $\un{H}$ as braided Hopf
algebra (Proposition \ref{pr4.2}).
So the true meaning of the map ${\cal Q}$ is that
$\un{H}$  and $\un{H^*}$ are self dual (in a categorical sense)
provided ${\cal Q}$ is bijective, i.e. $H$ is factorizable.

Let $H$ be a finite dimensional quasi-Hopf algebra and $D(H)$
the Drinfeld double of $H$. Similarly to the Hopf case
we will show in Section 2 that $D(H)$ is always factorizable.
The description of $D(H)$ when $H$ is
quasi-triangular was given in \cite{bc}. In this case $D(H)$
is a biproduct (in the sense of \cite{bn2}) of a braided Hopf
algebra $B^i$ and $H$, and, as a vector space, $B^i$ is isomorphic to
$H^*$. Furthermore, in Section 5 we will see that the
Drinfeld double $D(H)$ has a very simple description
when $H$ itself is factorizable. In fact,
we will give a quasi-Hopf version of a result claimed in \cite{rs}
and proved in \cite{sch}. We will show that if $H$ is factorizable then
$D(H)$ is isomorphic as a quasi-Hopf algebra to a twist of a usual
(componentwise) tensor product quasi-Hopf algebra $H\ot H$.
To this end we need the
alternative definition for the space of coinvariants
of a right quasi-Hopf $H$-bimodule and the
second structure theorem for right quasi-Hopf $H$-bimodules
proved in \cite{bc1}. Finally, in Section 6 we will prove that 
any finite dimensional factorizable quasi-Hopf algebra is unimodular. 
In particular, we deduce that the Drinfeld double $D(H)$ is 
a unimodular quasi-Hopf algebra. 

As we will see, the theory of quasi-Hopf algebras is technically
more complicated than the classical Hopf algebra theory.
This happens because of the nature
of the problems which occur in the quasi-Hopf algebra theory.
When we pass from Hopf algebras to quasi-Hopf algebras the
appearance of the reassociator $\Phi $ and of the 
elements $\a $ and $\b $ in the definition of the antipode 
increases the complexity of formulas, and therefore
of computations and proofs.
\section{Preliminaries}\selabel{1}
\subsection{Quasi-Hopf algebras}
We work over a commutative field $k$. All algebras, linear
spaces etc. will be over $k$; unadorned $\ot $ means $\ot_k$.
Following Drinfeld \cite{d1}, a quasi-bialgebra is
a four-tuple $(H, \Delta , \va , \Phi )$ where $H$ is
an associative algebra with unit,
$\Phi$ is an invertible element in $H\ot H\ot H$, and
$\Delta :\ H\ra H\ot H$ and $\va :\ H\ra k$ are algebra
homomorphisms satisfying the identities
\begin{eqnarray}
&&(id \ot \Delta )(\Delta (h))=
\Phi (\Delta \ot id)(\Delta (h))\Phi ^{-1},\label{q1}\\
&&(id \ot \va )(\Delta (h))=h\ot 1,
\mbox{${\;\;\;}$}
(\va \ot id)(\Delta (h))=1\ot h,\label{q2}
\end{eqnarray}
for all $h\in H$, and
$\Phi$ has to be a $3$-cocycle, in the sense that
\begin{eqnarray}
&&(1\ot \Phi)(id\ot \Delta \ot id)
(\Phi)(\Phi \ot 1)=
(id\ot id \ot \Delta )(\Phi )
(\Delta \ot id \ot
id)(\Phi ),\label{q3}\\
&&(id \ot \va \ot id)(\Phi )=1\ot 1\ot 1.\label{q4}
\end{eqnarray}
The map $\Delta $ is called the
coproduct or the
comultiplication, $\va $ the counit and $\Phi $ the
reassociator.
As for Hopf algebras we denote $\Delta (h)=h_1\ot h_2$,
but since $\Delta $ is only quasi-coassociative we adopt the
further convention (summation understood):
$$
(\Delta \ot id)(\Delta (h))=h_{(1, 1)}\ot h_{(1, 2)}\ot h_2,
\mbox{${\;\;\;}$}
(id\ot \Delta )(\Delta (h))=h_1\ot h_{(2, 1)}\ot h_{(2,2)},
$$
for all $h\in H$. We will
denote the tensor components of $\Phi$
by capital letters, and the ones of
$\Phi^{-1}$ by small letters,
namely
\begin{eqnarray*}
&&\Phi=X^1\ot X^2\ot X^3=T^1\ot T^2\ot T^3=
V^1\ot V^2\ot V^3=\cdots\\
&&\Phi^{-1}=x^1\ot x^2\ot x^3=t^1\ot t^2\ot t^3=
v^1\ot v^2\ot v^3=\cdots
\end{eqnarray*}
$H$ is
called a quasi-Hopf
algebra if, moreover, there exists an
anti-morphism $S$ of the algebra
$H$ and elements $\a , \b \in
H$ such that, for all $h\in H$, we
have:
\begin{eqnarray}
&&S(h_1)\a h_2=\va(h)\a
\mbox{${\;\;\;}$ and ${\;\;\;}$}
h_1\b S(h_2)=\va (h)\b,\label{q5}\\
&&X^1\b S(X^2)\a X^3=1
\mbox{${\;\;\;}$ and${\;\;\;}$}
S(x^1)\a x^2\b S(x^3)=1.\label{q6}
\end{eqnarray}
Our definition of a quasi-Hopf algebra is different to the
one given by Drinfeld \cite{d1} in the sense that we do not
require the antipode to be bijective. Nevertheless, in the finite
dimensional or quasi-triangular case the antipode is automatically bijective,
cf. \cite{bc1} and \cite{bn3}. This is way we omit the bijectivity of $S$ in
the definition of a quasi-Hopf algebra.\\
For a quasi-Hopf algebra the antipode is
determined uniquely up
to a transformation $\a \mapsto \a _{\mathbb{U}}:=\mathbb{U}\a $,
$\b \mapsto \b _{\mathbb{U}}:=\b \mathbb{U}^{-1}$,
$S(h)\mapsto S_{\mathbb{U}}(h):=\mathbb{U}S(h)\mathbb{U}^{-1}$,
where $\mathbb{U}\in H$ is invertible. In this case we will denote
by $H^{\mathbb{U}}$ the new quasi-Hopf algebra $(H, \Delta , \va, \Phi ,
S_{\mathbb{U}}, \a _{\mathbb{U}}, \b _{\mathbb{U}})$. \\
The axioms for a quasi-Hopf algebra imply that
$\va \circ S=\va $ and $\va (\a )\va (\b )=1$,
so, by rescaling $\a $ and $\b $, we may assume without loss of generality
that $\va (\a )=\va (\b )=1$. The identities
(\ref{q2}), (\ref{q3}) and (\ref{q4}) also imply that
\begin{equation}\label{q7}
(\va \ot id\ot id)(\Phi )=
(id \ot id\ot \va )(\Phi )=1\ot 1\ot 1.
\end{equation}
Next we recall that the definition of a
quasi-Hopf
algebra is "twist coinvariant" in the following sense.
An invertible element $F\in H\ot H$ is called a gauge
transformation or twist if
$(\va \ot id)(F)=(id\ot \va)(F)=1$.
If $H$ is a quasi-Hopf algebra and $F=F^1\ot F^2\in H\ot H$
is a gauge transformation with inverse $F^{-1}=G^1\ot
G^2$, then we can define a
new quasi-Hopf algebra $H_F$ by keeping
the multiplication, unit, counit
and antipode of $H$ and
replacing the comultiplication, reassociator and the
elements $\alpha$
and $\beta$ by
\begin{eqnarray}
&&\Delta _F(h)=F\Delta (h)F^{-1},\label{g1}\\
&&\Phi_F=(1\ot F)(id \ot \Delta )(F) \Phi (\Delta \ot
id)(F^{-1})(F^{-1}\ot 1),\label{g2}\\
&&\a_F=S(G^1)\a G^2,
\mbox{${\;\;\;}$}%
\b_F=F^1\b S(F^2).\label{g3}
\end{eqnarray}
It is well-known that the antipode of a Hopf algebra
is an
anti-coalgebra morphism. For a quasi-Hopf algebra, we have
the
following statement: there exists a gauge transformation
$f\in H\ot H$ such that
\begin{equation} \label{ca}
f\Delta (S(h))f^{-1}=(S\ot S)(\Delta
^{op}(h))
\mbox{,${\;\;\;}$for all $h\in H$}
\end{equation}
where $\Delta ^{op}(h)=h_2\ot h_1$. $f$ can be
computed explicitly. First set
\begin{eqnarray*}
A^1\ot A^2\ot A^3\ot A^4&=&
(\Phi \ot 1) (\Delta \ot id\ot id)(\Phi
^{-1}),\\
B^1\ot B^2\ot B^3\ot B^4&=&
(\Delta \ot id\ot id)(\Phi )(\Phi ^{-1}\ot 1)
\end{eqnarray*}
and then define $\gamma, \delta\in H\ot H$ by
\begin{equation} \label{gd}%
\gamma =S(A^2)\a A^3\ot S(A^1)\a A^4~~{\rm and}~~
\delta =B^1\b S(B^4)\ot B^2\b S(B^3).
\end{equation}
$f$ and $f^{-1}$ are then given by
the formulas
\begin{eqnarray}
f&=&(S\ot S)(\Delta ^{op}(x^1)) \gamma
\Delta (x^2\b
S(x^3)),\label{f}\\
f^{-1}&=&\Delta (S(x^1)\a x^2)
\delta (S\ot S)(\Delta ^{op}(x^3)).\label{g}
\end{eqnarray}
Moreover, $f$ satisfies the following relations:
\begin{equation} \label{gdf}
f\Delta (\a )=\gamma ,
\mbox{${\;\;\;}$}
\Delta (\b )f^{-1}=\delta .
\end{equation}
Furthermore the corresponding
twisted reassociator (see (\ref{g2})) is given by
\begin{equation} \label{pf}
\Phi _f=(S\ot S\ot S)(X^3\ot X^2\ot X^1).
\end{equation}
In a Hopf algebra $H$, we obviously have the
identity
$$
h_1\ot h_2S(h_3)=h\ot 1,~{\rm for~all~}h\in H.
$$
We will need the generalization of this formula to
quasi-Hopf algebras. Following
\cite{hn1}, \cite{hn2}, we define
\begin{equation} \label{qr}
p_R=p^1\ot p^2=x^1\ot x^2\b S(x^3),
\mbox{${\;\;\;}$} q_R=q^1\ot
q^2=X^1\ot S^{-1}(\a X^3)X^2.
\end{equation}
\begin{equation}\label{ql}
p_L=p^1_L\ot p^2_L=X^2\smi (X^1\b )\ot X^3,%
\mbox{${\;\;\;}$}%
q_L=q^1_L\ot q^2_L=S(x^1)\a x^2\ot x^3.%
\end{equation}
For all $h\in H$, we then have
\begin{equation} \label{qr1}
\Delta (h_1)p_R[1\ot S(h_2)]=p_R[h\ot 1],
\mbox{${\;\;\;}$}%
[1\ot S^{-1}(h_2)]q_R\Delta (h_1)=(h\ot 1)q_R,
\end{equation}
\begin{equation}\label{ql1}
\Delta (h_2)p_L[\smi (h_1)\ot 1]=p_L(1\ot h),%
\mbox{${\;\;\;}$}%
[S(h_1)\ot 1]q_L\Delta (h_2)=(1\ot h)q_L,
\end{equation}
and
\begin{equation} \label{pqr}
\Delta (q^1)p_R[1\ot S(q^2)]=1\ot 1,
\mbox{${\;\;\;}$}
[1\ot S^{-1}(p^2)]q_R\Delta (p^1)=1\ot 1,
\end{equation}
\begin{equation}\label{pql}
[S(p^1_L)\ot 1]q_L\Delta (p^2_L)=1\ot 1,%
\mbox{${\;\;\;}$}%
\Delta (q^2_L)p_L[\smi (q^1_L)\ot 1]=1\ot 1,%
\end{equation}
\begin{eqnarray}
&&\hspace*{-2cm}\Phi (\Delta \ot id)(p_R)(p_R\ot id)\nonumber\\
&=&(id\ot \Delta )(\Delta (x^1)p_R)(1\ot f^{-1})(1\ot
S(x^3)\ot S(x^2))\label{pr},\\
&&\hspace*{-2cm}
(1\ot q_L)(id \ot \Delta )(q_L)\Phi =
(S(x^2)\ot S(x^1)\ot 1)(f\ot 1)(\Delta \ot id)
(q_L\Delta (x^3)),\label{ql2}
\end{eqnarray}
where $f=f^1\ot f^2$ is the twist defined in (\ref{f}).\\
Note that some of the above formulas
use the bijectivity of the antipode $S$. Nevertheless,
we will restrict to apply them only in this case.
\subsection{Quasi-triangular quasi-Hopf algebras}
A quasi-Hopf algebra $H$ is quasi-triangular
(QT for short) if  there exists an element $R\in H\ot H$ such that
\begin{eqnarray}
(\Delta \ot id)(R)&=&\Phi _{312}R_{13}
\Phi ^{-1}_{132}R_{23}\Phi ,\label{qt1}\\
(id \ot \Delta )(R)&=&\Phi ^{-1}_{231}R_{13}
\Phi _{213}R_{12}\Phi ^{-1},\label{qt2}\\
\Delta ^{\rm op}(h)R&=&R\Delta (h),~{\rm for~all~}h\in H,\label{qt3}\\
(\va \ot id)(R)&=&(id\ot \va)(R)=1,\label{qt4}
\end{eqnarray}
where, if $\sigma $ denotes a permutation of $\{1, 2, 3\}$, we set $\Phi _{\sigma (1)
\sigma (2)\sigma (3)}=X^{\sigma ^{-1}(1)}\ot
X^{\sigma ^{-1}(2)}\ot X^{\sigma ^{-1}(3)}$, and $R_{ij}$ means $R$ acting non-trivially
in the $i^{th}$ and $j^{th}$ positions of $H\ot H\ot H$.\\
In {\cite{bn3}} it is shown that, consequently, $R$ is
invertible. Furthermore, the element
\begin{equation} \label{elmu}
u=S(R^2p^2)\a R^1p^1
\end{equation}
(with $p_R=p^1\ot p^2$ defined as in (\ref{qr})) is
invertible in $H$, and
\begin{equation} \label{inelmu}
u^{-1}=X^1R^2p^2S(S(X^2R^1p^1)\a X^3),
\end{equation}
\begin{equation} \label{sqina}
\va (u)=1~~{\rm and}~~
S^2(h)=uhu^{-1}
\end{equation}
for all $h\in H$. Consequently the antipode $S$ is bijective, so,
as in the Hopf algebra case, the assumptions about invertibility
of $R$ and bijectivity of $S$ can be dropped. Moreover, the
$R$-matrix $R=R^1\ot R^2$ satisfy the identities (see \cite{ac},
\cite{hn2}, \cite{bn3}):
\begin{eqnarray}
&&f_{21}Rf^{-1}=(S\ot S)(R), \label{ext}\\
&&S(R^2)\a R^1=S(\a )u.\label{exta}
\end{eqnarray}
\subsection{Hopf algebras in braided categories}
For further use we briefly recall some concepts concerning braided categories
and braided Hopf algebras. For more details the reader is invited to consult
\cite{Kassel95} or \cite{maj}.

A monoidal category means a category ${\cal C}$ with objects
$U, V, W$ etc., a functor $\ot : {\cal C}\times {\cal C}\ra {\cal C}$
equipped with an associativity natural transformation consisting
of functorial isomorphisms
$a_{U, V, W}: (U\ot V)\ot W\ra U\ot (V\ot W)$ satisfying a pentagon identity,
and a compatible unit object $\un {1}$ and associated functorial isomorphisms
(the left and right unit constraints, $l_V: V\cong V\ot \un{1}$ and
$r_V: V\cong \un{1}\ot V$, respectively).

Let ${\cal C}$ be a monoidal category.
An object $V\in {\cal C}$ has a left dual
or is left rigid if there is
an object $V^*$ and morphisms
$ev_V: V^*\ot V\ra \un{1}$,
$coev_V: \un{1}\ra V\ot V^*$ such that
\begin{equation}\label{rig1}
l^{-1}_V\circ (id_V\ot ev_V)
\circ a_{V, V^*, V}\circ (coev_V\ot id_V)\circ r_V=id_V,
\end{equation}
\begin{equation}\label{rig2}
r^{-1}_{V^*}\circ (ev_V\ot id_{V^*})
\circ a^{-1}_{V^*, V, V^*}\circ
(id_{V^*}\ot coev_V)\circ l_{V^*}=id_{V^*}.
\end{equation}
If every object in the category has a left dual,
then we say that ${\cal C}$ is a left rigid
monoidal category or a monoidal category with left duality.

A braided category is a monoidal category equipped with a
commutativity natural transformation consisting of functorial
isomorphisms $c_{U, V}: U\ot V\ra V\ot U$
compatible with the unit and the associativity
structures in a natural way (for a complete definition see
\cite{Kassel95}, \cite{maj}).

Suppose that $(H, \Delta, \va , \Phi)$ is a
quasi-bialgebra. If $U, V, W$ are left $H$-modules,
define $a_{U, V, W}: (U\ot V)\ot W\ra U\ot (V\ot W)$
by
\begin{equation}\label{as}
a_{U, V, W}((u\ot v)\ot w)=
\Phi \cd (u\ot (v\ot w)).
\end{equation}
Then the category ${}_H{\cal M}$ of left $H$-modules becomes a monoidal
category with tensor product $\ot $ given via $\Delta $, associativity
constraints $a_{U, V, W}$, unit $k$ as a trivial $H$-module and the usual
left and right unit constraints.
If $H$ is a quasi-Hopf algebra then the
category of finite dimensional left
$H$-modules is left rigid; the left dual of $V$ is
$V^*$ with the $H$-module structure given by
$(h\cd v^*)(v)=v^*(S(h)\cd v)$, for all $v\in V$,
$v^*\in V^*$, $h\in H$ and with
\begin{equation}\label{qrig}
ev_V(v^*\ot v)=v^*(\a \cd v),
\mbox{${\;\;\;}$}
coev_V(1)=\b \cd {}_iv\ot {}^iv,
\end{equation}
where $\{{}_iv\}$ is a basis in $V$ with dual basis $\{{}^iv\}$.
Now, if $H$ is a QT quasi-Hopf algebra with
$R$-matrix $R=R^1\ot R^2$, then ${}_H{\cal M}$ is a
braided category with braiding
\begin{equation}\label{br}
c_{U, V}(u\ot v)=R^2\cd v\ot R^1\cd u,
\mbox{${\;}$ for all $u\in U$, $v\in V$.}
\end{equation}
Finally, the definition of a Hopf
algebra $B$ in a braided category ${\cal C}$ is obtained
in the obvious way in analogy
with the standard axioms \cite{sw}. Thus, a bialgebra in ${\cal C}$
is $(B, m_B, \eta _B, \Delta _B, \va _B)$ where $B$ is an object in ${\cal C}$
and the morphism $m_B:\ B\ot B\ra B$ forms a multiplication
that is associative up to the isomorphism $a$.
Similarly for the coassociativity of the
comultiplication $\Delta _B: B\ra B\ot B$.
The identity in the algebra $B$ is expressed as usual
by $\eta _B:\ \un {1}\ra B$ such that
$m_B\circ (\eta _B\ot id)=m_B\circ (id \ot \eta _B)=id$.
The counit axiom is $(\va _B\ot id)\circ \Delta _B=
(id \ot \va _B)\circ \Delta _B=id$.
In addition, $\Delta _B$ is required to be an algebra morphism where
$B\ot B$ has the multiplication
\begin{eqnarray*}
m_{B\ot B}:&&\hspace*{-3mm}
(B\ot B)\ot (B\ot B)\stackrel{a}{\longrightarrow}B\ot (B\ot (B\ot B))
\stackrel{id\ot a^{-1}}{\longrightarrow}B\ot ((B\ot B)\ot B)\\
            &&\hspace*{-3mm}
\stackrel{id\ot (c\ot id)}{\longrightarrow}B\ot ((B\ot B)\ot B)
\stackrel{id\ot a}{\longrightarrow}B\ot (B\ot (B\ot B))
            \stackrel{a^{-1}}{\longrightarrow}(B\ot B)\ot (B\ot B)
\stackrel{m_B\ot m_B}{\longrightarrow}B\ot B.
\end{eqnarray*}
A Hopf algebra $B$ is a bialgebra with a
morphism $S_B: B\ra B$ in ${\cal C}$ (the antipode)
satisfying the usual axioms
$m_B\circ (S_B\ot id)\circ \Delta _B=\eta _B\circ \va _B=
m_B\circ (id \ot S_B)\circ \Delta _B$.
\section{Factorizable QT quasi-Hopf algebras}\selabel{2}
\setcounter{equation}{0}
In this Section we will introduce the notion of
factorizable quasi-Hopf algebra and we will show
that the quantum double is an example of this type.\\
If $(H, R)$ is a QT quasi-Hopf algebra then we will see that
the $k$-linear map ${\cal Q}:\ H^*\ra H$, given
for all $\chi \in H^*$ by
\begin{equation}\label{qf1}
{\cal Q}(\chi)=<\chi , S(X^2_2\tprb )f^1R^2r^1U^1X^3>
X^1S(X^2_1\tpra )f^2R^1r^2U^2,
\end{equation}
where $r^1\ot r^2$ is another copy of $R$ and
\begin{equation}\label{eu}
U=g^1S(q^2)\ot g^2S(q^1)
\end{equation}
(here $f^{-1}=g^1\ot g^2$ and $q_R=q^1\ot q^2$
are the elements defined by (\ref{g}) and (\ref{qr}), respectively)
has most of the properties satisfied by the map
$H^*\ni \chi \mapsto <\chi , R^2r^1>R^1r^2\in H$ defined
in the Hopf case. For this reason we will propose the
following.

\begin{definition}
A QT quasi-Hopf algebra $(H, R)$ is called factorizable if the
map ${\cal Q}$ defined by (\ref{qf1}) is bijective.
\end{definition}

We will see in the next Section that the above definition
has a categorical explanation. In fact in this way we were
able to found a suitable definition for the map ${\cal Q}$.\\
In the sequel we will need a second formula for
the map ${\cal Q}$. Also, another $k$-linear map
$\ov {\cal Q}:\  H^*\ra H$ is required.

\begin{proposition}\label{pr2.2}
Let $(H, R)$ be a QT quasi-Hopf algebra.
\begin{itemize}
\item[i)] The map ${\cal Q}$ defined by (\ref{qf1})
has a second formula given for all $\chi \in H^*$ by
\begin{equation}\label{qf2}
{\cal Q}(\chi )=<\chi , \tqra X^1R^2r^1p^1>
\tilde{q}^2_1X^2R^1r^2p^2S(\tilde{q}^2_2X^3),
\end{equation}
where $q_L=\tqra \ot \tqrb $ and $p_R=p^1\ot p^2$
are the elements defined by (\ref{ql}) and (\ref{qr}), respectively.
\item[ii)] Let $\ov {\cal Q}:\  H^*\ra H$ be the $k$-linear map
defined for all $\chi \in H^*$ by
\begin{equation}\label{oqf}
\ov{\cal Q}(\chi )=<\chi , \smi (X^3)q^2R^1r^2X^2_2\tprb >
q^1R^2r^1X^2_1\tpra \smi (X^1),
\end{equation}
where $q_R=q^1\ot q^2$ and $p_L=\tpra \ot \tprb $ are the elements
defined by (\ref{qr}) and (\ref{ql}), respectively.
Then ${\cal Q}$ is bijective if and only if $\ov{\cal Q}$ is bijective.
\end{itemize}
\end{proposition}
\begin{proof}
i) It is not hard to see that (\ref{q3}) and (\ref{q7}) imply
\begin{equation}\label{f1tp}
X^1S(X^2_1\tpra )\ot X^2_2\tprb \ot X^3=
S(x^1\tpra )\ot x^2\tilde{p}^2_1\ot x^3\tilde{p}^2_2.
\end{equation}
We need also the formulas
\begin{eqnarray}
&&p_R=\Delta (S(\tpra ))U(\tprb \ot 1),\label{uptp}\\
&&U^1\ot U^2S(h)=S(h_1)_1U^1h_2\ot S(h_1)_2U^2\label{eu1}
\end{eqnarray}
which can be found in \cite{hn3}. Now, we claim that
\begin{equation}\label{ru}
R^1U^1\ot R^2U^2=\tilde{q}^1_2R^1p^1\ot \tilde{q}^1_1R^2p^2S(\tqrb ).
\end{equation}
Indeed, we calculate:
\begin{eqnarray*}
\tilde{q}^1_2R^1p^1\ot \tilde{q}^1_1R^2p^2S(\tqrb )
&\pile{{\rm (\ref{qt3})}\\ =}&
R^1\tilde{q}^1_1p^1\ot R^2\tilde{q}^1_2p^2S(\tqrb )\\
{\rm (\ref{uptp})}&=&R^1(\tqra S(\tpra ))_1U^1\tprb \ot
R^2(\tqra S(\tpra ))_2U^2S(\tqrb )\\
{\rm (\ref{eu1})}&=&R^1(\tqra S(\tilde{q}^2_1\tpra ))_1U^1
\tilde{q}^2_2\tprb \ot
R^2(\tqra S(\tilde{q}^2_1\tpra ))_2U^2\\
{\rm (\ref{pql})}&=&R^1U^1\ot R^2U^2,
\end{eqnarray*}
as needed. Now, if we denote by $\tQra \ot \tQrb $ another
copy of $q_L$ then for all $\chi \in H^*$ we have
\begin{eqnarray*}
{\cal Q}(\chi)&\pile{{\rm (\ref{qf1}, \ref{f1tp})}\\ =}&
<\chi , S(x^2\tilde{p}^2_1)f^1R^2r^1U^1x^3\tilde{p}^2_2>
S(x^1\tpra )f^2R^1r^2U^2\\
{\rm (\ref{ru}, \ref{qr1}, \ref{qt3})}&=&
<\chi , S(x^2\tilde{p}^2_1)f^1\tilde{q}^1_1(x^3\tilde{p}^2_2)_{(1, 1)}
R^2r^1p^1>\\
&&S(x^1\tpra )f^2\tilde{q}^1_2(x^3\tilde{p}^2_2)_{(1, 2)}
R^1r^2p^2S(\tqrb (x^3\tilde{p}^2_2)_2)\\
{\rm (\ref{ql2})}&=&<\chi , S(\tilde{p}^2_1)
\tQra X^1(\tilde{p}^2_2)_{(1, 1)}R^2r^1p^1>\\
&&S(\tpra )\tqra \tilde{Q}^2_1X^2(\tilde{p}^2_2)_{(1, 2)}
R^1r^2p^2S(\tqrb \tilde{Q}^2_2X^3(\tilde{p}^2_2)_2)\\
{\rm (\ref{q1}, \ref{ql1})}&=&
<\chi , \tQra X^1R^2r^1p^1>S(\tpra )\tqra \tilde{p}^2_1
\tilde{Q}^2_1X^2R^1r^2p^2
S(\tqrb \tilde{p}^2_2\tilde{Q}^2_2X^3)\\
{\rm (\ref{pql})}&=&
<\chi , \tQra X^1R^2r^1p^1>\tilde{Q}^2_1X^2R^1r^2p^2
S(\tilde{Q}^2_2X^3),
\end{eqnarray*}
so we have proved the relation (\ref{qf2}). \\
ii) For all $\chi \in H^*$ we have
\begin{eqnarray*}
{\cal Q}(\chi )&\pile{{\rm (\ref{qf1}, \ref{eu})}\\ =}&
<\chi, S(X^2_2\tprb )f^1R^2r^1g^1S(q^2)X^3>
X^1S(X^2_1\tpra )f^2R^1r^2g^2S(q^1)\\
{\rm (twice\hspace*{2mm}\ref{ext})}&=&
<\chi , S(q^2r^1R^2X^2_2\tprb )X^3>
X^1S(q^1r^2R^1X^2_1\tpra )\\
{\rm (\ref{qf2})}&=&S(\ov{\cal Q}(\chi \circ S)).
\end{eqnarray*}
Since the antipode $S$ is bijective we conclude that
${\cal Q}$ is bijective if and only if
$\ov{\cal Q}$ is bijective, so our proof is complete.
\end{proof}
In the Hopf case perhaps the most important example of
factorizable Hopf algebra is the Drinfeld double.
We will see that this is also true in the quasi-Hopf
case. Firstly, from \cite{hn1}, \cite{hn2}, \cite{bc},
we recall the definition of the Drinfeld double $D(H)$ of a
finite dimensional quasi-Hopf algebra $H$. We notice that the quasi-Hopf
algebra $D(H)$ was first described by Majid \cite{maj2} in the form
of an implicit Tannaka-Krein reconstruction theorem.\\
Let $\{{}_ie\}_{i=\ov{1, n}}$ be a basis of $H$ and $\{{}^ie\}_{i=\ov{1, n}}$
the corresponding dual basis of $H^*$. It is well known that $H^*$
is a coassociative coalgebra with comultiplication
$$
\Delta (\chi )=\chi _1\ot \chi _2:=\chi ({}_ie{}_je){}^ie\ot {}^je
$$
and counit $\va $. Moreover, $H^*$ is a $H$-bimodule, by
$$
<h^{'}\rh \chi \lh h^{''}, h>=<\chi , h^{''}hh^{'}>
$$
for all $h, h^{'}, h^{''}\in H$ and $\chi \in H^*$. The convolution
is a multiplication on $H$. It is not associative but ensures us that
$H^*$ is an algebra in the monoidal category of $H$-bimodules.
We also introduce $\ov{S} : H^*\ra H^*$ as the coalgebra antimorphism
dual to $S$, i. e. $<\ov{S}(\chi ), h>=<\chi , S(h)>$. Now, consider
$\Omega \in H^{\ot 5}$ given by
\begin{eqnarray}
\Omega &=&\Omega ^1\ot \Omega ^2\ot \Omega ^3\ot \Omega ^4\ot \Omega ^5\nonumber \\
&=&X^1_{(1, 1)}y^1x^1\ot X^1_{(1, 2)}y^2x^2_1\ot X^1_2y^3x^2_2\ot
\smi (f^1X^2x^3)\ot \smi (f^2X^3)\label{ome}
\end{eqnarray}
where $f=f^1\ot f^2$ is the Drinfeld twist defined in (\ref{f}). The
quantum double $D(H)$ is defined as follows. As $k$ vector spaces
$D(H)=H^*\ot H$ and the multiplication is given by
\begin{equation}\label{mdd}
(\chi \Join h)(\psi \Join h^{'})=
[(\Omega ^1 \rh \chi \lh \Omega ^5)
(\Omega ^2h_{(1, 1)}\rh \psi \lh \smi (h_2)\Omega ^4]
\Join \Omega ^3h_{(1, 2)}h^{'}
\end{equation}
for all $\chi , \psi \in H^*$ and $h, h^{'}\in H$. The unit is $\va \Join 1$.
By the above, it is easy to see that
\begin{equation}\label{mdd1}
(\va \Join h)(\chi \Join h^{'})=h_{(1, 1)}\rh \chi \lh \smi (h_2)\Join
h_{(1, 2)}h^{'}
\mbox{${\;\;}$and${\;\;}$}
(\chi \Join h)(\va \Join h^{'})=\chi \Join hh^{'}
\end{equation}
for any $h, h^{'}\in H$, $\chi \in H^*$. The explicit formulas
for the comultiplication, the counit and the antipode are
\begin{eqnarray}
\Delta _D(\chi \Join h)&=&(\va \Join X^1Y^1)(p^1_1x^1\rh \chi _2\lh
Y^2\smi (p^2)\Join p^1_2x^2h_1)\nonumber\\
&&\ot (X^2_1\rh \chi _1\lh \smi (X^3)\Join X^2_2Y^3x^3h_2),\label{cdd}\\
\va _D(\chi \Join h)&=&\va (h)\chi (\smi (\a )),\label{codd}\\
S_D(\chi \Join h)&=&(\va \Join S(h)f^1)(p^1_1U^1\rh \ov{S}^{-1}(\chi )\lh
f^2\smi (p^2)\Join p^1_2U^2),\label{add}\\
\a _D&=&\va \Join \a ,\\
\b _D&=&\va \Join \b .\label{abdd}
\end{eqnarray}
Here $p_R=p^1\ot p^2$, $f=f^1\ot f^2$ and $U=U^1\ot U^2$
are the elements defined by (\ref{qr}), (\ref{f}) and (\ref{eu}),
respectively. Thus, $D(H)$ is a quasi-Hopf algebra and $H$ is a
quasi-Hopf subalgebra via the canonical morphism $i_D:\ H\ra D(H)$,
$i_D(h)=\va \Join h$. Moreover, $D(H)$ is QT, the $R$-matrix being
defined by
\begin{equation}\label{madd}
{\cal R}=(\va \Join \smi (p^2){}_iep^1_1)\ot ({}^ie\Join p^1_2).
\end{equation}
We are now able to prove the following result.

\begin{proposition}\label{pr2.3}
Let $H$ be a finite dimensional quasi-Hopf algebra and
$D(H)$ its Drinfeld double. Then $D(H)$ is a factorizable
quasi-Hopf algebra.
\end{proposition}
\begin{proof}
We will show that in the Drinfeld double case
the map $\ov{\cal Q}$ defined by (\ref{qf2}) is bijective, so
by Proposition \ref{pr2.2} it follows that $D(H)$ is factorizable.
For this we will compute first the element
${\cal R}^2\mathbf{R}^1\ot {\cal R}^1\mathbf{R}^2$, where we denote by
$\mathbf{R}^1\ot \mathbf{R}^2$ another copy of the $R$-matrix ${\cal R}$
of $D(H)$. In fact, if we denote by $P^1\ot P^2$ another
copy of the element $p_R$ then we compute
\begin{eqnarray*}
{\cal R}^2\mathbf{R}^1\ot {\cal R}^1\mathbf{R}^2
&=&({}^ie\Join p^1_2)(\va \Join \smi (P^2){}_jeP^1_1)\ot
(\va \Join \smi (p^2){}_iep^1_1)({}^je\Join P^1_2)\\
{\rm (\ref{mdd1})}&=&
{}^ie\Join p^1_2\smi (P^2){}_jeP^1_1
\ot (\smi (p^2){}_iep^1_1)_{(1, 1)}\rh {}^je\lh
\smi ((\smi (p^2){}_iep^1_1)_2)\\
&&\Join (\smi (p^2){}_iep^1_1)_{(1, 2)}P^1_2\\
&=&{}^ie\Join p^1_2\smi ((\smi (p^2){}_iep^1_1)_2P^2){}_je
(\smi (p^2){}_iep^1_1)_{(1, 1)}P^1_1
\ot {}^je\Join (\smi (p^2){}_iep^1_1)_{(1, 2)}P^1_2\\
{\rm (\ref{qr1})}&=&
{}^ie\Join \smi ((\smi (p^2){}_ie)_2P^2){}_je
(\smi (p^2){}_ie)_{(1, 1)}P^1_1p^1_1
\ot {}^je\Join (\smi (p^2){}_ie)_{(1, 2)}P^1_2p^1_2.
\end{eqnarray*}
Now, $H$ is a quasi-Hopf subalgebra of $D(H)$, so
we have to calculate the element
$$
b^1\ot b^2:=
(\va \Join \smi (X^3)q^2){\cal R}^1\mathbf{R}^2(\va \Join X^2_2\tprb )
\ot (\va \Join q^1){\cal R}^2\mathbf{R}^1(\va \Join X^2_1\tpra \smi (X^1)).
$$
By dual basis and (\ref{mdd1}) we have
\begin{eqnarray*}
b^1\ot b^2&=&(\va \Join \smi (X^3))
(e^j\Join (q^2\smi (q^1_2p^2)e_i)_{(1, 2)}((q^1_1)_{(1, 1)}P^1)_2p^1_2)
(\va \Join X^2_2\tprb )\\
&&\ot e^i\Join \smi ((q^2\smi (q^1_2p^2)e_i)_2(q^1_1)_{(1, 2)}
P^2S((q^1_1)_2))e_j\\
&&(q^2\smi (q^1_2p^2)e_i)_{(1, 1)}((q^1_1)_{(1, 1)}P^1)_1
p^1_1)(\va \Join X^2_1\tpra \smi (X^1))\\
{\rm (\ref{qr1}, \ref{pqr})}&=&
(\va \Join \smi (X^3))(e^j\Join (e_i)_{(1, 2)}P^1_2)(\va \Join X^2_2\tprb )\\
&&\ot (e^i\Join \smi ((e_i)_2P^2)e_j(e_i)_{(1, 1)}P^1_1)(\va \Join
X^2_1\tpra \smi (X^1))\\
{\rm (\ref{mdd1}, \ref{f1tp})}&=&
(\va \Join \smi (x^3\tilde{p}^2_2)e_i)(e^j\Join P^1_2x^2\tilde{p}^2_1)
\ot (e^i\Join \smi (P^2)e_jP^1_1x^1\tpra ).
\end{eqnarray*}
Now we want an explicit formula for the element $S_D(b^1)\ot b^2$. To
this end we need the following relations:
\begin{eqnarray}
&&S(P^1_2x^2\tilde{p}^2_1)_1f^1_1p^1\ot
S(P^1_2x^2\tilde{p}^2_1)_2f^1_2p^2S(f^2)S^2(P^1_1x^1\tpra )\nonumber\\
&&\hspace*{1cm}=g^1S(P^1y^3x^2_2\tilde{p}_{(1, 2)})\ot g^2S(S(y^1x^1\tpra )
\a y^2x^2_1\tilde{p}^2_{(1, 1)}),\label{form1}\\
&&S(P^1)_2U^2\ot S(P^1)_1U^1P^2=g^2\ot g^1.\label{form2}
\end{eqnarray}
The first one follows applying (\ref{ca}, \ref{g2}, \ref{pf}, \ref{q1}, \ref{q5})
and then $f^1\b S(f^2)=S(\a )$ and (\ref{q1}, \ref{q5}). The second one can
be proved more easily by using of (\ref{eu}, \ref{ca}) and (\ref{pqr}),
we leave the details to the reader.
Therefore, if we denote by $G^1\ot G^2$ another copy of $f^{-1}$ then
from the definition (\ref{add}) of $S_D$ we obtain
\begin{eqnarray*}
S_D(b^1)\ot b^2&=&(\va \Join S(P^1_2x^2\tilde{p}^2_1))S_D(e^j\Join 1)
(\va \Join S(e_i)x^3\tilde{p}^2_2)
\ot (e^i\Join \smi (P^2)e_jP^1_1x^1\tpra )\\
&=&(\va \Join S(P^1_2x^2\tilde{p}^2_1)f^1)(p^1_1U^1\rh \ov{S}^{-1}(e^j)\lh
f^2\smi (p^2)\Join p^1_2U^2S(e_i)x^3\tilde{p}^2_2)\\
&&\ot (e^i\Join \smi (P^2)e_jP^1_1x^1\tpra )\\
&=&(\va \Join S(P^1_2x^2\tilde{p}^2_1)f^1)(\ov{S}^{-1}(e^j)
\Join p^1_2U^2S(e_i)x^3\tilde{p}^2_2)\\
&&\ot (e^i\Join \smi (p^1_1U^1P^2)e_j\smi (f^2\smi (p^2))P^1_1x^1\tpra )\\
{\rm (\ref{mdd1})}&=&
(\ov{S}^{-1}(e^j)\Join S(P^1_2x^2\tilde{p}^2_1)_{(1, 2)}f^1_{(1, 2)}
p^1_2U^2S(e_i)x^3\tilde{p}^2_2)\\
&&\ot (e^i\Join \smi (S(P^1_2x^2\tilde{p}^2_1)_{(1, 1)}f^1_{(1, 1)}p^1_1U^1P^2)
e_j\smi (f^2\smi (S(P^1_2x^2\tilde{p}^2_1)_2f^1_2p^2))P^1_1x^1\tpra )\\
{\rm (\ref{form1})}&=&
(\ov{S}^{-1}(e^j)\Join g^1_2S(P^1y^3x^2_2\tilde{p}^2_{(1, 2)})_2
U^2S(e_i)x^3\tilde{p}^2_2)\\
&&\ot (e^i\Join \smi (g^1_1S(P^1y^3x^2_2\tilde{p}^2_{(1, 2)})_1U^1P^2)
e_jS^{- 2}(g^2S(S(y^1x^1\tpra )\a y^2x^2_1\tilde{p}^2_{(1, 1)})\\
{\rm (\ref{f1tp}, \ref{form2}, \ref{ql})}&=&
(\ov{S}^{-1}(e^j)\Join g^1_2S(\tqrb X^2_{(2 , 2)}\tilde{p}^2_2)_2
G^2S(e_i)X^3)\\
&&\ot (e^i\Join \smi (g^1_1S(\tqrb X^2_{(2, 2)}
\tilde{p}^2_2)_1G^1)e_jS^{- 2}(g^2S(X^1S(X^2_1\tpra )\tqra X^2_{(2, 1)}
\tilde{p}^2_1))\\
{\rm (\ref{ql1}, \ref{pql})}&=&
(\ov{S}^{-1}(e^j)\Join g^1_2S(X^2)_2G^2S(e_i)X^3)
\ot (e^i\Join \smi (g^1_1S(X^2)_1G^1)e_jS^{- 2}(g^2S(X^1))\\
{\rm (\ref{ca})}&=&
(\ov{S}^{-1}(e^j)\Join S(e_i))\ot
(X^2_1\smi (g^1_2G^2)\rh e^i\lh \smi (X^3)\\
&&\Join X^2_2\smi (g^1_1G^1)
e_jS^{- 2}(g^2S(X^1))).
\end{eqnarray*}
We are able now to prove that $\ov{\cal Q}$ is injective. To this end,
let $\mathbf{D}\in (D(H))^*$ such that $\ov{\cal Q}(\mathbf{D}\circ S_D)=0$.
That means $\mathbf{D}(S_D(b^1))b^2=0$ and it is equivalent to
$$
\mathbf{D}(\ov{S}^{-1}({}^je)\Join S({}_ie))
<{}^ie, \smi (X^3)hX^2_1\smi (g^1_2G^2)>
<\chi , X^2_2\smi (g^1_1G^1){}_jeS^{- 2}(g^2S(X^1))>=0,
$$
for all $h\in H$ and $\chi \in H^*$. In particular,
\begin{eqnarray*}
&&\mathbf{D}(\ov{S}^{-1}({}^je)\Join S({}_ie))
<{}^ie, \smi (X^3)(\smi (x^3)h\smi (F^2f^1_2)x^2_1)X^2_1\smi (g^1_2G^2)>\\
&&<S^{- 2}(S(x^1)f^2)\rh \chi \lh \smi (F^1f^1_1)x^2_2,
X^2_2\smi (g^1_1G^1){}_jeS^{- 2}(g^2S(X^1))>=0,
\end{eqnarray*}
for all $h\in H$ and $\chi \in H^*$, and therefore
$$
\mathbf{D}(\ov{S}^{-1}(\chi )\Join S(h))=0
~~~\forall ~\chi \in H^*~~\mbox{and}~~h\in H.
$$
Since the antipode $S$ is bijective ($H$ is finite dimensional) we conclude
that $\mathbf{D}=0$ and using the bijectivity of $S_D$ it follows that
$\ov{\cal Q}$ is injective. Finally, $\ov{\cal Q}$ is bijective because of finite
dimensionality of $D(H)$, so the proof is finished.
\end{proof}

\section{Transmutation theory for dual quasi-Hopf algebras}\selabel{3}
\setcounter{equation}{0}
As we pointed out, our definition of a factorizable quasi-Hopf algebra
has a categorical interpretation. In order to this end, we first need
to associate to any co-quasi-triangular dual quasi-Hopf algebra $A$ a
braided commutative Hopf algebra $\un{A}$ in the category of right $A$-comodules,
${\cal M}^A$. For this we will use the dual reconstruction theorem
due to Majid \cite{maj1}. We notice that this reconstruction theorem
is the formal dual case of the reconstruction theorem used in \cite{bn2}
but, even in the finite dimensional case,
the resulting object $\un{A}$ cannot be viewed
as the formal dual of the object obtained in \cite{bn2}. Thus, we will present
all the details concerning how we can get the structure of
$\un{A}$ as a Hopf algebra in ${\cal M}^A$.\\
Throughout, $A$ will be a dual quasi-bialgebra or a dual
quasi-Hopf algebra. Following \cite{maj}, a dual quasi-bialgebra $A$ is a
coassociative coalgebra $A$ with comultiplication
$\Delta $ and counit $\va $ together with coalgebra morphisms
$m_A :\ A\ot A\ra A$ (the multiplication;
we write $m_A(a\ot b)=ab$)
and $\eta _A:\ k\ra A$ (the unit; we write $\eta _A(1)=1$), and an
invertible element $\v \in (A\ot A\ot A)^*$ (the reassociator),
such that for all $a, b, c, d\in A$ the following relations
hold (summation understood):
\begin{eqnarray}
&&a_1(b_1c_1)\v (a_2, b_2, c_2)=
\v (a_1, b_1, c_1)(a_2b_2)c_2,\label{dq1}\\
&&
1a=a1=a,\label{dq2}\\
&&\v (a_1, b_1, c_1d_1)\v (a_2b_2, c_2,
d_2)=\v (b_1, c_1, d_1)\v (a_1, b_2c_2, d_2)
\v (a_2, b_3, c_3),\label{dq3}\\
&&\v (a, 1, b)=\va (a)\va (b)\label{dq4}.
\end{eqnarray}
$A$ is called a dual quasi-Hopf algebra if, moreover,
there exist an anti-morphism $S$ of the coalgebra $A$ and elements
$\a , \b \in H^*$ such that, for all $a\in A$:
\begin{eqnarray}
&&S(a_1)\a (a_2)a_3=\a (a)1,
\mbox{${\;\;\;}$}
a_1\b (a_2)S(a_3)=\b (a)1,\label{dq5}\\
&&\v (a_1\b (a_2),
S(a_3), \a (a_4)a_5)=
\v ^{-1}(S(a_1), \a (a_2)a_3, \b (a_4)S(a_5))=\va
(a).\label{dq6}
\end{eqnarray}
It follows from the axioms that $S(1)=1$ and $\a (1)\b (1)=1$,
so we can assume that $\a (1)=\b (1)=1$. Moreover (\ref{dq3}) and
(\ref{dq4}) imply
\begin{equation}\label{dq7}
\v (1, a, b)=\v (a, b, 1)=\va
(a)\va (b),
\mbox{${\;\;}$$\forall a, b\in A$.}
\end{equation}
For dual quasi-Hopf algebras the antipode is an
anti-algebra morphism up to a conjugation by a twist.
Let $\gamma , \delta \in (A\ot A)^*$ be
defined by
\begin{eqnarray}
&&\gamma (a, b)=\v (S(b_2), S(a_2), a_4)\a (a_3)
\v ^{-1}(S(b_1)S(a_1), a_5, b_4)\a (b_3),\label{dug}\\
&&\delta (a, b)=\v (a_1b_1, S(b_5), S(a_4))\b (a_3)
\v ^{-1}(a_2, b_2, S(b_4))\b (b_3),\label{dud}
\end{eqnarray}
for all $a, b\in A$. If we define $f, f^{-1}\in (A\ot A)^*$,
\begin{eqnarray}
&&f(a, b)=\v ^{-1}(S(b_1)S(a_1), a_3b_3, S(a_5b_5))\b (a_4b_4)
\gamma (a_2, b_2),\label{df}\\
&&f^{-1}(a, b)=\v ^{-1}(S(a_1b_1), a_3b_3, S(b_5)S(a_5))\a (a_2b_2)
\d (a_4, b_4),\label{dg}
\end{eqnarray}
then $f$ and $f^{-1}$ are inverses in the convolution algebra and
\begin{equation}\label{dca}
f(a_1, b_1)S(a_2b_2)f^{-1}(a_3, b_3)=S(b)S(a),
\end{equation}
for all $a, b\in A$. Moreover, the following relations hold:
\begin{equation}\label{dgdf}
\gamma (a, b)=f(a_1, b_1)\a (a_2b_2)
\mbox{${\;\;}$and${\;\;}$}
\delta (a, b)=\b (a_1b_1)f^{-1}(a_2, b_2).
\end{equation}
Suppose that $A$ is a dual quasi-bialgebra or a dual
quasi-Hopf algebra. A right $A$-comodule $M$ is a
$k$-vector space together with a linear map $\r _M:\ M\ra M\ot A$
required to satisfy
$$
(\r _M\ot id_A)\circ \r _M=(id_M\ot \Delta )\circ \r _M
\mbox{${\;\;}$and${\;\;}$}
(id_M\ot \va )\circ \r _M=id_M.
$$
As usual, we denote $\r _M(m)=m_{(0)}\ot m_{(1)}$. The category of
right $A$-comodules is denoted by ${\cal M}^A$ and it
is a monoidal category. The tensor product is given
via $m_A$, i.e. for any $M, N\in {\cal M}^A$,
$M\ot N\in {\cal M}^A$ via the structure map

\begin{equation}\label{fox4}
\r _{M\ot N}(m\ot n)=m_{(0)}\ot n_{(0)}\ot m_{(1)}n_{(1)}.
\end{equation}

The associativity constraints ${\bf a}_{M, N, P} : (M\ot N)\ot P\ra
M\ot (N\ot P)$ are defined by

\begin{equation}\label{fox5}
{\bf a}_{M, N, P}(m, n, p)=\v (m_{(1)}, n_{(1)}, p_{(1)})
m_{(0)}\ot (n_{(0)}\ot p_{(0)}),
\end{equation}

for all $M, N, P\in {\cal M}^A$. The unit is $k$ as a
trivial right $A$-comodule,
and the left and right unit constraints are the usual ones.\\
If $A$ is a dual quasi-Hopf algebra then any finite dimensional
object $M$ of ${\cal M}^A$ has a left dual, i.e. the category
of finite dimensional right $A$-comodules is left rigid.
Indeed, the left dual of
$M$ is $M^*$ with the right $A$-comodule structure
$$
\rho _{M^*}(m^*)=<m^*, {}_im_{(0)}>{}^im\ot S({}_im_{(1)}),
$$
for all $m^*\in M^*$, where $({}_im)_i$ is a basis in $M$ with dual basis
$({}^im)_i$. The evaluation and coevaluation maps are defined by
\begin{eqnarray}
&&ev_M : M^*\ot M\ra k
\mbox{,${\;\;}$}
ev_M(m^*\ot m)=\a (m_{(1)})m^*(m_{(0)}),\label{dev}\\
&&coev_M : k\ra M\ot M^*
\mbox{,${\;\;}$}
coev_M(1)=\b ({}_im_{(1)}){}_im_{(0)}\ot {}^im.\label{dcoev}
\end{eqnarray}
A dual quasi-bialgebra or dual quasi-Hopf algebra is called
co-quasi-triangular (CQT for short) if there exists a $k$-bilinear
form $\sigma : A\ot A\ra k$ such that the following relations hold:
\begin{eqnarray}
&&\sigma (ab, c)=\v (c_1, a_1, b_1)\sigma (a_2, c_2)
\v ^{- 1}(a_3, c_3, b_2)\sigma (b_3, c_4)\v (a_4, b_4, c_5),\label{cqt1}\\
&&\sigma (a, bc)=\v ^{- 1}(b_1, c_1, a_1)\sigma (a_2, c_2)\v (b_2, a_3, c_3)
\sigma (a_4, b_3)\v ^{- 1}(a_5, b_4, c_4),\label{cqt2}\\
&&\sigma (a_1, b_1)a_2b_2=b_1a_1\sigma (a_2, b_2)\label{cqt3}\\
&&\sigma (a, 1)=\sigma (1, a)=\va (a),\label{cqt4}
\end{eqnarray}
for all $a, b, c\in A$.\\
As in the Hopf case, if $A$ is a CQT dual quasi-Hopf algebra then
we can prove that the bilinear form $\sigma $ is
convolution invertible, and that the antipode $S$ is bijective.

\begin{proposition}
Let $(A, \sigma )$ be a CQT dual quasi-Hopf algebra. Then:
\begin{itemize}
\item[i)] $\sigma $ is convolution invertible. More exactly, its inverse
(denoted by $\sigma ^{- 1}$) is given by
\begin{equation}\label{insig}
\sigma^{- 1}(a, b)=\v (a_1, S(a_3), b_4a_{10})\b (a_2)\v (b_1, S(a_6), a_8)
\sigma (S(a_5), b_2)\v ^{-1}(S(a_4), b_3, a_9)\a (a_7),
\end{equation}
for all $a, b\in A$.
\item[ii)] The element $u\in A^*$, given by
\begin{equation}\label{delmu}
u(a)=\v ^{- 1}(a_7, S(a_3), S^2(a_1))\sigma (a_6, S(a_4))\a (a_5)\b (S(a_2))
\end{equation}
for all $a\in A$, is invertible. Its inverse is given for all $a\in A$, by
\begin{equation}\label{indelmu}
u^{-1}(a)=\v (a_1, S^2(a_8), S(a_6))\b (a_4)\sigma (S^2(a_9), a_2)\a (S(a_7))
\v ^{-1}(S^2(a_{10}), a_3, S(a_5)).
\end{equation}
\item[iii)] For all $a\in A$,
\begin{equation}\label{dsqant}
S^2(a)=u(a_1)a_2u^{-1}(a_3).
\end{equation}
In particular, the antipode $S$ is bijective.
\end{itemize}
\end{proposition}
\begin{proof}
If $A$ is finite dimensional then the proof follows from
\cite{bn3} by duality. This is why we restrict to give a sketch of the proof,
leaving other details to the reader.\\
i) Follows by \cite[Lemma 2.2]{bn3} by duality.\\
ii) Firstly, one can prove that
\begin{equation}\label{fox1}
\sigma (S(a_1), S(b_1))\gamma (a_2, b_2)=\gamma (b_1, a_1)\sigma (a_2, b_2)
\end{equation}
for all $a, b\in A^*$ and then that
\begin{equation}\label{fox2}
f(b_1, a_1)\sigma (a_2, b_2)f^{-1}(a_3, b_3)=\sigma (S(a), S(b)).
\end{equation}
Note that these formulas are the formal
duals of \cite[Lemma 2.3]{bn3}. Secondly, using (\ref{fox2}) and the
equalities
\begin{equation}\label{fox3}
u(a_2)S^2(a_1)=u(a_1)a_2
\mbox{${\;\;}$and${\;\;}$}
\a (S(a_1))u(a_2)=\sigma (a_3, S(a_1))\a (a_2)
\end{equation}
one can show that $u\circ S^2=u$ (see \cite[Lemmas 2.4 and 2.5]{bn3} for the dual case).
Now, using (\ref{fox3}), (\ref{cqt2}) and (\ref{cqt4}) it can be proved
that $u^{-1}$ defined
by (\ref{indelmu}) is a left inverse of $u$. It is also a right inverse since
$$
u(a_1)u^{-1}(a_2)=u^{-1}(S^2(a_1))u(a_2)=u^{-1}(S^2(a_1))u(S^2(a_2))=\va (S^2(a))=\va (a),
$$
because of (\ref{fox3}) and $u\circ S^2=u$ (for the dual case see \cite[Theorem 2.6]{bn3}).
\end{proof}
By the above Proposition, if $(A, \sigma )$ is a CQT dual quasi-Hopf algebra
then it follows that ${\cal M}^A$ is a braided category. For any
$M, N\in {\cal M}^A$, the braiding is given by
$$
c_{M, N}(m\ot n)=\sigma (m_{(1)}, n_{(1)})n_{(0)}\ot m_{(0)}.
$$
We will use now the dual
braided reconstruction theorem in order to obtain the structure
of $\un{A}$ as a braided Hopf algebra in ${\cal M}^A$.
Let ${\cal C}$ and ${\cal D}$ be two monoidal
categories with ${\cal D}$ braided.
If $\mathbb{F}, \mathbb{G}: {\cal C}\ra {\cal D}$ are two functors then we denote
by $Nat(\mathbb{F}, \mathbb{G})$ the set of natural transformations
$\xi : \mathbb{F}\ra \mathbb{G}$, by $\mathbb{F}\ot M : {\cal C}\ra {\cal D}$ the functor
$(\mathbb{F}\ot M)(N)=\mathbb{F}(N)\ot M$, where $N\in {\cal C}$, $M\in {\cal D}$,
and by $Hom(M, M^{'})$ the set of morphism between $M$ and $M^{'}$ in ${\cal D}$. Suppose that
there is an object $B\in {\cal D}$ such that for all $M\in {\cal D}$
$$
Hom(B, M)\cong Nat(\mathbb{F}, \mathbb{F}\ot M)
$$
by functorial bijections $\theta _M$ and let
$\mu =\{\mu _N: \mathbb{F}(N)\ra \mathbb{F}(N)\ot B\mid N\in {\cal C}\}$
be the natural transformation corresponding to the identity
morphism $id_B$. Then, using $\mu $ and the braiding in ${\cal D}$ we have induced maps
$$
\Theta ^s_M: Hom(B^{\ot s}, M)\cong Nat(\mathbb{F}^s, \mathbb{F}^s\ot M)
$$
and we assume that these are bijections. This is the representability assumption for
comodules and is always satisfied if ${\cal D}$ is co-complete and if the image of
$\mathbb{F}$ is rigid, cf. \cite{maj1}. Then, using $(\Theta ^2_B)^{- 1}$,
$\mu _{\un{1}}$, $\theta ^{-1}_{B\ot B}$ and $\theta ^{-1}_{\un{1}}$ we can define
a multiplication, a unit, a comultiplication and a counit for $B$.

\begin{theorem}\cite{maj1}\label{thmaj}
Let ${\cal C}$ and ${\cal D}$ be monoidal categories with
${\cal D}$ braided and $\mathbb{F}: {\cal C}\ra {\cal D}$
a monoidal functor satisfying the representability assumption for
comodules. Then $B$ as above is a bialgebra in ${\cal D}$. If ${\cal C}$
is rigid then $B$ is a Hopf algebra in ${\cal D}$.
\end{theorem}

Let now $(A, \sigma )$ be a CQT dual quasi-Hopf algebra,
$A_R$ the $k$-vector space $A$
viewed as a right $A$-comodule via $\Delta $ and $\un{A}$ the same
$k$-vector space $A$ but viewed now as an object of ${\cal M}^A$ via
the right adjoint coaction:
\begin{equation}\label{draa}
\r _{\un{A}}(a)=a_2\ot S(a_1)a_3,
\end{equation}
for all $a\in A$. We apply now the Theorem \ref{thmaj} in the
case that ${\cal C}={\cal D}={\cal M}^A$ and $\mathbb{F}=id$. The first
step is to show that $\un{A}$ is the representability object which we need.\\
Dual to the quasi-Hopf case, since the antipode $S$ is bijective,
we define the elements $p_L, q_L\in (A\ot A)^*$ given by
\begin{equation}\label{dpl}
p_L(a, b)=\v (\smi (a_3), a_1, b)\b (\smi (a_2))
\mbox{${\;\;}$and${\;\;}$}
q_L(a, b)=\v ^{-1}(S(a_1), a_3, b)\a (a_2)
\end{equation}
for all $a, b\in A$. Then, for all $a, b\in A$, the following
relations hold:
\begin{equation}\label{dpl1}
p_L(a_2, b_2)\smi (a_3)(a_1b_1)=p_L(a, b_1)b_2,
\mbox{${\;\;}$}
q_L(a_2, b_1)S(a_1)(a_3b_2)=q_L(a, b_2)b_1,
\end{equation}
\begin{equation}\label{dpql}
p_L(S(a_1), a_3b_2)q_L(a_2, b_1)=\va (a)\va (b),
\mbox{${\;\;}$}
q_L(\smi (a_3), a_1b_1)p_L(a_2, b_2)=\va (a)\va (b).
\end{equation}

\begin{lemma}\label{lm3.1}
Let $A$ be a dual quasi-Hopf algebra and $M\in {\cal M}^A$.
If we define
\begin{eqnarray}
&&\theta _M :\ Hom(\un {A}, M)\ra Nat(id, id\ot M),\nonumber\\
&&\theta _M(\chi )_N(n)=
p_L(S(n_{(1)}), n_{(3)})n_{(0)}\ot \chi (n_{(2)}),\label{theta}
\end{eqnarray}
for all $\chi \in Hom(\un {A}, M)$, $N\in {\cal M}^A$
and $n\in N$, then $\theta _M$ is well
defined and a bijection. Its inverse,
$\theta _M^{-1}:\ Nat(id, id\ot M)\ra Hom(\un {A}, M)$, is given
for all $\xi \in Nat(id, id\ot M)$ by
\begin{equation}\label{intheta}
\theta ^{-1}_M(\xi )(a)=q_L(a_1, (a_2)_{<1>_{(1)}})
\va ((a_2)_{<0>})(a_2)_{<1>_{(0)}},
\end{equation}
for all $a\in A$, where we denote
$\xi _{A_R}(a)=a_{<0>}\ot a_{<1>}$.
\end{lemma}
\begin{proof}
We have to prove first that $\theta _M$ is well defined, that means
$\theta _M(\chi )_N$ is a right $A$-colinear map and
$\theta _M(\chi )$ is a natural transformation. Since $\chi :\ \un{A}\ra M$
is a morphism in ${\cal M}^A$ we have
\begin{equation}\label{fox6}
\chi (a)_{(0)}\ot \chi (a)_{(1)}=\chi (a_2)\ot S(a_1)a_3,
\end{equation}
for all $a\in A$. Now, if $n\in N$ then:
\begin{eqnarray*}
\r _{N\ot M}(\theta _M(\chi )_N(n))&=&p_L(S(n_{(1)}), n_{(3)})
\r _{N\ot M}(n_{(0)}\ot \chi (n_{(2)}))\\
{\rm (\ref{fox4})}&=&p_L(S(n_{(2)}), n_{(4)})n_{(0)}\ot
\chi (n_{(3)})_{(0)}\ot n_{(1)}\chi (n_{(3)})_{(1)}\\
{\rm (\ref{fox6})}&=&p_L(S(n_{(2)}), n_{(6)})n_{(0)}
\ot \chi (n_{(4)})\ot
n_{(1)}(S(n_{(3)})n_{(5)})\\
{\rm (\ref{dpl1})}&=&p_L(S(n_{(1)}), n_{(3)})n_{(0)}
\ot \chi (n_{(2)})\ot n_{(4)}\\
{\rm (\ref{theta})}&=&\theta _M(\chi )_N(n_{(0)})\ot n_{(1)}
=(\theta _M(\chi )_N\ot id_A)(\r _N(n)),
\end{eqnarray*}
as needed. It is not hard to see that $\theta _M(\chi )$ is a natural
transformation, so we are left to show that $\theta _M^{-1}$
is also well defined, and that $\theta _M$ and $\theta _M^{-1}$
are inverses. The first assertion follows from the following. Since
$\xi _{A_R}$ is a right $A$-comodule map we have
\begin{equation}\label{fox7}
(a_1)_{<0>}\ot (a_1)_{<1>}\ot a_2=
a_{<0>_1}\ot a_{<1>_{(0)}}\ot a_{<0>_2}a_{{<1>}_{(1)}},
\end{equation}
for all $a\in A$. On the other hand, for all $a^{*}\in A^*$ the map
$\l _{a^*}:\ A_R\ra A_R$, $\l _{a^*}(a):=a^*(a_1)a_2$, is right
$A$-colinear. Since $\xi $ is functorial under the morphism
$\l _{a^*}$ we obtain that
$$
a^*(a_{<0>_1})a_{<0>_2}\ot a_{<1>}=
a^*(a_1)(a_2)_{<0>}\ot (a_2)_{<1>},
$$
for all $a^*\in A^*$ and $a\in A$, and this is equivalent to
\begin{equation}\label{fox8}
a_{<0>_1}\ot a_{<0>_2}\ot a_{<1>}=
a_1\ot (a_2)_{<0>}\ot (a_2)_{<1>},
\end{equation}
for all $a\in A$. Then for $a\in A$ we have
\begin{eqnarray*}
(\theta _M^{-1}(\xi )\ot id_A)(\r _{\un{A}}(a))&=&
\theta ^{-1}_M(\xi )(a_2)\ot S(a_1)a_3\\
&=&q_L(a_2, (a_3)_{<1>_{(1)}})\va ((a_3)_{<0>})
(a_3)_{<1>_{(0)}}\ot S(a_1)a_4\\
{\rm (\ref{fox7})}&=&q_L(a_2, (a_3)_{<1>_{(1)}})\va ((a_3)_{<0>_1})
(a_3)_{<1>_{(0)}}\ot S(a_1)((a_3)_{<0>_2}(a_3)_{<1>_{(2)}})\\
&=&q_L(a_2, (a_3)_{<1>_{(1)}})\va ((a_3)_{<0>_2})
(a_3)_{<1>_{(0)}}\ot S(a_1)((a_3)_{<0>_1}(a_3)_{<1>_{(2)}})\\
{\rm (\ref{fox8})}&=&q_L(a_2, (a_4)_{<1>_{(1)}})
\va ((a_4)_{<0>})(a_4)_{<1>_{(0)}}\ot S(a_1)(a_3(a_4)_{<1>_{(2)}})\\
{\rm (\ref{dpl1})}&=&q_L(a_1, (a_2)_{<1>_{(2)}})\va ((a_2)_{<0>})
(a_2)_{<1>_{(0)}}\ot (a_2)_{<1>_{(1)}}\\
&=&(\r _M\circ \theta ^{-1}_M(\xi ))(a),
\end{eqnarray*}
so $\theta ^{-1}_M(\xi )$ is a right $A$-comodule map. We show now
that $\theta ^{-1}_M$ is a left inverse for $\theta _M$. Indeed,
from definitions we have
$$
\theta _M(\xi )_{A_R}(a)=p_L(S(a_2), a_4)a_1\ot \xi (a_3):=a_{<0>}\ot a_{<1>},
$$
for all $a\in A$, and therefore
\begin{eqnarray*}
(\theta ^{-1}_M\circ \theta _M)(\chi )(a)&=&
q_L(a_1, \chi (a_4)_{(1)})\va (a_2)p_L(S(a_3), a_5)\chi (a_4)_{(0)}\\
&=&q_L(a_1, S(a_3)a_4)p_L(S(a_2), a_5)\chi (a_2)\\
{\rm (\ref{dpql})}&=&\va (a_1)\va (a_3)\chi (a_2)=\chi (a),
\end{eqnarray*}
for all $\chi \in Hom(\un {A}, M)$ and $a\in A$. In order to prove that
$\theta ^{-1}_M$ is a right inverse for $\theta _M$ observe first
that for any $N\in {\cal M}^A$ and $n^*\in N^*$, the map $\l _{n^*}:\ N\ra A_R$,
$\l _{n^*}(n)=n^*(n_{(0)})n_{(1)}$, is right $A$-colinear. The
fact that $\xi $ is functorial under the morphism $\l_{n^*}$ means
$$
n^*(n_{[0]_{(0)}})n_{[0]_{(1)}}\ot n_{[1]}=
n^*(n_{(0)})n_{(1)_{<0>}}\ot n_{(1)_{<1>}}
$$
where we denote $\xi _N(n):=n_{[0]}\ot n_{[1]}$.
Since it is true for any $n^*\in N^*$ we obtain
\begin{equation}\label{fox9}
n_{[0]_{(0)}}\ot n_{[0]_{(1)}}\ot n_{[1]}=
n_{(0)}\ot n_{(1)_{<0>}}\ot n_{(1)_{<1>}}
\end{equation}
for all $n\in N$. Now, for all $n\in N$ we compute:
\begin{eqnarray*}
(\theta _M\circ \theta ^{-1}_M)(\xi )_N(n)&=&\theta _M(\theta _M^{-1}(\xi ))_N(n)
=p_L(S(n_{(1)}), n_{(3)})n_{(0)}\ot \theta ^{-1}_M(\xi )(n_{(2)})\\
&=&p_L(S(n_{(1)}), n_{(4)})q_L(n_{(2)}, (n_{(3)})_{<1>_{(1)}})\va ((n_{(3)})_{<0>})
n_{(0)}\ot (n_{(3)})_{<1>_{(0)}}\\
{\rm (\ref{fox7})}&=&p_L(S(n_{(1)}), (n_{(3)})_{<0>_2}(n_{(3)})_{<1>_{(2)}})
q_L(n_{(2)}, (n_{(3)})_{<1>_{(1)}})\\
&&\va ((n_{(3)})_{<0>_1})
n_{(0)}\ot (n_{(3)})_{<1>_{(0)}}\\
{\rm (\ref{fox9})}&=&p_L(S(n_{[0]_{(1)}}), n_{[0]_{(3)}}n_{[1]_{(2)}})
q_L(n_{[0]_{(2)}}, n_{[1]_{(1)}})n_{[0]_{(0)}}\ot n_{[1]_{(0)}}\\
{\rm (\ref{dpql})}&=&\va (n_{[0]_{(1)}})\va (n_{[1]_{(1)}})
n_{[0]_{(0)}}\ot n_{[1]_{(0)}}=n_{[0]}\ot n_{[1]}=\xi _N(n),
\end{eqnarray*}
as needed, and this finishes the proof.
\end{proof}
We are now able to begin our reconstruction. The natural transformation
$\mu \in Nat(id, id\ot \un{A})$ corresponding to the
identity morphism $id_{\un{A}}$ is given by
$$
\mu _N(n)=\theta _{\un{A}}(id_{\un{A}})_N(n)=p_L(S(n_{(1)}), n_{(3)})
n_{(0)}\ot n_{(2)}
$$
for all $N\in {\cal M}^A$ and $n\in N$. By \cite[Lemma 2.4]{maj1} the
multiplication of $\un{A}$ is characterized as being the unique morphism
$\un{m} :\ \un{A}\ot \un{A}\ra \un{A}$ in ${\cal M}^A$ such that
\begin{eqnarray*}
\mu _{M\ot N}&=&(id_{M\ot N}\ot \un{m})\circ
{\bf a}^{-1}_{M, N, \un{A}\ot \un{A}}\circ
(id_M\ot {\bf a}_{N, \un{A}, \un{A}})\circ
(id_M\ot (c_{\un{A}, N}\ot id_{\un{A}}))\\
&&\circ (id_M\ot {\bf a}^{-1}_{\un{A}, N, \un{A}})\circ
{\bf a}_{M, \un{A}, N\ot \un{A}}\circ (\mu _M\ot \mu _N),
\end{eqnarray*}
for any $M, N\in {\cal M}^A$. Using the braided categorical
structure of ${\cal M}^A$ and the definition of $\mu $ it is
not hard to see that $\un{m}$ is the unique morphism in
${\cal M}^A$ which satisfies
\begin{eqnarray*}
&&p_L(S(m_{(1)}n_{(1)}), m_{(3)}n_{(3)})(m_{(0)}\ot n_{(0)})\ot
m_{(2)}n_{(2)}=p_L(S(m_{(3)}), m_{(15)})p_L(S(n_{(5)}), n_{(13)})\\
&&\hspace*{1cm}
\v (m_{(2)}, S(m_{(4)})m_{(14)}, n_{(14)})
\v ^{- 1}(S(m_{(5)})m_{(13)}, n_{(4)}, S(n_{(6)})n_{(12)})
\sigma (S(m_{(6)})m_{(12)}, n_{(3)})\\
&&\hspace*{1cm}
\v (n_{(2)}, S(m_{(7)})m_{(11)}, S(n_{(7)})n_{(11)})
\v ^{- 1}(m_{(1)}, n_{(1)}, m_{(9)}n_{(9)})\\
&&\hspace*{1cm}
(m_{(0)}\ot n_{(0)})\ot (S(m_{(8)})m_{(10)})\un {\cdot}
(S(n_{(8)}n_{(10)})
\end{eqnarray*}
for all $M, N\in {\cal M}^A$ and $m\in M$, $n\in N$, where we denote
by $a\un {\cdot}b:=\un{m}(a\ot b)$. We can easily check that
the above equality is equivalent to
\begin{eqnarray}
&&p_L(S(a_1b_1), a_3b_3)a_2b_2=
p_L(S(a_3), a_{15})p_L(S(b_5), b_{13})\v (a_2, S(a_4)a_{14}, b_{14})\nonumber\\
&&\hspace*{1cm}
\v ^{- 1}(S(a_5)a_{13}, b_4, S(b_6)b_{12})
\sigma (S(a_6)a_{12}, b_3)
\v (b_2, S(a_7)a_{11}, S(b_7)b_{11})\nonumber\\
&&\hspace*{1cm}
\v ^{- 1}(a_1, b_1, a_9b_9)
(S(a_8)a_{10})\un{\cdot}(S(b_8)b_{10})\label{complicat}
\end{eqnarray}
for all $a, b\in A$. Now, the explicit formula for the multiplication
$\un{\cdot}$ is the following:
\begin{eqnarray}
\hspace*{-1.5cm}a\un{\cdot}b&=&\v (S(a_1), a_{10}, S(b_1)b_{12})f(b_6, a_3)
\sigma (a_8, S(b_3))\v ^{- 1}(S(a_2), S(b_5), a_6b_9)
\sigma (a_4, b_7)\nonumber\\
\hspace*{-1.5cm}&&\v ^{- 1}(a_9, S(b_2), b_{11})
\v (S(b_4), a_7, b_{10})a_5b_8,\label{simplu}
\end{eqnarray}
for all $a, b\in A$. Indeed, it is easy to see that the multiplication $\un{\cdot}$
defined by (\ref{simplu}) is a right $A$-colinear map.
A straightforward but tedious computation ensures
that $\un{\cdot}$ satisfies the relation (\ref{complicat}),
we leave all these details to the reader. It is not hard to see that the unit of
$\un{A}$ is $1$, the unit of $A$.\\
Following \cite{maj1},
the comultiplication of $\un{A}$ is obtained as
$\un{\Delta}=\theta ^{-1}_{\un{A}\ot \un{A}}(\xi )$, where $\xi $ is defined
by the following composition
$$
\xi _N:\ N\stackrel{\mu _N}{\longrightarrow}N\ot \un{A}
\stackrel{\mu _N}{\longrightarrow}(N\ot \un{A})\ot \un{A}
\stackrel{{\bf a}_{N, \un{A}, \un{A}}}{\longrightarrow}N\ot (\un{A}\ot \un{A}),
$$
for all $N\in {\cal M}^A$. Explicitly, for all $n\in N$,
\begin{equation}\label{fox10}
\xi _N(n)=\v (n_{(1)}, S(n_{(3)})n_{(5)}, S(n_{(8)})n_{(10)})
p_L(S(n_{(7)}), n_{(11)})p_L(S(n_{(2)}), n_{(6)})
n_{(0)}\ot (n_{(4)}\ot n_{(9)}).
\end{equation}
The counit $\un{\va}$ is obtained as
$\un{\va}(a)=\theta ^{-1}_k(l)(a)$, where $l$ is the left unit constraint.

\begin{proposition}\label{pr3.2}
Let $A$ be a dual quasi-Hopf algebra. Then the comultiplication
of $\un{A}$ is given for all $a\in A$ by
\begin{equation}\label{fox11}
\un{\Delta}(a)=\v ^{-1}(S(a_1), a_5, S(a_7))\b (a_6)
\v (S(a_2)a_4, S(a_8), a_{10})a_3\ot a_9 .
\end{equation}
The counit of $\un{\Delta}$ is $\un{\va}=\a $.
\end{proposition}
\begin{proof}
Let us start by nothing that (\ref{dq3}) and the definitions (\ref{dpl}) of
$p_L$ and $q_L$ imply:
\begin{eqnarray}
q_L(a_1, b_1c_1)\v (a_2, b_2, c_2)&=&\a (a_3)
\v ^{-1}(S(a_2), a_4, b_1)\v ^{-1}(S(a_1), a_5b_2, c),\label{fox12}\\
\v ^{-1}(a, b_1, S(b_3)c_1)p_L(S(b_2), c_2)&=&
\v (a_1b_1, S(b_5), c)\v ^{-1}(a_2, b_2, S(b_4))\b (b_3),\label{fox13}
\end{eqnarray}
for all $a, b, c\in A$. On the other hand, from
(\ref{fox10}) we can easily see that
\begin{equation}\label{fox14}
\xi _{A_R}(a)=a_{<0>}\ot a_{<1>}=p_L(S(a_3), a_7)p_L(S(a_8), a_{12})
\v (a_2, S(a_4)a_6, S(a_9)a_{11})a_1\ot (a_5\ot a_{10}).
\end{equation}
Now, for all $a\in A$ we compute:
\begin{eqnarray*}
\un{\Delta }_{\un{A}}(a)&=&\theta ^{-1}_{\un{A}\ot \un{A}}(\xi )
=q_L(a_1, (a_2)_{<1>_{(1)}})\va ((a_2)_{<0>})(a_2)_{<1>_{(0)}}\\
{\rm (\ref{fox14})}&=&q_L(a_1, (a_5\ot a_{10})_{(1)})
p_L(S(a_8), a_{12})p_L(S(a_3), a_7)\\
&&\v (a_2, S(a_4)a_6, S(a_9)a_{11})(a_5\ot a_{10})_{(0)}\\
{\rm (\ref{fox4})}&=&q_L(a_1, (S(a_5)a_7)(S(a_{12})a_{14}))p_L(S(a_{10}), a_{16})
p_L(S(a_3), a_9)\\
&&\v (a_2, S(a_4)a_8, S(a_{11})a_{15})a_6\ot a_{13}\\
{\rm (\ref{fox12})}&=&\a (a_3)\v ^{-1}(S(a_2), a_4, S(a_8)a_{10})
\v ^{-1}(S(a_1), a_5(S(a_7)a_{11}), S(a_{14})a_{16})\\
&&p_L(S(a_{13}), a_{17})p_L(S(a_6), a_{12})a_9\ot a_{15}\\
{\rm (\ref{dpl1})}&=&\a (a_3)\v ^{-1}(S(a_2), a_4, S(a_6)a_{8})
\v ^{-1}(S(a_1), a_{10}, S(a_{12})a_{14})\\
&&p_L(S(a_{11}), a_{15})p_L(S(a_5), a_{9})a_7\ot a_{13}\\
{\rm (\ref{fox13})}&=&\a (a_4)\v ^{-1}(S(a_3), a_5, S(a_7)a_{9})
\v (S(a_2)a_{11}, S(a_{15}), a_{17})
\v ^{-1}(S(a_{1}), a_{12}, S(a_{14}))\\
&&\b (a_{13})p_L(S(a_{6}), a_{10})a_8\ot a_{16}\\
{\rm (\ref{dpl})}&=&q_L(a_3, S(a_5)a_{7})
\v (S(a_2)a_{9}, S(a_{13}), a_{15})
\v ^{-1}(S(a_{1}), a_{10}, S(a_{12}))\\
&&\b (a_{11})p_L(S(a_{4}), a_{8})a_6\ot a_{14}\\
{\rm (\ref{dpql})}&=&\v ^{-1}(S(a_1), a_5, S(a_7))\b (a_6)
\v (S(a_2)a_4, S(a_8), a_{10})a_3\ot a_9,
\end{eqnarray*}
for all $a\in A$. The counit of $\un{\Delta}$ is
$\un{\va}(a)=\theta ^{-1}_k(l)(a)=q_L(a, 1)=\a (a)$
for all $a\in A$, so $\un{\va}=\a $.
\end{proof}
Let now $M$ be a finite dimensional right $A$-comodule and
$M^*$ its left dual. According to \cite[Proposition 2.9]{maj1},
the reconstructed antipode $\un{S}$ of $\un{A}$ is
characterized as being the unique morphism in ${\cal C}$
satisfying
\begin{eqnarray*}
&&(id_M\circ \un{S})\circ \mu _M=
l^{-1}_{M\ot \un{A}}\circ (id_{M\ot \un{A}}\ot ev_M)\circ
(a_{M, \un{A}, M^*}\ot id_M)\circ (a^{-1}_{M, \un{A}, M^*}\ot id_M)\\
&&\hspace*{1cm}\circ ((id_M\ot c^{-1}_{\un{A}, M^*})\ot id_M)\circ
((id_M\ot \mu _{M^*})\ot id_M)\circ (coev_M\ot id_M)\circ r_M,
\end{eqnarray*}
for any finite dimensional object $M$ of ${\cal M}^A$, where
$l$, $r$, $a$, $c$, $ev$ and $coev$ are the left unit constraints,
the right unit constraints, the associativity
constraints, the braiding of ${\cal M}^A$, and the evaluation and
coevaluation map, respectively. This comes out explicitly as
\begin{eqnarray*}
&&p_L(S(m_{(1)}), m_{(3)})m_{(0)}\ot \un{S}(m_{(2)})=
\b (m_{(3)})p_L(S^2(m_{(12)}), S(m_{(4)}))\\
&&\hspace*{1cm}\sigma ^{-1}(S^2(m_{(11)})S(m_{(5)}), S(m_{(13)}))
\v ^{-1}(m_{(2)}, S^2(m_{(10)})S(m_{(6)}), S(m_{(14)}))\\
&&\hspace*{1cm}\v (m_{(1)}[S^2(m_{(9)})S(m_{(7)})], S(m_{(15)}), m_{(17)})
\a (m_{(16)})m_{(0)}\ot S(m_{(8)}),
\end{eqnarray*}
for all finite dimensional right $A$-comodule $M$ and $m\in M$.
It follows that the above relation is equivalent to
\begin{eqnarray*}
p_L(S(a_1), a_3)\un{S}(a_2)&=&\b (a_3)p_L(S^2(a_{12}), S(a_4))
\sigma ^{-1}(S^2(a_{11})S(a_5), S(a_{13}))\\
&&\v ^{-1}(a_2, S^2(a_{10})S(a_6), S(a_{14}))
\v (a_1[S^2(a_9)S(a_7)], S(a_{15}), a_{17})\a (a_{16})S(a_8)\\
{\rm (\ref{dq3}, \ref{dq5})}&=&\b (a_2)p_L(S^2(a_{11}), S(a_3))
\sigma ^{-1}(S^2(a_{10})S(a_4), S(a_{12}))\\
&&\v (S^2(a_8)S(a_6), S(a_{14}), a_{16})
\v (a_1, [S^2(a_9)S(a_5)]S(a_{13}), a_{17})\a (a_{15})S(a_7)\\
{\rm (\ref{dpl1}, \ref{dpl})}&=&p_L(S(a_1), a_{13})
p_L(S^2(a_8), S(a_2))
\sigma ^{-1}(S^2(a_7)S(a_3), S(a_9))\\
&&\v (S^2(a_6)S(a_4), S(a_{10}), a_{12})
\a (a_{11})S(a_5)
\end{eqnarray*}
for all $a\in A$, and therefore
\begin{equation}
\un{S}(a)=p_L(S^2(a_7), S(a_1))
\sigma ^{-1}(S^2(a_6)S(a_2), S(a_8))
\v (S^2(a_5)S(a_3), S(a_9), a_{11})\a (a_{10})S(a_4)
\end{equation}
for all $a\in A$ (it is not hard to see that $\un{S}$ defined above is
right $A$-colinear). We summarize all this in the following.

\begin{theorem}\label{te3}
Let $(A, \sigma )$ be a CQT dual quasi-Hopf algebra. Then there is a
braided Hopf algebra $\un{A}$ in the category ${\cal M}^A$.
$\un{A}$ coincides with $A$ as $k$-linear space, and it is an
object in ${\cal M}^A$ by the right coadjoint action
$$
\rho _{\un{A}}(a)=a_2\ot S(a_1)a_3.
$$
The algebra structure, the coalgebra structure and the antipode
are transmuted to
\begin{eqnarray*}
a\un{\cdot}b&=&\v (S(a_1), a_{10}, S(b_1)b_{12})f(b_6, a_3)
\sigma (a_8, S(b_3))\v ^{- 1}(S(a_2), S(b_5), a_6b_9)
\sigma (a_4, b_7)\\
&&\v ^{- 1}(a_9, S(b_2), b_{11})\v (S(b_4), a_7, b_{10})a_5b_8,\\
\un{\Delta}(a)&=&\v ^{-1}(S(a_1), a_5, S(a_7))\b (a_6)
\v (S(a_2)a_4, S(a_8), a_{10})a_3\ot a_9,\\
\un{S}(a)&=&p_L(S^2(a_7), S(a_1))
\sigma ^{-1}(S^2(a_6)S(a_2), S(a_8))
\v (S^2(a_5)S(a_3), S(a_9), a_{11})\a (a_{10})S(a_4),
\end{eqnarray*}
for all $a, b\in \un{A}$. The unit element is $1$ of $A$, and
the counit is $\un{\va}=\a $. As in the Hopf case,
we will cal $\un{A}$ the associated function algebra
braided group of $A$.
\end{theorem}

\begin{remark}
The braided group $\un{A}$ is braided commutative
in the sense of \cite{maj1}. More precisely, $\un{A}$ has
a second multiplication (denoted by $\un{m}^{op}$) also making
$\un{A}$ into a braided bialgebra, and there exists a convolution
invertible morphism
$\mathfrak{R}: \un{A}\ot \un{A}\ra k$ relating $\un{m}^{op}$
by conjugation to $\un{m}$. Moreover, $\mathfrak{R}$ makes
$\un{A}$ with its to products into a CQT Hopf algebra in
${\cal M}^A$ in some sense, analogues to the definition of a
ordinary CQT Hopf algebra, see \cite{maj3}, \cite{maj}.
Now, $\un{A}$ is braided commutative mean that
$\un{m}^{op}=\un{m}$ and $\mathfrak{R}=\un{\va}\ot \un{\va}$.\\
We would like to stress that the opposite multiplication $\un{m}^{op}$ is
characterized by
\[\left.
\begin{array}{rcccccccl}
M\ot N&\stackrel{\mu _M\ot \mu _N}{\longrightarrow}&
(M\ot \un{A})\ot (N\ot \un{A})& &
M\ot N&\hspace*{-2mm}\stackrel{\mu _M\ot id_N}{\longrightarrow}&
(M\ot \un{A})\ot N\\
      &\stackrel{{\bf a}_{M, \un{A}, N\ot \un{A}}}{\longrightarrow}&
M\ot (\un{A}\ot (N\ot \un{A}))& &
      &\stackrel{{\bf a}_{M, \un{A}, N}}{\longrightarrow}&
M\ot (\un{A}\ot N)\\
      &\stackrel{id_M\ot {\bf a}^{-1}_{\un{A}, N, \un{A}}}{\longrightarrow}&
M\ot ((\un{A}\ot N)\ot \un{A})& &
      &\stackrel{id_M\ot c^{-1}_{N, \un{A}}}{\longrightarrow}&
M\ot (N\ot \un{A})\\
      &\stackrel{id_M\ot (c_{\un{A}, N}\ot id_{\un{A}})}{\longrightarrow}&
M\ot ((N\ot \un{A})\ot \un{A})& = &
      &\stackrel{id_M\ot (\mu _N\ot id_{\un{A}})}{\longrightarrow}&
M\ot((N\ot \un{A})\ot \un{A})\\
      &\stackrel{id_M\ot {\bf a}_{N, \un{A}, \un{A}}}{\longrightarrow}&
M\ot (N\ot (\un{A}\ot \un{A}))& &
      &\stackrel{id_M\ot {\bf a}_{N, \un{A}, \un{A}}}{\longrightarrow}&
M\ot (N\ot (\un{A}\ot \un{A}))\\
      &\stackrel{{\bf a}^{-1}_{M, N, \un{A}\ot \un{A}}}{\longrightarrow}&
(M\ot N)\ot (\un{A}\ot \un{A})& &
      &\stackrel{{\bf a}^{-1}_{M, N, \un{A}\ot \un{A}}}{\longrightarrow}&
(M\ot N)\ot (\un{A}\ot \un{A})\\
      &\stackrel{id_{M\ot N}\ot \un{m}^{op}}{\longrightarrow}&
(M\ot N)\ot \un{A}& &
      &\stackrel{id_{M\ot N}\ot \un{m}}{\longrightarrow}&
(M\ot N)\ot \un{A}.
\end{array}\right. \]
The above equality and $\un{m}^{op}=\un{m}$ reduce to an intrinsic form of
braided commutativity of $\un{A}$, we leave the details to the reader.
\end{remark}
\section{The categorical interpretation}
\setcounter{equation}{0}
In this Section our goal is to give a categorical interpretation for our
definition of a factorizable quasi-Hopf algebra. In the Hopf case it was given
by Majid \cite{maj}. If $(H, R)$ is a Hopf algebra then we can associate to $H$
a braided cocommutative Hopf algebra $\un{H}$ in the
braided category ${}_H{\cal M}$. As we have already
seen in the previous Section,
we can associate to any CQT (dual quasi-)Hopf
algebra $(A, \sigma )$ a braided commutative
Hopf algebra $\un{A}$ in the category of right $A$-comodules
${\cal M}^A$. Now, let $(H, R)$ be a finite dimensional
factorizable Hopf algebra and $(A, \sigma )$
the CQT Hopf algebra dual to $H$. If $\un{A}$ is viewed as a
braided Hopf algebra in $_{H}{\cal M}$ then $\un{H}$ and
$\un{A}$ are isomorphic as braided Hopf algebras.
Moreover, the isomorphism is given by the
canonical map ${\cal Q}$ considered in
Section 2. Also, $\un{A}$ is always isomorphic
to the categorical left dual of $\un{H}$.\\
We will generalize the above results to the quasi-Hopf case.
In \cite{bpv} it was introduced another multiplication
on $H$, denoted by $\bullet $, given by the formula
\begin{eqnarray}
h\bullet h^{'}&=&X^1hS(x^1X^2)\a
x^2X^3_1h^{'}S(x^3X^3_2)\label{ma}\\
{\rm (\ref{q3}, \ref{q5})}&=&X^1x^1_1hS(X^2x^1_2)
\a X^3x^2h^{'}S(x^3)\label{alma}
\end{eqnarray}
for all $h, h'\in H$ and it was proved that,
if we denote by $H_0$ this structure, then $H_0$
becomes an algebra within the monoidal category
of left $H$-modules, with unit $\beta $ and left $H$-action
given by
\begin{equation}\label{s1}
h\tr h^{'}=h_1h^{'}S(h_2),
\end{equation}
for all $h, h'\in H$. If $(H, R)$ is quasi-triangular then
$H_0$ is a Hopf algebra with bijective antipode
in ${}_H{\cal M}$, with the additional structures (see \cite{bn2}):
\begin{eqnarray}
&&\hspace*{-2cm}
\un{\Delta }(h)=h_{\un{1}}\ot h_{\un{2}}:=x^1X^1h_1g^1S(x^2R^2y^3X^3_2)
\ot x^3R^1\tr
y^1X^2h_2g^2S(y^2X^3_1),\label{und}\\
&&\hspace*{-2cm}\un{\va }(h)=\va (h),\label{unva}\\
&&\hspace*{-2cm}\un{S}(h)=X^1R^2p^2
S(q^1(X^2R^1p^1\tr h)S(q^2)X^3),\label{unant}
\end{eqnarray}
for all $h\in H$, where $R=R^1\ot R^2$ is the $R$-matrix $R$ of $H$, and
$f^{-1}=g^1\ot g^2$, $p_R=p^1\ot p^2$ and $q_R=q^1\ot q^2$
are the elements defined by (\ref{g}) and (\ref{qr}), respectively.
Thus, in the quasi-Hopf case, $\un{H}=H_0$ as an algebra with
the additional structures (\ref{und}, \ref{unva}) and (\ref{unant}).
As in the Hopf case, we will call $\un{H}$ the associated enveloping algebra
braided group of $H$. Note that, all the above structures
were obtained by using the braided
reconstruction theorem also due to Majid \cite{maj} (see \cite{bn2}
for full details).\\
Suppose now that $(H, R)$ is a finite dimensional QT quasi-Hopf algebra.
Then $H^*$, the linear dual of $H$, it is in an obvious way a CQT dual quasi-Hopf
algebra, so it makes sense to consider $\un{H^*}$, the function algebra braided group
associated to $H^*$. It is a braided Hopf algebra in the category
of right $H^*$-comodules, hence it is a braided Hopf algebra in the category
of left $H$-modules. By Theorem \ref{te3}, $\un{H^*}$ is a braided Hopf algebra
in ${}_H{\cal M}$. From (\ref{draa}), $\un{H^*}$ is a left $H$-module via
\begin{equation}\label{acovhs}
h\blacktriangleright \chi =h_2\rh \chi \lh S(h_1)
\end{equation}
for all $h\in H$ and $\chi \in H^*$. By Theorem \ref{te3},
the structure of $\un{H^*}$ as a Hopf algebra in
${}_H{\cal M}$ is given by:
\begin{eqnarray}
&&\hspace*{-1.5cm}
\chi \un{\cdot }\psi =[x^3_1Y^2r^1y^1X^2\rh \chi
\lh S(x^1X^1)f^2R^1]
[x^3_2Y^3y^3X^3_2\rh \psi
\lh S(x^2Y^1r^2y^2X^3_1)f^1R^2],\label{movhs}\\
&&\hspace*{-1.5cm}
1_{\un{H^*}}=\va ,\\
&&\hspace*{-1.5cm}
\un{\Delta }_{\un{H^*}}(\chi )=\chi _1\lh S(x^1)\ot
x^3_2X^3\rh \chi _2\lh
x^2X^1\b S(x^3_1X^2),\label{covhs}\\
&&\hspace*{-1.5cm}
\un{\va }_{\un{H^*}}(\chi )=\chi (\a ),\\
&&\hspace*{-1.5cm}
\un{S}(\chi )=q^1_2\ov{R}^1_2\tprb \rh
\chi S\lh q^2\ov{R}^2
S(q^1_1\ov{R}^1_1\tpra ),
\end{eqnarray}
for all $\chi , \psi \in H^*$.
Here $p_R=p^1\ot p^2$ and
$q_R=q^1\ot q^2$ are the elements
defined by (\ref{qr}),
$f=f^1\ot f^2$ is the Drinfeld's twist
defined by (\ref{f}),
$R^{-1}=\ov{R}^1\ot \ov{R}^2$, and
$q_L=\tqra \ot \tqrb $ is the element
given by (\ref{ql}), respectively.\\
We would like to stress that formula
(\ref{qf1}) was chosen in such a way that it
provides a left $H$-module morphism from
$\un{H^*}$ to $\un{H}$. Indeed, for all
$\chi \in H^*$ and $h\in H$ we have:
\begin{eqnarray*}
h\tr {\cal Q}(\chi )&\pile{{\rm (\ref{qf1})}\\ =}&
<\chi , S(X^2_2\tprb )f^1R^2r^1U^1X^3>
h_1X^1S(X^2_1\tpra )f^2R^1r^2U^2S(h_2)\\
{\rm (\ref{eu1}, \ref{qt3}, \ref{ca})}&=&
<\chi , S((h_{(2, 1)}X^2)_2\tprb )f^1R^2r^1U^1h_{(2, 2)}X^3>
h_1X^1S((h_{(2, 1)}X^2)_1\tpra )f^2R^1r^2U^2\\
{\rm (\ref{q1}, \ref{ql1})}&=&
<\chi , S(X^2_2\tprb h_1)f^1R^2r^1U^1X^3h_2>
X^1S(X^2_1\tpra )f^2R^1r^2U^2\\
{\rm (\ref{qf1}, \ref{acovhs})}&=&
{\cal Q}(h_2\rh \chi \lh S(h_1))
={\cal Q}(h\blacktriangleright \chi ).
\end{eqnarray*}
It is quite remarkable that (\ref{qf1}) is a braided Hopf algebra
morphism, too.

\begin{proposition}\label{pr4.1}
Let $(H, R)$ be a finite dimensional QT quasi-Hopf algebra,
$\un{H}$ the associated enveloping algebra
braided group of $H$ and $\un{H^*}$ the function algebra braided group
associated to $H^*$. Then the map ${\cal Q}$ defined by (\ref{qf1})
is a braided Hopf algebra morphism from $\un{H^*}$ to $\un{H}$.
\end{proposition}
\begin{proof}
We have already seen that ${\cal Q}$ is a morphism in ${}_H{\cal M}$.
Hence, it remains to show that ${\cal Q}$ is an
algebra and a coalgebra morphism.
To this end, we will use the second formula (\ref{qf2}) for the map
${\cal Q}$. From (\ref{q3}, \ref{q5}) it follow that
\begin{eqnarray}
&&\tqra X^1\ot \tilde{q}^2_1X^2\ot \tilde{q}^2_2X^3=
S(x^1)\tqra x^2_1\ot \tqrb x^2_2\ot x^3,\label{fff1}\\
&&x^1\ot x^2_1p^1\ot x^2_2p^2S(x^3)=X^1p^1_1\ot X^2p^1_2
\ot X^3p^2.\label{fff2}
\end{eqnarray}
We set $R=R^1\ot R^2=r^1\ot r^2=\mathbf{R}^1\ot \mathbf{R}^2
=\mathfrak{R}^1\ot \mathfrak{R}^2=\mathfrak{r}^1\ot \mathfrak{r}^2
={\cal R}^1\ot {\cal R}^2$, $q_L=\tqra \ot \tqrb
=\tQra \ot \tQrb $ and $p_R=p^1\ot p^2=P^1\ot P^2$.
Now, for all $\chi , \psi \in H^*$ we compute:
\begin{eqnarray*}
{\cal Q}(\chi \un{\cdot }\psi )&=&
<\chi , S(x^1X^1)f^2R^1\tilde{q}^1_1Z^1_1\mathfrak{R}^2_1
\mathfrak{r}^1_1p^1_1x^3_1Y^2r^1y^1X^2>\\
&&<\psi , S(x^2Y^1r^2y^2X^3_1)
f^1R^2\tilde{q}^1_2Z^1_2
\mathfrak{R}^2_2\mathfrak{r}^1_2p^1_2x^3_2Y^3y^3X^3_2>
\tilde{q}^2_1Z^2\mathfrak{R}^1\mathfrak{r}^2p^2
S(\tilde{q}^2_2Z^3)\\
{\rm (\ref{qt3}, \ref{fff2})}&=&
<\chi , S(x^1X^1)f^2\tilde{q}^1_2
[Z^1\mathfrak{R}^2\mathfrak{r}^1x^3_{(1, 1)}p^1]_2R^1Y^2r^1y^1X^2>\\
&&<\psi , S(x^2Y^1r^2y^2X^3_1)f^1\tilde{q}^1_1
[Z^1\mathfrak{R}^2\mathfrak{r}^1x^3_{(1, 1)}p^1]_1
R^2Y^3y^3X^3_2>\\
&&\tilde{q}^2_1Z^2\mathfrak{R}^1\mathfrak{r}^2x^3_{(1, 2)}
p^2S(\tilde{q}^2_2Z^3x^3_2)\\
{\rm (\ref{qt3}, \ref{q1}, \ref{ql2})}&=&
<\chi , S(X^1)\tqra \widetilde{Q}^2_1T^2Z^1_2\mathfrak{R}^2_2
\mathfrak{r}^1_2p^1_2R^1Y^2r^1y^1X^2>\\
&&<\psi , S(Y^1r^2y^2X^3_1)\tQra T^1Z^1_1
\mathfrak{R}^2_1\mathfrak{r}^1_1p^1_1R^2Y^3y^3X^3_2>\\
&&\tilde{q}^2_1\widetilde{Q}^2_{(2, 1)}T^3_1Z^2\mathfrak{R}^1
\mathfrak{r}^2p^2S(\tilde{q}^2_2\widetilde{Q}_{(2, 2)}T^3_2Z^3)\\
{\rm (\ref{q3}, \ref{q1}, \ref{fff1}, \ref{qt3})}&=&
<\chi , S(X^1)\tqra V^1\widetilde{Q}^2_1x^2_{(2, 1)}Z^2
\mathfrak{R}^2_2R^1\mathfrak{r}^1_1p^1_1Y^2r^1y^1X^2>\\
&&<\psi , S(x^1Y^1r^2y^2X^3_1)\tQra x^2_1Z^1
\mathfrak{R}^2_1R^2\mathfrak{r}^1_2p^1_2Y^3y^3X^3_2>\\
&&\tilde{q}^2_1V^2\widetilde{Q}^2_2x^2_{(2, 2)}Z^3\mathfrak{R}^1
\mathfrak{r}^2p^2S(\tilde{q}^2_2V^3x^3)\\
{\rm (\ref{q1}, \ref{qt3}, \ref{fff2})}&=&
<\chi , S(X^1)\tqra V^1\widetilde{Q}^2_1Z^2\mathfrak{R}^2_2
R^1\mathfrak{r}^1_1T^2_1Y^2r^1y^1p^1_1X^2>\\
&&<\psi , S(T^1Y^1r^2y^2(p^1_2X^3)_1)\tQra Z^1
\mathfrak{R}^2_1R^2\mathfrak{r}^1_2T^2_2Y^3y^3(p^1_2X^3)_2>\\
&&\tilde{q}^2_1V^2\widetilde{Q}^2_2Z^3\mathfrak{R}^1\mathfrak{r}^2
T^3p^2S(\tilde{q}^2_2V^3)\\
{\rm (\ref{q3}, \ref{qt3}, \ref{qt1}, \ref{qt2})}&=&
<\chi , S(X^1)\tqra V^1\widetilde{Q}^2_1\mathfrak{R}^2Z^3
x^3R^1W^2\mathfrak{r}^1z^1T^2r^1Y^1_1y^1p^1_1X^2>\\
&&<\psi , S(T^1r^2Y^1_2y^2(p^1_2X^3)_1)
\tQra Z^1\mathbf {R}^2x^2R^2W^3z^3{\cal R}^1T^3_1Y^2y^3
(p^1_2X^3)_2>\\
&&\tilde{q}^2_1V^2\widetilde{Q}^2_2\mathfrak{R}^1Z^2
\mathbf{R}^1x^1W^1\mathfrak{r}^2z^2{\cal R}^2T^3_2Y^3p^2
S(\tilde{q}^2_2V^3)\\
{\rm (\ref{q3}, \ref{qt3}, \ref{qt2})}&=&
<\chi , S(X^1)\tqra V^1\widetilde{Q}^2_1\mathfrak{R}^2Z^3x^3W^3_2
R^1T^2r^1D^1z^1_1Y^1_1y^1p^1_1X^2>\\
&&<\psi , S(W^1T^1_1r^2_1D^2z^1_2Y^1_2y^2(p^1_2X^3)_1)\tQra
Z^1\mathbf{R}^2x^2W^3_1R^2T^3z^3{\cal R}^1Y^2y^3(p^1_2X^3)_2>\\
&&\tilde{q}^2_1V^2\widetilde{Q}^2_2\mathfrak{R}^1Z^2
\mathbf{R}^1x^1W^2T^1_2r^2_2D^3z^2{\cal R}^2Y^3p^2S(\tilde{q}^2_2V^3)\\
{\rm (\ref{qt2}, \ref{q3}, \ref{q1})}&=&
<\chi , S(X^1)\tqra V^1\widetilde{Q}^2_1\mathfrak{R}^2Z^3T^3r^1
Y^1p^1_1X^2>\\
&&<\psi , S(T^1r^2_1C^1Y^2_1(p^1_2X^3)_1)\tQra
Z^1\mathbf{R}^2T^2_2r^2_{(2, 2)}{\cal R}^1C^2Y^2_2(p^1_2X^3)_2>\\
&&\tilde{q}^2_1V^2\widetilde{Q}^2_2\mathfrak{R}^1Z^2\mathbf{R}^1
T^2_1r^2_{(2, 1)}{\cal R}^2C^3Y^3p^2S(\tilde{q}^2_2V^3)\\
{\rm (\ref{fff1}, \ref{qt3}, \ref{q1}, \ref{q5})}&=&
<\chi , S(y^1X^1)\tqra \mathfrak{R}^2r^1y^2_1Y^1p^1_1X^2>\\
&&<\psi , S(x^1C^1(Y^2p^1_2)_1X^3_1)\a x^2
\mathbf{R}^2{\cal R}^1C^2(Y^2p^1_2)_2X^3_2>\\
&&\tqrb \mathfrak{R}^1r^2y^2_2x^3\mathbf{R}^1{\cal R}^2
C^3Y^3p^2S(y^3)\\
{\rm (\ref{fff2}, \ref{qt3}, \ref{q1}, \ref{q5})}&=&
<\chi , S(y^1X^1)\tqra \mathfrak{R}^2r^1y^2_1z^1X^2>\\
&&<\psi , S(C^1p^1_1X^3_1)\tQra
\mathbf{R}^2{\cal R}^1C^2p^1_2
X^3_2>\tqrb \mathfrak{R}^1r^2y^2_2z^2\tQrb \mathbf{R}^1
{\cal R}^2C^3p^2S(y^3z^3).
\end{eqnarray*}
On the other hand, if we denote by
$P^1\ot P^2$ another copy of $p_R$ then by (\ref{alma}, \ref{qr},
\ref{qf2}) we have:
\begin{eqnarray*}
{\cal Q}(\chi )\bullet {\cal Q}(\psi )&=&
<\chi , \tqra Y^1R^2r^1P^1>
<\psi , \tQra Z^1\mathfrak{R}^2\mathfrak{r}^1p^1>\\
&& q^1y^1_1\tilde{q}^2_1Y^2R^1r^2P^2
S(q^2y^1_2\tilde{q}^2_2Y^3)y^2
\widetilde{Q}^2_1Z^2\mathfrak{R}^1\mathfrak{r}^2
p^2S(y^3\widetilde{Q}^2_2Z^3)\\
{\rm (\ref{fff1}, \ref{qt3}, \ref{fff2})}&=&
<\chi , S(X^1P^1_1)\tilde{q}^1R^2r^1X^2P^1_2>
<\psi , \tQra Z^1\mathfrak{R}^2\mathfrak{r}^1p^1>\\
&&q^1y^1_1\tilde{q}^2R^1r^2
X^3P^2S(q^2y^1_2)y^2\tilde{Q}^2_1Z^2
\mathfrak{R}^1\mathfrak{r}^2p^2S(y^3\widetilde{Q}^2_2Z^3)\\
{\rm (\ref{ql1}, \ref{qt3}, \ref{q1}, \ref{qr1})}&=&
<\chi , S(X^1(q^1_1P^1)_1y^1_1)\tqra R^2r^1
X^2(q^1_1P^1)_2y^1_2>\\
&&<\psi , \tQra Z^1\mathfrak{R}^2\mathfrak{r}^1p^1>
\tqrb R^1r^2X^3q^1_2P^2S(q^2)
y^2\widetilde{Q}^2_1Z^2\mathfrak{R}^1\mathfrak{r}^2p^2S(y^3
\widetilde{Q}^2_2Z^3)\\
{\rm (\ref{pqr}, \ref{q3}, \ref{fff1}, \ref{qt3}, \ref{fff2})}&=&
<\chi , S(y^1X^1)\tqra R^2r^1y^2_1x^1X^2>
<\psi , S(Y^1p^1_1)\tQra \mathfrak{R}^2\mathfrak{r}^1Y^2p^1_2>\\
&&\tqrb R^1r^2y^2_2x^2X^3_1\tQrb \mathfrak{R}^1\mathfrak{r}^2
Y^3p^2S(y^3x^3X^3_2)\\
{\rm (\ref{ql1}, \ref{qt3}, \ref{q1}, \ref{qr1})}&=&
<\chi , S(y^1X^1)\tqra R^2r^1y^2_1x^1X^2>\\
&&<\psi , S(Y^1p^1_1X^3_1)\tQra \mathfrak{R}^2\mathfrak{r}^1
Y^2p^1_2X^3_2>
\tqrb R^1r^2y^2_2x^2\tQrb \mathfrak{R}^1\mathfrak{r}^2Y^3p^2S(y^3x^3).
\end{eqnarray*}
By the above it follows that ${\cal Q}$ is
multiplicative. Since
${\cal Q}(1_{\un{H^*}})={\cal Q}(\va )=\b =1_{\un{H}}$, we
conclude that ${\cal Q}$ is an algebra map. Thus, one has only to show
that ${\cal Q}$ is a coalgebra map. To this end, observe first that
(\ref{q3}, \ref{q5}) imply
\begin{equation}\label{fff3}
X^1_1p^1\ot X^1_2p^2S(X^2)\ot X^3=
x^1\ot x^2S(x^3_1\tpra )\ot x^3_2\tprb .
\end{equation}
Also, it is not hard to see that (\ref{und}, \ref{s1}, \ref{qt1},
\ref{qt3}) and (\ref{ext}) imply
\begin{eqnarray}
\un{\Delta }_{\un{H}}(h)=x^1X^1h_1r^2g^2S(x^2Y^1R^2y^2X^3_1)\ot
x^3_1Y^2R^1y^1X^2h_2r^1g^1S(x^3_2Y^3y^3X^3_2).\label{fff4}
\end{eqnarray}
Therefore, by (\ref{fff4}) and (\ref{qf2}), for any $\chi \in H^*$ we have
\begin{eqnarray*}
\un{\Delta }_{\un{H}}({\cal Q}(\chi ))&=&
<\chi , \tqra Z^1\mathfrak{R}^2\mathfrak{r}^1p^1>
x^1X^1\tilde{q}^2_{(1, 1)}Z^2_1\mathfrak{R}^1_1\mathfrak{r}^2_1
p^2_1S(\tilde{q}^2_2Z^3)_1r^2g^2S(x^2Y^1R^2y^2X^3_1)\\
&&\ot x^3_1Y^2R^1y^1X^2\tilde{q}^2_{(1, 2)}Z^2_2
\mathfrak{R}^1_2\mathfrak{r}^2_2p^2_2S(\tilde{q}^2_2Z^3)_2r^1g^1
S(x^3_2Y^3y^3X^3_2)\\
{\rm (\ref{qt3}, \ref{ca}, \ref{pr}, \ref{q1})}&=&
<\chi , \tqra Z^1\mathfrak{R}^2\mathfrak{r}^1V^1(T^1_1p^1)_1P^1>\\
&&x^1X^1(\tilde{q}^2_1Z^2)_1\mathfrak{R}^1_1r^2\mathfrak{r}^2_2V^3
T^1_2p^2S(x^2Y^1R^2y^2(X^3\tilde{q}^2_2)_1Z^3_1T^2)\\
&&\ot x^3_1Y^2R^1y^1X^2(\tilde{q}^2_1Z^2)_2\mathfrak{R}^1_2r^1
\mathfrak{r}^2_1V^2(T^1_1p^1)_2P^2S(x^3_2Y^3y^3
(X^3\tilde{q}^2_2)_2Z^3_2T^3)\\
{\rm (\ref{fff3}, \ref{qt2}, \ref{qt1}, \ref{fff1})}&=&
<\chi , S(v^1)\tqra v^2_1\mathfrak{R}^2t^3\mathfrak{r}^1_2
\mathbf{R}^1z^1_1P^1>\\
&&x^1X^1\tilde{q}^2_1v^2_{(2, 1)}\mathfrak{R}^1_1t^1
\mathfrak{r}^2z^2S(x^2Y^1R^2y^2X^3_1v^3_1z^3_1\tpra )\\
&&\ot x^3_1Y^2R^1y^1X^2\tilde{q}^2_2v^2_{(2, 2)}
\mathfrak{R}^1_2t^2\mathfrak{r}^1_1\mathbf{R}^2z^1_2P^2
S(x^3_2Y^3y^3X^3_2v^3_2z^3_2\tprb )\\
{\rm (\ref{qt1}, \ref{qt3}, \ref{qt1})}&=&
<\chi , S(v^1)\tqra v^2_1T^1\mathfrak{R}^2\mathfrak{r}^1
t^1V^1{\cal R}^2\mathbf{R}^1z^1_1P^1>\\
&&x^1X^1\tilde{q}^2_1v^2_{(2, 1)}T^2\mathfrak{R}^1
\mathfrak{r}^2t^2r^2V^3z^2S(x^2Y^1R^2y^2X^3_1v^3_1z^3_1\tpra )\\
&&\ot x^3_1Y^2R^1y^1X^2\tilde{q}^2_2v^2_{(2, 2)}T^3t^3r^1V^2
{\cal R}^1\mathbf{R}^2z^1_2P^2
S(x^3_2Y^3y^3X^3_2v^3_2z^3_2\tprb )\\
{\rm (\ref{q1}, \ref{fff1}, \ref{ql1}, \ref{q3})}&=&
<\chi , S(v^1)\tqra v^2_1X^1_1\mathfrak{R}^2\mathfrak{r}^1
t^1V^1{\cal R}^2\mathbf{R}^1z^1_1P^1>\\
&&x^1\tqrb v^2_2X^1_2\mathfrak{R}^1
\mathfrak{r}^2t^2r^2V^3z^2
S(x^2Y^1R^2y^2v^3_{(2, 1)}X^3_1z^3_1\tpra )\\
&&\ot x^3_1Y^2R^1y^1v^3_1X^2t^3r^1V^2
{\cal R}^1\mathbf{R}^2z^1_2P^2
S(x^3_2Y^3y^3v^3_{(2, 2)}X^3_2z^3_2\tprb )\\
{\rm (\ref{qt3}, \ref{q3}, \ref{ql1})}&=&
<\chi , S(v^1)\tqra v^2_1\mathfrak{R}^2\mathfrak{r}^1
t^1Z^1{\cal R}^2\mathbf{R}^1P^1>\\
&&x^1\tqrb v^2_2\mathfrak{R}^1
\mathfrak{r}^2t^2X^1r^2z^2
S(x^2Y^1R^2y^2(v^3t^3)_{(2, 1)}X^3_1z^3_1\tpra )\\
&&\ot x^3_1Y^2R^1y^1(v^3t^3)_1X^2r^1z^1Z^2
{\cal R}^1\mathbf{R}^2P^2
S(x^3_2Y^3y^3(v^3t^3)_{(2, 2)}X^3_2z^3_2\tprb Z^3)\\
{\rm (\ref{q1}, \ref{qt3}, \ref{q3}, \ref{qt2})}&=&
<\chi , S(v^1)\tqra v^2_1\mathfrak{R}^2\mathfrak{r}^1
t^1Z^1{\cal R}^2\mathbf{R}^1P^1>\\
&&x^1\tqrb v^2_2\mathfrak{R}^1
\mathfrak{r}^2t^2T^1X^1R^2_1V^2y^1_2z^2
S(x^2(v^3t^3)_1Y^1T^2_1X^2R^2_2V^3y^2z^3_1\tpra )\\
&&\ot x^3_1(v^3t^3)_{(2, 1)}Y^2T^2_2X^3R^1V^1y^1_1z^1Z^2
{\cal R}^1\mathbf{R}^2P^2
S(x^3_2(v^3t^3)_{(2, 2)}Y^3T^3y^3z^3_2\tprb Z^3)\\
{\rm (\ref{q3}, \ref{ql}, \ref{q5})}&=&
<\chi , S(v^1)\tqra v^2_1\mathfrak{R}^2\mathfrak{r}^1
t^1Z^1{\cal R}^2\mathbf{R}^1P^1>
x^1\tqrb v^2_2\mathfrak{R}^1
\mathfrak{r}^2t^2X^1\b S(x^2(v^3t^3)_1X^2)\\
&&\ot x^3_1(v^3t^3)_{(2, 1)}X^3_1Z^2{\cal R}^1
\mathbf{R}^2P^2S(x^3_2(v^3t^3)_{(2, 2)}X^3_2Z^3)\\
{\rm (\ref{q3}, \ref{q5}, \ref{ql1}, \ref{qt3})}&=&
<\chi , S(t^1x^1)\tqra \mathfrak{R}^2\mathfrak{r}^1
t^2_1x^2_{(1, 1)}z^1Z^1{\cal R}^2\mathbf{R}^1P^1>
\tqrb \mathfrak{R}^1\mathfrak{r}^2
t^2_2x^2_{(1, 2)}z^2\b S(t^3x^2_2z^3)\\
&&\ot x^3_1Z^2{\cal R}^1
\mathbf{R}^2P^2S(x^3_2Z^3)\\
{\rm (\ref{q1}, \ref{q5}, \ref{qt3}, \ref{qf2})}&=&
{\cal Q}(\chi _1\lh S(x^1))\ot
\chi _2(x^2Z^1{\cal R}^2\mathbf{R}^1P^1)x^3_1Z^2
{\cal R}^1\mathbf{R}^2P^2S(x^3_2Z^3).
\end{eqnarray*}
On the other hand, by (\ref{covhs}) we have
\begin{eqnarray*}
({\cal Q}\ot {\cal Q})(\un{\Delta }_{\un{H^*}}(\chi ))&=&
{\cal Q}(\chi _1\lh S(x^1))\ot {\cal Q}(x^3_2X^3\rh \chi _2\lh
x^2X^1\b S(x^3_1X^2))\\
{\rm (\ref{qf2}, \ref{qr1}, \ref{qt3})}&=&
{\cal Q}(\chi _1\lh S(x^1))\ot
<\chi _2, x^2X^1\b S(x^3_1X^2)\tQra Z^1(x^3_2X^3)_{(1, 1)}
{\cal R}^2\mathbf{R}^1P^1>\\
&&\widetilde{Q}^2_1Z^2(x^3_2X^3)_{(1, 2)}
{\cal R}^1\mathbf{R}^2P^2S(\widetilde{Q}^2_2Z^3(x^3_2X^3)_2)\\
{\rm (\ref{q1}, \ref{ql1}, \ref{ql})}&=&
{\cal Q}(\chi _1\lh S(x^1))\ot
<\chi _2, x^2S(\tpra )\tQra \tilde{p}^2_1Z^1
{\cal R}^2\mathbf{R}^1P^1>\\
&&x^3_1(\widetilde{Q}^2\tilde{p}^2_2)_1Z^2
{\cal R}^1\mathbf{R}^2P^2
S(x^3_2(\widetilde{Q}^2\tilde{p}^2_2)_2Z^3)\\
{\rm (\ref{pql})}&=&
{\cal Q}(\chi _1\lh S(x^1))\ot
<\chi _2, x^2Z^1{\cal R}^2\mathbf{R}^1P^1>x^3_1Z^2
{\cal R}^1\mathbf{R}^2P^2S(x^3_2Z^3).
\end{eqnarray*}
So ${\cal Q}$ is a coalgebra map since $(\un{\va}_{\un{H}}
\circ {\cal Q})(\chi )=\chi (\a )=\un{\va}_{\un{H^*}}$
and this finishes our proof.
\end{proof}
Let ${\cal C}$ be a braided category with left duality.
For any two objects $M, N\in {\cal C}$
there exists a canonical isomorphism in ${\cal C}$,
$M^*\ot N^*\stackrel{\sigma ^*_{M, N}}{\longrightarrow}
(M\ot N)^*$. In fact,
$\sigma ^*_{M, N}=\phi ^*_{N, M}\circ c^{-1}_{N^*, M^*}$,
where $\phi ^*_{N, M}: N^*\ot M^*\ra (M\ot N)^*$ is given by the
following compositions:
\begin{eqnarray}
N^*\ot M^*&&
\stackrel{id_{N^*\ot M^*}\ot l_{N^*\ot M^*}}{\longrightarrow}
(N^*\ot M^*)\ot \un{1}
\stackrel{id_{N^*\ot M^*}\ot coev_{M\ot N}}{\longrightarrow}
(N^*\ot M^*)\ot ((M\ot N)\ot (M\ot N)^*)\nonumber\\
&&\stackrel{a^{-1}_{N^*\ot M^*, M\ot N, (M\ot N)^*}}{\longrightarrow}
((N^*\ot M^*)\ot (M\ot N))\ot (M\ot N)^*\nonumber\\
&&\stackrel{a_{N^*, M^*, M\ot N}\ot id_{(M\ot N)^*}}{\longrightarrow}
(N^*\ot (M^*\ot (M\ot N)))\ot (M\ot N)^*\nonumber\\
&&\stackrel{(id_{N^*}\ot a^{-1}_{M^*, M, N})\ot id_{(M\ot N)^*}}{\longrightarrow}
(N^*\ot ((M^*\ot M)\ot N))\ot (M\ot N)^*\nonumber\\
&&\stackrel{(id_{N^*}\ot (ev_M\ot id_N))\ot id_{(M\ot N)^*}}{\longrightarrow}
(N^*\ot (\un{1}\ot N))\ot (M\ot N)^*\nonumber\\
&&\stackrel{(id_{N^*}\ot r^{-1}_{N})\ot id_{(M\ot N)^*}}{\longrightarrow}
(N^*\ot N)\ot (M\ot N)^*\nonumber\\
&&\stackrel{ev_N\ot id_{(M\ot N)^*}}{\longrightarrow}
\un{1}\ot (M\ot N)^*
\stackrel{r^{-1}_{(M\ot N)^*}}{\longrightarrow}
(M\ot N)^*.\label{phir}
\end{eqnarray}
The morphism $\phi ^*_{N, M}$  is an isomorphism.
One can compute its inverse in the same manner as above, see
\cite{ag}, \cite{ta}, \cite{bcp}.
Hence, $\sigma ^*_{M, N}$ is an isomorphism and
$\sigma ^{*-1}_{M, N}=c_{N^*, M^*}\circ \phi ^{*-1}_{N, M}$. \\
Also, following \cite{Kassel95}, for any morphism $\nu : M\ra N$ in ${\cal C}$,
we can define the transpose of $\nu $ as being
\begin{eqnarray}
\nu ^* : &&N^*\stackrel{l_{N^*}}{\longrightarrow}N^*\ot \un{1}
\stackrel{id_{N^*}\ot coev_M}{\longrightarrow}N^*\ot (M\ot M^*)
\stackrel{id_{N^*}\ot (\nu \ot id_{M^*})}{\longrightarrow}
N^*\ot (N\ot M^*)\nonumber \\
&&\stackrel{a^{-1}_{N^*, N, M^*}}{\longrightarrow}(N^*\ot N)\ot M^*
\stackrel{ev_N\ot id_{M^*}}{\longrightarrow}\un {1}\ot M^*
\stackrel{r^{-1}_{M^*}}{\longrightarrow}M^*.\label{rt}
\end{eqnarray}
Let now $(B, \un{m}_B, \un{\Delta }_B, \un{S}_B)$ be a braided Hopf algebra
in ${\cal C}$ and $B^*$ the categorical left dual of $B$ in ${\cal C}$.
Then $B^*$ is also a braided Hopf algebra in ${\cal C}$ with
multiplication $\un{m}_{B^*}$, comultiplication $\un{\Delta }_{B^*}$,
antipode $\un{S}_{B^*}$, unit $\un{u}_{B^*}$ and counit
$\un{\va }_{B^*}$ defined by (see \cite{maj}, p. 489)
\begin{eqnarray}
&&\un{m}_{B^*} : B^*\ot B^*\stackrel{\phi ^*_{B, B}}{\longrightarrow}(B\ot B)^*
\stackrel{\un{\Delta }_B^*}{\longrightarrow}B^*,\label{mbra}\\
&&\un{\Delta }_{B^*} : B^*\stackrel{\un{m}_B^*}{\longrightarrow}(B\ot B)^*
\stackrel{\phi ^{*-1}_{B, B}}{\longrightarrow}B^*\ot B^*,\label{combra}\\
&&\un{S}_{B^*}=\un{S}_B^*,\mbox{${\;}$}
\un{u}_{B^*}=\un{\va }_B^*, \mbox{${\;}$}
\un{\va }_{B^*}=\un{u}_B^*.\label{rebra}
\end{eqnarray}
Suppose now that $(H, R)$ is a QT quasi-Hopf algebra, and that $M$, $N$
are two finite dimensional left $H$-modules. Denote by
$\{{}_im\}_{i=\ov{1, s}}$ and $\{{}^im\}_{i=\ov{1, s}}$ dual bases
in $M$ and $M^*$, and by $\{{}_jn\}_{j=\ov{1, t}}$
and ${\{}^jn\}_{j=\ov{1, t}}$
dual bases in $N$ and $N^*$, respectively. By \cite{bcp}, in
this particular case we have that the morphism (\ref{phir}) is given by
\begin{equation}\label{sphi}
\phi ^*_{N, M}(n^*\ot m^*)(m\ot n)=
<m^*, f^1\cd m><n^*, f^2\cd n>
\end{equation}
for all $m^*\in M^*$, $n^*\in N^*$, $m\in M$ and $n\in N$. Its
inverse is defined by
\begin{equation}\label{sat}
\phi ^{*-1}_{N, M}(\mu )=<\mu , g^1\cd {}_im\ot g^2\cd {}_jn>
{}^jn\ot {}^im
\end{equation}
for any $\mu \in (M\ot N)^*$.
Also, the morphism $\nu ^*$ defined by (\ref{rt}) coincide with the usual
transpose map of $\nu $, i.e. $\nu ^*(n^*)=n^*\circ \nu $.\\
Therefore, if $(H, R)$ is finite dimensional then the categorical
left dual of $\un{H}$ has a braided Hopf algebra structure in ${}_H{\cal M}$.
We denote $\un{H}^*$ with this dual Hopf algebra structure by $(\un{H})^*$.
By the above, $(\un{H})^*$ is a left $H$-module via
\begin{equation}\label{hmscd}
(h\succ \chi )(h^{'})=\chi (S(h)\tr h^{'}),
\mbox{${\;\;}$$\forall $ $h, h^{'}\in H$, $\chi \in H^*$.}
\end{equation}
By (\ref{mbra}-\ref{rebra}) the structure of $(\un{H})^*$
as a Hopf algebra in ${}_H{\cal M}$ is given by the formulas
\begin{eqnarray}
&&(\chi * \psi )(h)=<\chi , f^2\tr h_{\un{2}}>
<\psi , f^1\tr h_{\un{1}}>,\label{mscd}\\
&&1_{(\un{H})^*}=\va ,\label{uscd}\\
&&\un{\Delta }_{(\un{H})^*}(\chi )=<\chi , (g^1\tr {}_ie)
\bullet (g^2\tr {}_je)>{}^je\ot {}^ie,\label{coscd}\\
&&\un{\va }_{(\un{H})^*}(\chi )=\chi (\b ),\label{cscd}\\
&&\un{S}_{(\un{H})^*}(\chi )=\chi \circ \un{S},\label{ascd}
\end{eqnarray}
where $\{{}_ie\}_{i=\ov{1, n}}$ and $\{{}^ie\}_{i=\ov{1, n}}$
are dual bases in $H$ and $H^*$.\\
Following \cite{bcp}, if $(H, R)$ is a finite dimensional QT quasi-Hopf algebra
then the usual liniar dual of $H$ has in ${}_H{\cal M}$
three braided Hopf algebra
structures. Two of them are the left and the right categorical dual
of $\un{H}$ in the sense of
Takeuchi, \cite{ta} (in the Hopf case, the same point of view was used in
\cite{ag}), and the third one is obtained in \cite{bc} by using
the structure of a quasi-Hopf algebra with a projection given in \cite{bn2}.
Now, by the above, $H^*$ has a fourth braided Hopf algebra structure in
${}_H{\cal M}$.

\begin{proposition}\label{pr4.2}
Let $(H, R)$ be a finite dimensional QT quasi-Hopf algebra,
$\un{H}$ the associated enveloping algebra
braided group of $H$, $(\un{H})^*$ the dual Hopf algebra structure
of $\un{H}$ in ${}_H{\cal M}$, and $\un{H^*}$ the function algebra braided group
associated to $H^*$. Then the map
$\l : (\un{H})^*\ra \un{H^*}$ given for all $\chi \in H^*$ by
\begin{equation}\label{lam}
\l (\chi  )=\smi (g^1)\rh \chi \circ S\lh g^2
\end{equation}
is a braided Hopf algebra isomorphism. Here $g^1\ot g^2$ is the
inverse of the Drinfeld twist $f$, see (\ref{g}).
\end{proposition}
\begin{proof}
The map $\l $ is left $H$-linear since
\begin{eqnarray*}
(h\blacktriangleright \l (\chi ))(h^{'})&=&<\l (\chi ), S(h_1)h^{'}h_2>\\
&=&<\chi , g^1S(h_2)S(h^{'})S(g^2S(h_1))>\\
{\rm (\ref{ca}, \ref{s1})}&=&<\chi , S(h)\tr (g^1S(g^2h^{'}))>\\
{\rm (\ref{hmscd})}&=&<(h\succ \chi )\circ S, g^2h^{'}\smi (g^1)>
=\l (h\succ \chi )(h^{'})
\end{eqnarray*}
for all $h, h^{'}\in H$ and $\chi \in H^*$. Thus, we only have to prove
that $\l $ is an algebra and coalgebra morphism, and that it is bijective.
Firstly, for all $\chi , \psi \in H^*$ and $h\in H$ we compute
\begin{eqnarray*}
\l (\chi * \psi )(h)&=&
<\chi , f^2\tr (g^1S(g^2h))_{\un{2}}><\psi , f^1\tr (g^1S(g^2h))_{\un{1}}>\\
{\rm (\ref{und}, \ref{s1}, \ref{ca})}&=&
<\chi , f^2x^3R^1\tr y^1X^2g^1_2G^2S(y^2X^3_1g^2_1h_1)>\\
&&<\psi , f^1_1x^1X^1g^1_1G^1S(f^1_2x^2R^2y^3X^3_2g^2_2h_2)> \\
{\rm (\ref{g2}, \ref{pf}, \ref{q1}, \ref{ca})}&=&
<\chi , f^2x^3R^1\tr G^2_{(1, 1)}y^1g^1S(G^2_{(1, 2)}y^2g^2_1
\mathfrak{G}^1S(X^1_2)F^1h_1X^2)>\\
&&<\psi , f^1_1x^1G^1S(f^1_2x^2R^2G^2_2y^3g^2_2\mathfrak{G}^2S(X^1_1)
F^2h_2X^3)>\\
{\rm (\ref{g2}, \ref{pf}, \ref{s1}, \ref{qt3})}&=&
<\chi , f^2x^3G^2_2R^1\mathfrak{G}^1\tr g^1S(g^2S(X^1_2y^2)F^1h_1X^2y^3)>\\
&&<\psi , f^1_1x^1G^1S(f^1_2x^2G^2_1R^2\mathfrak{G}^2S(X^1_1y^1)F^2h_2X^3)>\\
{\rm (\ref{ext}, \ref{g2}, \ref{pf}, \ref{s1}, \ref{ca})}&=&
<\chi , g^1S(g^2S(X^1_2y^2R^1_1x^1_1)F^1h_1X^2y^3R^1_2x^1_2)>\\
&&<\psi , G^1S(G^2S(X^1_1y^1R^2x^2)F^2h_2X^3x^3)>,
\end{eqnarray*}
and, on the other hand, by (\ref{movhs}) we have
\begin{eqnarray*}
(\l (\chi )\un{\cdot }\l (\psi ))(h)&=&
<\chi , g^1S(g^2S(x^1X^1)f^2R^1h_1x^3_1Y^2r^1y^1X^2)>\\
&&<\psi , G^1S(G^2S(x^2Y^1r^2y^2X^3_1)f^1R^2h_2x^3_2Y^3y^3X^3_2)>\\
{\rm (\ref{ext}, \ref{q3}, \ref{qt3})}&=&
<\chi , g^1S(g^2S(Y^1_2R^1z^1x^1X^1)f^1h_1Y^2z^3r^1x^2_1y^1X^2)>\\
&&<\psi , G^1S(G^2S(Y^1_1R^2z^2r^2x^2_2y^2X^3_1)f^2h_2Y^3x^3y^3X^3_2)>\\
{\rm (\ref{q3}, \ref{qt1})}&=&
<\chi , g^1S(g^2S(Y^1_2z^2R^1_1x^1_1)f^1h_1Y^2z^3R^1_2x^1_2)>\\
&&<\psi , G^1S(G^2S(Y^1_1z^1R^2y^2)f^2h_2Y^3y^3)>,
\end{eqnarray*}
as needed. It is not hard to see that $\l (1_{(\un{H})^*})=1_{\un{H^*}}$,
so $\l $ is an algebra morphism. In order to prove that $\l $ is a
coalgebra map we need the following formula
\begin{equation}\label{fff5}
S(g^1)\a g^2=S(\b )
\end{equation}
which can be found in \cite{bn3}. Now, $\l $ is
a coalgebra morphism since
\begin{eqnarray*}
(\l \ot \l )(\un{\Delta }_{(\un{H})^*}(\chi ))&=&
<\chi , (g^1\tr {}_ie)\bullet (g^2\tr {}_je)>\l ({}^je)\ot \l ({}^ie)\\
{\rm (\ref{lam}, \ref{ma}, \ref{s1})}&=&
<\chi , X^1g^1_1\mathfrak{G}^1S(x^1X^2g^1_2\mathfrak{G}^2{}_ie)\a
x^2X^3_1g^2_1G^1S(x^3X^3_2g^2_2G^2{}_je)>{}^je\ot {}^ie\\
{\rm (\ref{g2}, \ref{pf}, \ref{q1}, \ref{q5})}&=&
<\chi , \mathfrak{G}^1S(g^1S(X^2x^3){}_ieX^3)\a g^2
S(\mathfrak{G}^2S(X^1_1x^1){}_jeX^1_2x^2)>{}^je\ot {}^ie\\
{\rm (\ref{fff5}, \ref{lam})}&=&
<\l (\chi ), S(X^1_1x^1){}_jeX^1_2x^2\b S(X^2x^3){}_ieX^3)>{}^je\ot
{}^ie\\
{\rm (\ref{q3}, \ref{q5}, \ref{covhs})}&=&
\l (\chi )_1\lh S(X^1_1x^1)\ot
X^3\rh \l (\chi )_2\lh X^1_2x^2\b S(X^2x^3)
=\un{\Delta }_{\un{H^*}}(\l (\chi ))
\end{eqnarray*}
for all $\chi \in H^*$, and since the definitions of counits
imply $\un{\va }_{\un{H^*}}\circ \l =\un{\va }_{(\un{H})^*}$.
It is easy to see that $\l $ is bijective with inverse
$\l ^{-1}(\chi )=S(f^2)\rh \chi \circ \smi \lh f^1$, for all
$\chi \in H^*$. Thus, the proof is complete.
\end{proof}
Summarizing the results of this Section we can give the true
meaning of the map ${\cal Q}: H^*\ra H$ defined in (\ref{qf1}).
It is a morphism of braided groups from $\un{H^*}$, the function algebra
braided group associated to $H^*$, to $\un{H}$,
the associated enveloping algebra braided group of $H$.
When $H$ is factorizable in the sense that
the map ${\cal Q}$ is bijective then
${\cal Q}: \un{H^*}\cong \un{H}$ as braided Hopf algebras.
In other words, the function algebra
braided group associated to $H^*$ and
the associated enveloping algebra braided group of $H$
are categorical self dual, cf. Proposition \ref{pr4.2}.
\section{$D(H)$ when $H$ is factorizable}
\setcounter{equation}{0}
Schneider's Theorem \cite{sch} asserts that the quantum
double of a finite dimensional factorizable Hopf algebra
$H$ is a $2$-cocycle twist of the usual (componentwise) tensor
product Hopf algebra $H\ot H$. We note that this result
also appears in \cite{rs} without an explicit proof.
The aim of this Section is to give a proof of a similar
result in the quasi-Hopf case. Our approach is based on the
methods developed in \cite{sch}, \cite{bc1}.\\
Throughout, $(H, R)$ will be a finite dimensional
QT quasi-Hopf algebra and $D(H)$ its quantum double. When
there is no danger of confusion
the elements $\chi \Join h$ of $D(H)$ will be simply denoted by
$\chi h$. Since $H^*$ can be viewed only as a
$k$-linear subspace of $D(H)$ we will denote by
\begin{eqnarray*}
\chi _{(1)}\ot \chi _{(2)}&:=&\Delta _D(\chi \Join 1_H)\\
{\rm (\ref{cdd})}&=&
(\va \Join X^1Y^1)(p^1_1x^1\rh \chi _2\lh Y^2\smi (p^2)\Join
p^1_2x^2)\ot X^2_1\rh \chi _1\lh \smi (X^3)\Join X^2_2Y^3x^3.
\end{eqnarray*}
In our notation, for all $\chi \in H^*$ and $h\in H$,
the comultiplication $\Delta _D$ of $D(H)$ comes out as
$$
\Delta _D(\chi h)=\chi _{(1)}h_1\ot \chi _{(2)}h_2.
$$
By \cite[Lemma 3.1]{bc}, there exists a quasi-Hopf algebra
projection $\pi : D(H)\ra H$ covering the canonical inclusion
$i_D: H\ra D(H)$. More precisely, if $R=R^1\ot R^2$ is the
$R$-matrix of $H$ then $\pi $ is defined by
\begin{equation}\label{pir}
\pi (\chi \Join h)=\chi (q^2R^1)q^1R^2h,
\end{equation}
where $q_R=q^1\ot q^2$ is the element defined by (\ref{qr}). We have
that $\pi $ is a quasi-Hopf algebra morphism and
$\pi \circ i_D=id_H$.\\
It is not hard to see that $\tilde{R}:=R^{-1}_{21}=\ov{R}^2\ot \ov{R}^1$
is another $R$-matrix for $H$. So, as in the Hopf case,
there is always a second projection
$\widetilde{\pi }: D(H)\ra H$ covering the canonical
inclusion $i_D$. Explicitly, the morphism $\widetilde{\pi }$
is given by
\begin{equation}\label{tpir}
\widetilde{\pi }(\chi \Join h)=\chi (q^2\ov{R}^2)q^1\ov{R}^1h
\end{equation}
for all $\chi \in H^*$ and $h\in H$.\\
Let $H$ and $A$ be two quasi-bialgebras. Recall that a
quasi-bialgebra map between $H$ and $A$ is an algebra map
$\nu : H\ra A$ which intertwines
the quasi-coalgebra structures, respects the counits,
and satisfies $(\nu \ot \nu \ot \nu )(\Phi )=\Phi _A$.
Following \cite{bn2}, when $H$ is a quasi-Hopf algebra
we define
\begin{equation}\label{coinvnu}
H^{co(\nu )}=\{h\in H\mid
h_1\ot \nu (h_2)=
x^1hS(x^3_2X^3)f^1\ot \nu (x^2X^1\b S(x^3_1X^2)f^2)\}
\end{equation}
as being the set of coinvariants of $H$ relative to $\nu $. Finally,
if $A$ is also a quasi-Hopf algebra
then $\nu $ is a quasi-Hopf algebra morphism if, in addition,
$\nu (\a )=\a _A$, $\nu (\b )=\b _A$ and $S_A\circ \nu =
\nu \circ S_A$.

\begin{lemma}\label{lm5.1}
Let $(H, R)$ be a finite dimensional
QT quasi-Hopf algebra, and $\pi $ and $\widetilde{\pi }$
the quasi-Hopf algebra morphisms defined by (\ref{pir}) and
(\ref{tpir}), respectively. Let $j: D(H)^{co(\pi )}\ra D(H)$
be the inclusion map and $\Psi : H^*\ra D(H)^{co(\pi )}$
defined by
\begin{equation}\label{psi}
\Psi (\chi )=\chi _{(1)}\b S(\pi (\chi _{(2)}))
\end{equation}
for all $\chi \in H^*$. Then the following assertions hold:
\begin{itemize}
\item[1)] $\Psi $ is well defined and bijective.
\item[2)] If $\ov{{\cal Q}}$ is the map defined by (\ref{oqf}) then
$S\circ \ov{\cal Q}=\widetilde{\pi }\circ j\circ \Psi $. In particular,
$\ov{\cal Q}$ is bijective if and only if
$\widetilde{\pi }\mid _{D(H)^{co(\pi )}}$ is bijective.
\end{itemize}
\end{lemma}
\begin{proof}
1) For all $\chi \in H^*$ we have
\begin{eqnarray*}
(id\ot \pi )\Delta _D(\Psi (\chi ))&=&\chi _{((1), (1))}\b _1
S(\pi (\chi _{(2)}))_1\ot \chi _{((1), (2))}\b _2S(\pi (\chi _{(2)}))_2
\end{eqnarray*}
where we use the Sweedler type notation
\begin{eqnarray*}
&&(\Delta _D\ot id)(\Delta _D(\chi ))=
\chi _{((1), (1))}\ot \chi _{((1), (2))}\ot \chi _{(2)},~~
(id \ot \Delta _D)(\Delta _D(\chi ))=\chi _{(1)}\ot
\chi _{((2), (1))}\ot \chi _{((2), (2))}.
\end{eqnarray*}
Now, since $H$ is a quasi-Hopf subalgebra of $D(H)$ and
$\pi $ is a quasi-Hopf algebra morphism such that
$\pi (h)=h$ for any $h\in H$, by similar computations as in
\cite[Lemma 4.2]{bn2} or \cite[Lemma 3.6]{bc1}, one can prove
that $\Psi (\chi )\in D(H)^{co(\pi)}$, so $\Psi $ is well
defined. We claim that the inverse of $\Psi $,
$\Psi ^{-1}: D(H)^{co(\pi )}\ra H^*$, is given for all
${\bf D}\in D(H)^{co(\pi )}$ by the formula
$$
\Psi ^{-1}({\bf D})=(id\ot \va )({\bf D}).
$$
Indeed, $\Psi ^{-1}$ is a left inverse since
\begin{eqnarray*}
&&(\Psi ^{-1}\circ \Psi )(\chi )=<id\ot \va ,
\chi _{(1)}\b S(\pi (\chi _{(2)}))>
=(id\ot \va )(\chi _{(1)})\va _D(\chi _{(2)})=
(id\ot \va )(\chi \ot 1)=\chi ,
\end{eqnarray*}
for all $\chi \in H^*$. It is also a right inverse. If
${\bf D}={}_i\chi {}_ih\in D(H)^{co(\pi )}$ then
$$
{}_i\chi _{(1)}{}_ih_1\ot \pi ({}_i\chi _{(2)}){}_ih_2
=x^1[{}_i\chi {}_ih]S(x^3_2X^3)f^1\ot x^2X^1\b S(x^3_1X^2)f^2
$$
in $D(H)\ot H$. Therefore,
\begin{eqnarray*}
(\Psi \circ \Psi ^{-1})({\bf D})&=&\va ({}_ih)\Psi ({}_i\chi )
=\va ({}_ih){}_i\chi _{(1)}\b S(\pi ({}_i\chi _{(2)}))\\
&=&{}_i\chi _{(1)}{}_ih_1\b S(\pi ({}_i\chi _{(2)}){}_ih_2)\\
&=&x^1[{}_i\chi {}_ih]S(x^3_2X^3)f^1\b S(x^2X^1\b S(x^3_1X^2)f^2)\\
&=&{}_i\chi {}_ih={\bf D},
\end{eqnarray*}
because of $f^1\b S(f^2)=S(\a )$, and (\ref{q5}, \ref{q6}).\\
${\;\;\:}$2) By (\ref{psi}, \ref{cdd}) and (\ref{pir}),
for any $\chi \in H^*$ we find that
\begin{eqnarray*}
\Psi (\chi )&=&(X^1Y^1)_{(1, 1)}p^1_1x^1\rh \chi \lh \smi (X^3)q^2R^1X^2_1Y^2
\smi ((X^1Y^1)_2p^2)\\
&&\Join (X^1Y^1)_{(1, 2)}p^1_2x^2\b S(q^1R^2X^2_2Y^3x^3)
\end{eqnarray*}
and, if we denote by $Q^1\ot Q^2$ another copy of $q_R$, and by
$r^1\ot r^2$ another copy of $R$, then by (\ref{tpir})
we compute that
\begin{eqnarray*}
(\widetilde{\pi }\circ j\circ \Psi )(\chi )&=&
<\chi , \smi (X^3)q^2R^1X^2_1Y^2\smi ((X^1Y^1)_2p^2)Q^2\ov{R}^2
(X^1Y^1)_{(1, 1)}p^1_1x^1>\\
&&Q^1\ov{R}^1(X^1Y^1)_{(1, 2)}p^1_2x^2\b S(q^1R^2X^2_2Y^3x^3)\\
{\rm (\ref{qt3}, \ref{qr1})}&=&
<\chi , \smi (X^3)q^2R^1X^2_1Y^2\smi (p^2)Q^2\ov{R}^2p^1_1x^1>
X^1Y^1Q^1\ov{R}^1p^1_2x^2\b S(q^1R^2X^2_2Y^3x^3)\\
{\rm (\ref{qt3}, \ref{pqr})}&=&
<\chi , \smi (X^3)q^2R^1X^2_1Y^2\ov{R}^2x^1>
X^1Y^1\ov{R}^1x^2\b S(q^1R^2X^2_2Y^3x^3)\\
{\rm (\ref{qt3}, \ref{qt1})}&=&
<\chi , \smi (X^3)q^2X^2_2R^1r^2Y^3\ov{R}^2>
X^1Y^1\ov{R}^1_1\b S(q^1X^2_1R^2r^1Y^2\ov{R}^1_2)\\
{\rm (\ref{q5}, \ref{qt4}, \ref{qt3})}&=&
<\chi , \smi (X^3)q^2R^1r^2X^2_2Y^3>
S(q^1R^2r^1X^2_1Y^2\smi (X^1Y^1\b ))\\
{\rm (\ref{ql}, \ref{oqf})}&=&(S\circ \ov{\cal Q})(\chi ),
\end{eqnarray*}
as needed. Since $H$ is finite dimensional the antipode $S$ is bijective, so
$\ov{\cal Q}$ is bijective if and only if
$\widetilde{\pi }\circ j$ is bijective. Thus, the proof is complete.
\end{proof}

For the next result we need the concept of right quasi-Hopf bimodule
introduced in \cite{hn3}, and the second Structure Theorem for right
quasi-Hopf bimodules proved in \cite{bc1}.\\
Let $H$ be a quasi-bialgebra, $M$ an $H$-bimodule and
$\rho : M\ra M\ot H$ an $H$-bimodule map. Then $(M, \rho )$
is called a right quasi-Hopf $H$-bimodule if the following relations
hold:
\begin{eqnarray}
&&(id\ot \va )\circ \rho =id,\label{qhbi1}\\
&&\Phi(\rho \ot id)(\r (m))=(id\ot \Delta )(\r (m))\Phi ,~
\mbox{$\forall$~$m\in M$.}\label{qhbi2}
\end{eqnarray}
A morphism between two right quasi-Hopf $H$-bimodules is an $H$-bimodule
map which is also right $H$-colinear (just like in the Hopf case).
${}_H{\cal M}^H_H$ is the category of right quasi-Hopf
$H$-bimodules and morphisms of right quasi-Hopf $H$-bimodules.\\
Let $H$ be a quasi-Hopf algebra and $M\in {}_H{\cal M}^H_H$. Following
\cite{bc1}, we define
\begin{equation}\label{ovcoinv}
M^{\ov{co(H)}}:=\{n\in M\mid \r (n)=x^1\cd n\cd S(x^3_2X^3)f^1
\ot x^2X^1\b S(x^3_1X^2)f^2\}
\end{equation}
and $\ov{E}: M\ra M$, by
\begin{equation}\label{ove}
\ov{E}(m)=m_{(0)}\cdot \b S(m_{(1)}),
\end{equation}
for all $m\in M$, where $\r (m):=m_{(0)}\ot m_{(1)}$.
From \cite[Lemma 3.6]{bc1} we have that
$Im(\ov{E})=M^{\ov{co(H)}}$, and that
$M^{\ov{co(H)}}$ is a left $H$-submodule of $M$, where $M$ is
considered a left $H$-module via the left adjoint action, that is
$h\tr m=h_1\cdot m\cdot S(h_2)$, for all $h\in H$ and $m\in M$.
Moreover, if $M^{\ov{co(H)}}\ot H$ is viewed as a right quasi-Hopf
$H$-bimodule via the structure
$$
h\cdot (n\ot h^{'})\cdot h^{''}=h_1\tr n\ot h_2h^{'}h^{''},~
\rho ^{'}(n\ot h)=x^1\tr n\ot x^2h_1\ot x^3h_2,
$$
then the map
$$
\ov{\nu }_M: M^{\ov{co(H)}}\ot H\ra M,~
\ov{\nu }_M(n\ot h)=X^1\cd n\cdot S(X^2)\a X^3h
$$
is an isomorphism of quasi-Hopf $H$-bimodules, cf. \cite[Theorem 3.7]{bc1}.
The inverse of $\ov{\nu }_M$ is given by the formula
$$
\ov{\nu }^{-1}_M(m)=\ov{E}(m_{(0)})\ot m_{(1)}.
$$
Suppose now that $H$ is a quasi-Hopf algebra and $\mu : M\ra N$
is a morphism between two right quasi-Hopf $H$-bimodules.
It is not hard to see that the restriction of $\mu $ defines a
left $H$-linear map between $M^{\ov{co(H)}}$ and $N^{\ov{co(H)}}$.
Moreover, if we denote this (co)restriction by $\mu _0$ then the
following diagram is commutative.
$$\begin{diagram}
M^{\ov{co(H)}}\ot H&\rTo{\ov{\nu }_M}&M\\
\dTo{\mu _0\ot id}& &\dTo{\mu }\\
N^{\ov{co(H)}}\ot H&\rTo{\ov{\nu }_N}&N
\end{diagram}$$
Consequently, the map $\mu $ is bijective if and only if
the map $\mu _0: M^{\ov{co(H)}}\ra N^{\ov{co(H)}}$ is bijective.

\begin{lemma}\label{lm5.2}
Let $D$, $A$ and $B$ three quasi-bialgebras and
$\vartheta , \upsilon , \kappa $ three quasi-bialgebra
morphisms as in the diagram below
$$\begin{diagram}
D&\rTo{\vartheta}&A\\
\dTo{\upsilon}\uTo_{\kappa}&\SE^{\zeta:=(\vartheta \ot \upsilon )\circ \Delta _D} & \\
B& & A\ot B.
\end{diagram}$$
Suppose that $\upsilon \circ \kappa =id_B$ and define $\zeta $ as above.
Then the following assertions hold:
\begin{itemize}
\item[1)] $D$ and $A\ot B$ are right quasi-Hopf $B$-bimodules
via the following structures
\begin{eqnarray*}
&&D\in {}_B{\cal M}^B_B : \left\{\begin{array}{lccccccr}
b\cd d\cd b^{'}=\kappa (b)d\kappa (b^{'})\\
\r _D(d)=d_1\ot \upsilon (d_2),
\end{array}\right. \\
&&A\ot B\in {}_B{\cal M}^B_B : \left\{\begin{array}{lccccccr}
b^{'}\cd (a\ot b)\cd b^{''}=\vartheta (\kappa (b^{'}_1))a
\vartheta (\kappa (b^{''}_1))\ot b^{'}_2bb^{''}_2\\
\r _{A\ot B}(a\ot b)=\vartheta (\kappa (x^1))a\vartheta (\kappa (X^1))
\ot x^2b_1X^2\ot x^3b_2X^3,
\end{array}\right.
\end{eqnarray*}
$a\in A$, $b, b^{'}, b^{''}\in B$, $d\in D$, and $\zeta $
becomes a quasi-Hopf $B$-bimodule morphism.
\item[2)] If $D$, $A$ and $B$ are quasi-Hopf algebras and
$\vartheta , \upsilon $ and $\kappa $ are quasi-Hopf algebra maps
then $D^{\ov{co(B)}}=D^{co(\upsilon)}$ and
\begin{eqnarray*}
&&(A\ot B)^{\ov{co(B)}}=\{\vartheta (\kappa (x^1))a\vartheta
(\kappa (S(x^3_2X^3)f^1))\ot x^2X^1\b S(x^3_1X^2)f^2\mid a\in A\}.
\end{eqnarray*}
\end{itemize}
\end{lemma}
\begin{proof}
Since no confusion is possible we will write without subscripts $D, A$ or $B$
in the tensor components of the reassociators of $D$, $A$ or $B$, respectively.
The same thing we will do when we write their inverses.\\
1) It is straightforward to show that with the above
structures $D$ is an object of ${}_B{\cal M}^B_B$, and that $A\ot B$ is a
$B$-bimodule. The map $\r _{A\ot B}$ is a $B$-bimodule map since
\begin{eqnarray*}
\r _{A\ot B}(b^{'}\cd (a\ot b)\cd b^{''})&=&
\rho _{A\ot B}(\vartheta (\kappa (b^{'}_1))a
\vartheta (\kappa (b^{''}_1))\ot b^{'}_2bb^{''}_2)\\
&=&\vartheta (\kappa (x^1b^{'}_1))a\vartheta (\kappa (b^{''}_1X^1))
\ot x^2b^{'}_{(2, 1)}b_1b^{''}_{(2, 1)}X^2\ot
\ot x^3b^{'}_{(2, 2)}b_2b^{''}_{(2, 2)}X^3\\
{\rm (\ref{q1})}&=&b^{'}_1\cdot (\vartheta (\kappa (x^1))a
\vartheta (\kappa (X^1))\ot x^2b_1X^2)\cd b^{''}_1\ot
b^{'}_2x^3b_2X^3b^{''}_2\\
&=&\Delta (b^{'})\r _{A\ot B}(a\ot b)\Delta (b^{''}),
\end{eqnarray*}
for all $a\in A$ and $b, b^{'}, b^{''}\in B$. Similar computations show
that
$$
\Phi ^{-1}(id\ot \Delta )(\rho _{A\ot B}(a\ot b))\Phi
=(\r _{A\ot B}\ot id)(\rho _{A\ot B}(a\ot b)),
$$
for all $a\in A$ and $b\in B$, so $A\ot B\in {}_B{\cal M}^B_B$.
Also, we can check directly that $\zeta $ becomes a morphism in
${}_B{\cal M}^B_B$, the details are left to the reader.\\
2) By definitions we have
\begin{eqnarray*}
D^{\ov{co(B)}}&=&\{d\in D\mid \r _D(d)=x^1\cd d\cd S(x^3_2X^3)f^1
\ot x^2X^1\b S(x^3_1X^2)f^2\}\\
&=&\{d\in D\mid d_1\ot \upsilon (d_2)=\kappa (x^1)d\kappa (S(x^3_2X^3)f^1)
\ot \upsilon (\kappa (x^2X^1\b S(x^3_1X^2)f^2))\}\\
&=&D^{co(\upsilon )}.
\end{eqnarray*}
Observe now that $a\ot b\in (A\ot B)^{\ov{co(B)}}$
if and only if
\begin{eqnarray}
&&\vartheta (\kappa (x^1))a\vartheta (\kappa (X^1))\ot x^2b_1X^2
\ot x^3b_2X^3\nonumber\\
&&\hspace*{1cm}=\vartheta (\kappa (x^1_1))a
\vartheta (\kappa ((S(x^3_2X^3)f^1)_1))
\ot x^1_2b(S(x^3_2X^3)f^1)_2\ot x^2X^1\b S(x^3_1X^2)f^2.\label{fff6}
\end{eqnarray}
If $a\ot b\in (A\ot B)^{\ov{co(B)}}$ applying $id\ot \va \ot id$ to the
equality (\ref{fff6}) we obtain
$$
a\ot b=\va (b)\vartheta (\kappa (x^1))a\vartheta (\kappa (S(x^3_2X^3)f^1))
\ot x^2X^1\b S(x^3_1X^2)f^2,
$$
so $(A\ot B)^{\ov{co(B)}}\subseteq
\{\vartheta (\kappa (x^1))a\vartheta
(\kappa (S(x^3_2X^3)f^1))\ot x^2X^1\b S(x^3_1X^2)f^2\mid a\in A\}$.
Conversely, if $\d =\d ^1\ot \d ^2$ is the element
defined by (\ref{gd})
and $F^1\ot F^2=\mathbf{F}^1\ot \mathbf{F}^2$
are other copies of $f$ then for all $a\in A$ we compute
\begin{eqnarray*}
&&\hspace*{-2.8cm}
\vartheta (\kappa (y^1))[\vartheta (\kappa (x^1))a\vartheta
(\kappa (S(x^3_2X^3)f^1))]\vartheta (\kappa (Y^1))\ot
y^2[x^2X^1\b S(x^3_1X^2)f^2]_1Y^2\ot
y^3[x^2X^1\b S(x^3_1X^2)f^2]_2Y^3\\
{\rm (\ref{gdf}, \ref{ca})}&=&
\vartheta (\kappa (y^1x^1))a\vartheta (\kappa (S(x^3_2X^3)f^1Y^1))\ot
y^2x^2_1X^1_1z^1\b S(x^3_{(1, 2)}X^2_2z^3_2Z^3)F^1f^2_1Y^2\\
&&\ot y^3x^2_2X^1_2z^2Z^1\b S(x^3_{(1, 1)}X^2_1z^3_1Z^2)F^2f^2_2Y^3\\
{\rm (twice~~\ref{q3})}&=&
\vartheta (\kappa (x^1_1y^1))a
\vartheta (\kappa (S((x^3y^3_2)_2X^3T^3)f^1Y^1))\ot
x^1_2y^2T^1\b S((x^3y^3_2)_{(1, 2)}X^2_2Z^3T^2)F^1f^2_1Y^2\\
&&\ot x^2y^3_1X^1Z^1\b S((x^3y^3_2)_{(1, 1)}X^2_1Z^2)F^2f^2_2Y^3\\
{\rm (\ref{q3}, \ref{q5}, \ref{q1})}&=&
\vartheta (\kappa (x^1_1y^1))a
\vartheta (\kappa (S(x^3_2z^3(y^3_{(2, 2)}X^3)_2T^3)f^1Y^1))\\
&&\ot x^1_2y^2T^1\b S(x^3_{(1, 2)}z^2(y^3_{(2, 2)}X^3)_1T^2)F^1f^2_1Y^2
\ot x^2y^3_1X^1\b S(x^3_{(1, 1)}z^1y^3_{(2, 1)}X^2)F^2f^2_2Y^3\\
{\rm (\ref{q1}, \ref{q5})}&=&
\vartheta (\kappa (x^1_1y^1))a
\vartheta (\kappa (S(x^3_2z^3X^3_2y^3_2T^3)f^1Y^1))
\ot x^1_2y^2T^1\b S(x^3_{(1, 2)}z^2X^3_1y^3_1T^2)F^1f^2_1Y^2\\
&&\ot x^2X^1\b S(x^3_{(1, 1)}z^1X^2)F^2f^2_2Y^3\\
{\rm (\ref{q1}, \ref{g2}, \ref{pf})}&=&
\vartheta (\kappa (x^1_1))[\vartheta (\kappa (y^1))a
\vartheta (\kappa (S(y^3_2T^3)\mathbf{F}^1))]
\vartheta (\kappa (g^1S(x^3_{(2, 2)}X^3_2)F^1f^1_1))\\
&&\ot x^1_2[y^2T^1\b S(y^3_1T^2)\mathbf{F}^2]
g^2S(x^3_{(2, 1)}X^3_1)F^2f^1_2
\ot x^2X^1\b S(x^3_1X^2)f^2\\
{\rm (\ref{ca})}&=&
\vartheta (\kappa (x^1_1))[\vartheta (\kappa (y^1))a
\vartheta (\kappa (S(y^3_2T^3)\mathbf{F}^1))]
\vartheta (\kappa ((S(x^3_2X^3)f^1)_1))\\
&&\ot x^1_2[y^2T^1\b S(y^3_1T^2)\mathbf{F}^2]
(S(x^3_2X^3_1)f^1)_2\ot x^2X^1\b S(x^3_1X^2)f^2,
\end{eqnarray*}
as needed. Therefore,
$\{\vartheta (\kappa (x^1))a\vartheta
(\kappa (S(x^3_2X^3)f^1))\ot x^2X^1\b S(x^3_1X^2)f^2\mid a\in A\}\subseteq
(A\ot B)^{\ov{co(B)}}$, and this finishes our proof.
\end{proof}

\begin{proposition}\label{pr5.3}
Let $D$ be a quasi-Hopf algebra, $A$ and $B$ two
quasi-bialgebras and $\vartheta : D\ra A$,
$\upsilon : D\ra B$ two quasi-bialgebra maps. Consider
$\zeta : D\ra A\ot B$ given by
$\zeta (d)=\vartheta (d_1)\ot \upsilon (d_2)$,
for all $d\in D$.
\begin{itemize}
\item[1)] Suppose that $(D, R)$ is quasitriangular and define
$\mathfrak{F}=\mathfrak{F}^1\ot \mathfrak{F}^2\in (A\ot B)^{\ot 2}$, by
\begin{equation}\label{sptwist}
\mathfrak{F}=\vartheta (Y^1_1x^1X^1y^1_1)\ot
\upsilon (Y^1_2x^2\ov{R}^1X^3y^2)\ot
\vartheta (Y^2x^3\ov{R}^2X^2y^1_2)\ot
\upsilon (Y^3y^3),
\end{equation}
where, as usual, $\ov{R}^1\ot \ov{R}^2$ is the inverse of the $R$-matrix
$R$ of $D$. Then $\mathfrak{F}$ is a twist on $A\ot B$ (here $A\ot B$ has
the componentwise quasi-bialgebra structure) and
$\zeta : D\ra (A\ot B)_{\mathfrak{F}}$ is a quasi-bialgebra morphism.
Moreover, if $A$ and $B$ are quasi-Hopf algebras
and $\vartheta $ and $\upsilon $ are quasi-Hopf algebra morphisms, then
$\zeta : D\ra (A\ot B)^{\mathfrak{U}}_{\mathfrak{F}}$ is a quasi-Hopf
algebra morphism, where
$\mathfrak{U}=\vartheta (R^2g^2)\ot \upsilon (R^1g^1)$.
\item[2)] Suppose that $A$ and $B$ are quasi-Hopf algebras, $\vartheta $
and $\upsilon $ are quasi-Hopf algebra morphisms,
and that there exists a quasi-Hopf
algebra map $\kappa : B\ra D$ such that $\upsilon \circ \kappa =id_B$.
Then $\zeta $ is a bijective map if and only if the restriction of
$\vartheta $ provides a bijection from $D^{co(\vartheta )}$ to $A$.
\end{itemize}
\end{proposition}
\begin{proof}
1) We have that $\zeta =(\vartheta \ot \upsilon )\circ \Delta _D$, so
clearly $\zeta $ is an algebra map. It also respects the
comultiplications. Indeed, applying (\ref{g1}),
twice (\ref{q1}), (\ref{qt3}), and then again (\ref{q1}) two times,
it is not hard to see that
$$
(\Delta _{(A\ot B)_{\mathfrak{F}}}\circ \zeta )(d)=
((\zeta \ot \zeta )\circ \Delta _D)(d)
$$
for all $d\in D$. Obviously, $\va _{A\ot B}\circ \zeta =\va _D$, so
$\zeta $ respects the counits. It remains to show that
$$(\zeta \ot \zeta \ot \zeta )(\Phi _D)=\Phi _{(A\ot B)_{\mathfrak{F}}}.
$$
This follows from a long, technical but straightforward computation, we leave
the details to the reader.\\
Suppose now that $A$ and $B$ are quasi-Hopf algebras and that $\vartheta $ and
$\upsilon $ are quasi-Hopf algebra morphisms. In this case,
$\zeta : D\ra (A\ot B)^{\mathfrak{U}}_{\mathfrak{F}}$
is also a quasi-bialgebra morphism
since $(A\ot B)^{\mathfrak{U}}_{\mathfrak{F}}=(A\ot B)_{\mathfrak{F}}$ as
quasi-bialgebras. Thus, we are left to show that
\begin{eqnarray*}
&&\zeta (\a )=\mathfrak{U}\a _{(A\ot B)_{\mathfrak{F}}},~~
\zeta (\b )=\b _{(A\ot B)_{\mathfrak{F}}}\mathfrak{U}^{-1},~~
(\zeta \circ S_D)(d)=\mathfrak{U}S_{A\ot B}(\zeta (d))\mathfrak{U}^{-1}
\end{eqnarray*}
for all $d\in D$. Take
$\mathfrak{F}^{-1}=\mathfrak{G}^1\ot \mathfrak{G}^2$ as being the inverse
of the twist $\mathfrak{F}$. By (\ref{g3}) and (\ref{sptwist}) we compute:
\begin{eqnarray*}
\a _{(A\ot B)_{\mathfrak{F}}}&=&S_{A\ot B}(\mathfrak{G}^1)\a _{A\ot B}
\mathfrak{G}^2\\
&=&\vartheta (S(Y^1_1x^1X^1y^1_1)\a Y^1_2x^2R^2X^3y^2)\ot
\upsilon (S(Y^2x^3R^1X^2y^1_2)\a Y^3y^3)\\
{\rm (\ref{q5}, \ref{qt2})}&=&
\vartheta (S(R^2_1X^2\ov{R}^2y^1_1)\a R^2_2X^3y^2)\ot
\upsilon (S(R^1X^1\ov{R}^1y^1_2)\a y^3)\\
{\rm (\ref{q5}, \ref{qt4}, \ref{qt3})}&=&
\vartheta (S(X^2y^1_2\ov{R}^2)\a X^3y^2)\ot
\upsilon (S(X^1y^1_1\ov{R}^1)\a y^3)\\
{\rm (\ref{gd}, \ref{gdf})}&=&
\vartheta (S(\ov{R}^2)\gamma ^1)\ot
\upsilon (S(\ov{R}^1)\gamma ^2)=
\vartheta (S(\ov{R}^2)f^1\a _1)\ot
\upsilon (S(\ov{R}^1)f^2\a _2)\\
{\rm (\ref{ext})}&=&
\vartheta (f^2\ov{R}^2\a _1)\ot
\upsilon (f^1\ov{R}^1\a _2)=
\mathfrak{U}^{-1}\zeta (\a ),
\end{eqnarray*}
as needed. In a similar manner one can prove that
$\b _{(A\ot B)_{\mathfrak{F}}}=\zeta (\b )\mathfrak{U}$, the
details are left to the reader. Finally, for all $d\in D$ we have
\begin{eqnarray*}
\mathfrak{U}S_{A\ot B}(\zeta (d))\mathfrak{U}^{-1}&=&
\vartheta (R^2g^2S(d_1)f^2\ov{R}^2)\ot
\upsilon (R^1g^1S(d_2)f^1\ov{R}^1)\\
{\rm (\ref{ca}, \ref{qt3})}&=&\vartheta (S(d)_1)\ot
\upsilon (S(d)_2)=\zeta (S(d)).
\end{eqnarray*}
2) We are in the same hypothesis as in the Lemma \ref{lm5.2}, so
$\zeta : D\ra A\ot B$ is a right quasi-Hopf $B$-bimodule morphism.
As we have already explained before Lemma \ref{lm5.2}, the morphism $\zeta $
is bijective if and only if $\zeta _0$, the restriction of $\zeta $,
defines an isomorphism between $D^{\ov{co(B)}}$ and
$(A\ot B)^{\ov{co(B)}}$. But $D^{\ov{co(B)}}=D^{co(\upsilon )}$, so
if $d\in D^{\ov{co(B)}}$ then
\begin{eqnarray*}
\zeta (d)&=&\vartheta (d_1)\ot \upsilon (d_2)\\
&=&\vartheta (\kappa (x^1))\vartheta (d)
\vartheta (\kappa (S(x^3_2X^3)f^1))\ot x^2X^1\b S(x^3_1X^2)f^2,
\end{eqnarray*}
because of $\upsilon \circ \kappa =id_B$. Hence, by Lemma \ref{lm5.2},
$\zeta $ is bijective if and only if the map
\begin{eqnarray*}
&&\zeta _0: D^{co(\upsilon)}\ra \{\vartheta (\kappa (x^1))a\vartheta
(\kappa (S(x^3_2X^3)f^1))\ot x^2X^1\b S(x^3_1X^2)f^2\mid a\in A\},\\
&&\zeta _0(d)=\vartheta (\kappa (x^1))\vartheta (d)
\vartheta (\kappa (S(x^3_2X^3)f^1))\ot x^2X^1\b S(x^3_1X^2)f^2
\end{eqnarray*}
is bijective. Now, it follows that $\zeta $ is bijective if and only
if the restriction of $\vartheta $ defines a bijection
between $D^{co(\upsilon )}$ and $A$.
\end{proof}

We can now state the structure theorem of $D(H)$ when $H$ is
factorizable. The next result generalizes \cite[Theorem 4.3]{sch}.

\begin{theorem}\label{te5.4}
Let $(H, R)$ be a finite dimensional QT quasi-Hopf algebra,
$D(H)\rTo^{\pi }_{\widetilde{\pi }}H$ the quasi-Hopf algebra
morphisms defined by (\ref{pir}) and (\ref{tpir}), respectively,
and define $\zeta : D(H)\ra H\ot H$, given by
$\zeta ({\bf D})=\widetilde{\pi }({\bf D}_1)\ot \pi ({\bf D}_2)$,
for all ${\bf D}\in D(H)$, and
\begin{equation}\label{psptwist}
\mathbf{F}=Y^1_1x^1X^1y^1_1\ot Y^1_2x^2R^2X^3y^2
\ot Y^2x^3R^1X^2y^1_2\ot Y^3y^3,
\end{equation}
where $R^1\ot R^2$ is the $R$-matrix $R$ of $H$. Then the following
assertions hold:
\begin{itemize}
\item[1)] $\zeta : D(H)\ra (H\ot H)_{\mathbf{F}}^{\mathbf{U}}$ is a
quasi-Hopf algebra morphism, where
$\mathbf{U}:=\ov{R}^1g^2\ot \ov{R}^2g^1$.
\item[2)] $\zeta $ is bijective if and only if
$(H, R)$ is factorizable.
\end{itemize}
\end{theorem}
\begin{proof}
We consider in Proposition \ref{pr5.3} $D=D(H)$, $A=B=H$,
$\vartheta =\widetilde {\pi }$, $\upsilon =\pi $ and
$\kappa =i_D$. So the map $\zeta $ in the statement is the map
$\zeta $ in Proposition \ref{pr5.3} specialized for our case.
Moreover, from definition (\ref{madd}) of the $R$-matrix ${\cal R}$
of $D(H)$ we have
\begin{eqnarray*}
\pi ({\cal R}^1)\ot \widetilde{\pi }({\cal R}^2)&=&
\smi (p^2){}_iep^1_1\ot <{}^ie, q^2\ov{R}^2>q^1\ov{R}^1p^1_2\\
&=&\smi (p^2)q^2\ov{R}^2p^1_1\ot q^1\ov{R}^1p^1_2\\
{\rm (\ref{qt3}, \ref{pqr})}&=&
\smi (p^2)q^2p^1_2\ov{R}^2\ot q^1p^1_1\ov{R}^1
=\ov{R}^2\ot \ov{R}^1.
\end{eqnarray*}
Since $\pi $ and $\widetilde{\pi }$ are algebra maps
we obtain that $\pi (\ov{{\cal R}}^1)\ot \widetilde{\pi }
(\ov{{\cal R}}^2)=R^2\ot R^1$, so the twist (\ref{psptwist})
is the twist $\mathfrak{F}$ defined in (\ref{sptwist}) specialized
for our situation. Also, the element $\mathbf{U}$ is the element
$\mathfrak{U}$ defined in Proposition \ref{pr5.3} specialized for
our context and this prove the first assertion.\\
Applying again Proposition \ref{pr5.3} we have that $\zeta $ is
bijective if and only if the restriction of $\widetilde{\pi }$ provides
a bijection from $D(H)^{co(\pi )}$ to $H$. By Lemma \ref{lm5.1} this
is equivalent to $\ov{\cal Q}$ bijective. Finally, by Proposition
\ref{pr2.2} we obtain that $\zeta $ is bijective if and only if
$(H, R)$ is factorizable, and this finishes our proof.
\end{proof}

\section{Factorizable implies unimodular}
\setcounter{equation}{0}
In \cite{rad1} it is proved that a
finite dimensional factorizable Hopf algebra is unimodular. In
this Section we will show that this also holds for a finite
dimensional factorizable quasi-Hopf algebra. In particular, we
obtain that for any finite dimensional quasi-Hopf algebra $H$
its  Drinfeld double $D(H)$ is always a 
unimodular quasi-Hopf algebra.\\
Throughout, $H$ will be a finite dimensional quasi-Hopf algebra.
Recall that $t\in H$ is called a left (respectively right)
integral in $H$ if $ht=\va (h)t$ (respectively $th=\va (h)t$) for
all $h\in H$. We denote by $\int _l^H$ ($\int _r^H$) the space of
left (right) integrals in $H$. It follows from the bijectivity of
the antipode that $S(\int _l^H)=\int _r^H$ and $S(\int _r^H)=\int
_l^H$. If there is a non-zero left integral in $H$ which is at
the same time a right integral, then $H$ is called unimodular.
Hausser and Nill \cite{hn3} proved that for a finite dimensional
quasi-Hopf algebra the space of left or right integrals has
dimension $1$.\\
Let $t$ be a non-zero integral in $H$. Since the space of left
integrals is a two-sided ideal it follows from the uniqueness of
integrals in $H$ that there exists $\mu \in H^*$ such that
\begin{equation}\label{fu1}
th=\mu (h)t,~~\forall~~t\in \int _l^H~~\mbox{and}~~h\in H.
\end{equation}
It was noted in \cite{hn3} that $\mu $ is an element of $Alg(H,
k)$, i. e. $\mu $ is an algebra morphism from $H$ to $k$.
Moreover, $Alg(H, k)$ is a group with multiplication given by
$\varrho\circ \varsigma =(\varrho \ot \varsigma )\circ \Delta $,
unit $\va $, and inverse $\varrho ^{-1}=\varrho \circ S=\varrho
\circ S^{-1}$. Observe that $\mu =\va $ if and only if $H$ is
unimodular. As in the case of a Hopf algebra we will call $\mu $
the distinguished group-like element of $H^*$.\\
Hasser and Nill \cite{hn3} also introduced left cointegrals on a
finite dimensional quasi-Hopf algebra. These cointegrals are the
elements $\l $ of the dual space $H^*$ which satisfy for all
$h\in H$,
\begin{equation}\label{fu2}
\l (V^2h_2U^2)V^1h_1U^1=\mu (x^1)\l (hS(x^2))x^3.
\end{equation}
Here $U=U^1\ot U^2$ is the element defined by (\ref{eu}), $\mu $
is the distinguished group-like element of $H^*$, and if
$p_R=p^1\ot p^2$ and $f=f^1\ot f^2$ are the elements defined by
(\ref{qr}) and (\ref{f}), respectively, then  $V=V^1\ot V^2$ is
given by
\begin{equation}\label{ev}
V=\smi (f^2p^2)\ot \smi (f^1p^1).
\end{equation}
Using another structure theorem for right quasi-Hopf
$H$-bimodules, Hausser and Nill prove that the space of left
cointegrals ${\cal L}$ is one dimensional, and that the dual
paring ${\cal L}\ot \int _r^H\ni \l \ot r\mapsto <\l , r>\in k$
is non-degenerated. Let $\l $ be a non-zero left cointegral and
$r$ a non-zero right integral in $H$ such that $\l (r)=1$.
Following \cite{hn3}, we call
\begin{equation}\label{fu3}
\un{g}:=\l (V^1r_1U^1)V^2r_2U^2
\end{equation}
the comodulus of $H$. It was proved in \cite{hn3} that $\un{g}$
is invertible, and that its inverse is given by
\begin{equation}\label{fu4}
\un{g}^{-1}=\l (S(V^2r_2U^2))S^2(V^1r_1U^1).
\end{equation}
The results in the next two Lemmas also appear in a recent 
preprint of Kadison \cite{Kadison}. We prefer here to give direct 
proofs because they provide new formulas, which are of independent 
interest. \\
The following result expresses $\un{g}$ and $\un{g}^{-1}$ in terms 
of left integrals.
\begin{lemma}\label{lm6.1}
Let $H$ be a finite dimensional quasi-Hopf algebra, $\l $ a left
cointegral on $H$ and $0\not=r\in \int _r^H$ such that $\l (r
)=1$. If we set $r=\smi (t)$ for a certain left integral $t$ in
$H$, then
\begin{eqnarray}
&&\un{g}=\l (\smi (q^2t_2p^2))\smi (q^1t_1p^1),\label{fu5}\\
&&\un{g}^{-1}=\l (q^1t_1p^1)S(q^2t_2p^2),\label{fu6}
\end{eqnarray}
where $p_R=p^1\ot p^2$ and $q_R=q^1\ot q^2$ are the elements
defined by (\ref{qr}).
\end{lemma}
\begin{proof}
Let $q_L=\tqra \ot \tqrb $ be the element defined by (\ref{ql}).
We prove first that
\begin{equation}\label{fu7}
q^1t_1\ot q^2t_2=\tqra t_1\ot \tqrb t_2,
\end{equation}
for all $t\in \int _l^H$. To this end, we need the following
relations
\begin{eqnarray}
&&q_R=(\tqrb \ot 1)V\Delta (\smi (\tqra )),\label{fu8}\\
&&p_R=\Delta (S(\tpra ))U(\tprb \ot 1),\label{fu9}\\
&&U^1\ot U^2S(h)=\Delta (S(h_1))U(h_2\ot 1),~~\forall~~h\in
H,\label{fu10}
\end{eqnarray}
which can be found in \cite{hn3} (here $\tpra \ot \tprb $ is the
element $p_L$ defined by (\ref{ql})). Now,  $t\in \int _l^H$ and
(\ref{fu8}) imply that
\begin{equation}\label{fu11}
q^1t_1\ot q^2t_2=V^1t_1\ot V^2t_2.
\end{equation}
Together with a quasi-Hopf algebra $H=(H, \Delta , \va , \Phi ,
S, \a , \b )$ we also have $H^{cop}$ as quasi-Hopf algebra, where
¨cop¨ means opposite comultiplication. The quasi-Hopf algebra
structure is obtained by putting $\Phi _{cop}=(\Phi
^{-1})^{321}=x^3\ot x^2\ot x^1$, $S_{cop}=\smi $, $\a _{cop}=\smi
(\a )$ and $\b _{cop}=\smi (\b )$. It is not hard to see that in
$H^{cop}$ we have $(q_R)_{cop}=\tqrb \ot \tqra $,
$(p_R)_{cop}=\tprb \ot \tpra $ and $f_{cop}=(\smi \ot \smi )(f)$,
and therefore $V_{cop}=S(\tpra )f^2\ot S(\tprb )f^1$. Specializing
(\ref{fu11}) for $H^{cop}$, we obtain
$$
\tqra t_1\ot \tqrb t_2=S(\tprb )f^1t_1\ot S(\tpra )f^2t_2.
$$
On the other hand, one can easily check that (\ref{ql}, \ref{g2},
\ref{pf}) and $\smi (f^2)\b f^1=\smi (\a )$ imply
$$
S(\tprb )f^1\ot S(\tpra )f^2=q^1g^1_1\ot \smi (g^2)q^2q^1_2,
$$
where, as usual, we denote $f^{-1}=g^1\ot g^2$. From the above,
we conclude that
\begin{eqnarray*}
\tqra t_1\ot \tqrb t_2&=&q^1g^1_1t_1\ot \smi (g^2)q^2g^1_2t_2\\
\mbox{$(t\in \int _l^H)$}&=&q^1t_1\ot q^2t_2,
\end{eqnarray*}
as needed. We claim now that
\begin{equation}\label{fu12}
U^1\ot U^2=\tilde{q}^1_1p^1\ot \tilde{q}^1_2p^2S(\tqrb ).
\end{equation}
Indeed, by (\ref{fu9}) we have
\begin{eqnarray*}
\tilde{q}^1_1p^1\ot \tilde{q}^1_2p^2S(\tqrb )&=&[\tqra S(\tpra
)]_1U^1\tprb \ot [\tqra S(\tpra )]_2U^2S(\tqrb )\\
{\rm (\ref{fu10})}&=&[\tqra S(\tilde{q}^2_1\tpra
)]_1U^1\tilde{q}^2_2\tprb \ot [\tqra S(\tilde{q}^2_1\tpra
)]_2U^2\\
{\rm (\ref{pql})}&=&U^1\ot U^2.
\end{eqnarray*}
We write $p_R=p^1\ot p^2=P^1\ot P^2$, $f=f^1\ot f^2=F^1\ot F^2$
and $f^{-1}=g^1\ot g^2$. Then the above relations allow us to
compute
\begin{eqnarray*}
V^1r_1U^1\ot V^2r_2U^2&=&V^1r_1\tilde{q}^1_1p^1\ot
V^2r_2\tilde{q}^1_2p^2S(\tqrb )\\
\mbox{$(r\in \int _r^H$, \ref{ev})}&=&\smi (f^2P^2)\smi
(t)_1p^1\ot \smi (f^1P^1)\smi (t)_2p^2\\
{\rm (\ref{ca}, \ref{qr})}&=&\smi (S(x^1)f^2t_2P^2)\ot \smi
(S(x^2)f^1t_1P^1)\b S(x^3)\\
{\rm (\ref{g2}, \ref{pf})}&=&\smi (F^2x^3g^2_2t_2P^2)\ot \smi
(f^2F^1_2x^2g^2_1t_1P^1)\b f^1F^1_1x^1\\
\mbox{($\smi (f^2)\b f^1=\smi (\a ), \ref{q5}, t\in \int
_l^H$)}&=&\smi (x^3t_2P^2)\ot \smi (S(x^1)\a x^2t_1P^1)\\
{\rm (\ref{ql}, \ref{fu7})}&=&\smi (\tqrb t_2P^2)\ot \smi (\tqra
t_1P^1)=\smi (q^2t_2P^2)\ot \smi (q^1t_1P^1).
\end{eqnarray*}
Thus, we have proved that
$$
V^1r_1U^1\ot V^2r_2U^2=\smi (q^2t_2P^2)\ot \smi (q^1t_1P^1).
$$
It follows now that the above equality and (\ref{fu3}, \ref{fu4})
imply (\ref{fu5}) and (\ref{fu6}), so our proof is complete.
\end{proof}
Recall from \cite[Remarks 2.6]{bc1} that the map
$$
\ov{\theta }:\ H^*\ra H, ~~ \ov{\theta }(\chi )=\chi
(q^2t_2p^2)q^1t_1p^1~~\forall~~\chi \in H^*,
$$
is bijective. Thus there is an unique $\l \in H^*$ such that
\begin{equation}\label{fu13}
\l (q^2t_2p^2)q^1t_1p^1=1.
\end{equation}

\begin{lemma}\label{lm6.2}
The linear map $\l $ defined above is a non-zero left cointegral
on $H$.
\end{lemma}
\begin{proof}
The fact that $\l $ is non-zero follows from $\ov{\theta }(\l
)=1$. Let $\l _0$ be a non-zero left cointegral on $H$. Then
\begin{eqnarray*}
\ov{\theta }(\l _0)&=&\l _0(q^2t_2p^2)q^1t_1p^1\\
{\rm (\ref{fu8}, \ref{fu9})}&=&\l_0(V^2[\smi (\tqra )tS(\tpra
)]_2U^2)\tqrb V^1[\smi (\tqra )tS(\tpra )]_1U^1\tprb \\
\mbox{$(t\in \int _l^H, \ref{fu1})$}&=&\mu ^{-1}(\tpra ) \l
_0(V^2t_2U^2)V^1t_1U^1\tprb \\
{\rm (\ref{fu2})}&=&\mu (x^1)\mu ^{-1}(\tpra )\l
_0(tS(x^2))x^3\tprb \\
{\rm (\ref{fu1}, \ref{ql})}&=&\mu (x^1)\mu (X^1\b S(X^2))\mu
(S(x^2))\l _0(t)x^3X^3=\ov{\theta }(\mu (\b )\l _0(t)\l ).
\end{eqnarray*}
Since $\ov{\theta }$ is bijective we deduce that $\l _0=\mu (\b
)\l _0(t)\l $, and since $0\not=\l _0\in {\cal L}$, by the
uniqueness of left cointegrals on $H$ we conclude that $\l $ is a
non-zero left cointegral on $H$.
\end{proof}
We finally need the following result.

\begin{lemma}\label{lm6.3}
Let $H$ be a finite dimensional quasi-Hopf algebra, $t$ a
non-zero left integral in $H$ and $\mu $ the distinguished
group-like element of $H^*$. Then for any $h\in H$ the following
relations hold:
\begin{eqnarray}
&&q^1t_1\ot \smi (h)q^2t_2=hq^1t_1\ot q^2t_2,\label{fu14}\\
&&t_1\ot t_2=\b q^1t_1\ot q^2t_2=q^1t_1\ot \smi (\b
)q^2t_2,\label{fu15}\\
&&t_1p^1\ot t_2p^2S(h\lh \mu )=t_1p^1h\ot t_2p^2,\label{fu16}
\end{eqnarray}
where for all $h\in H$ and $\chi \in H^*$ we define 
$h\lh \chi =\chi (h_1)h_2$. 
\end{lemma}
\begin{proof}
The relations (\ref{fu14}, \ref{fu15}) are proved in \cite[Lemma
2.1]{bc1}. The equality (\ref{fu16}) follows from the following
computation
\begin{eqnarray*}
t_1p^1\ot t_2p^2S(h\lh \mu )&=&\mu (h_1)t_1p^1\ot t_2p^2S(h_2)\\
{\rm (\ref{fu1})}&=&t_1h_{(1, 1)}p^1\ot t_2h_{(1, 2)}p^2S(h_2)\\
{\rm (\ref{qr1})}&=&t_1p^1h\ot t_2p^2,
\end{eqnarray*}
for all $h\in H$, and this finishes the proof.
\end{proof}
We can now prove the main result of this Section.

\begin{theorem}\label{te6.4}
Let $(H, R)$ be a finite dimensional QT quasi-Hopf algebra and
$\mu $ the distinguished group-like element of $H^*$. Then the
following assertions hold.
\begin{itemize}
\item[1)] If $q_R=q^1\ot q^2=Q^1\ot Q^2$ and $p_R=p^1\ot
p^2=P^1\ot P^2$ are the elements defined by (\ref{qr}) then
\begin{equation}\label{fu17}
\mu (Q^1)q^2t_2p^2S(Q^2(R^2P^2\lh \mu ))R^1P^1\ot
q^1t_1p^1=S(u)q^1t_1p^1\ot q^2t_2p^2,
\end{equation}
where $R=R^1\ot R^2$ is the $R$-matrix of $H$ and $u$ is the
element defined in (\ref{elmu}).
\item[2)] If $(H, R)$ is factorizable then $H$ is unimodular.
\end{itemize}
\end{theorem}
\begin{proof}
1) Let us start by nothing that $g^1S(g^2\a )=\b $, (\ref{ext},
\ref{exta}) and (\ref{sqina}) imply
\begin{equation}\label{fu18}
R^1\b S(R^2)=S(\b u).
\end{equation}
Now, from (\ref{fu16}) we have
\begin{eqnarray*}
&&\hspace*{-4cm}\mu (Q^1)q^2t_2p^2S(Q^2(R^2P^2\lh \mu ))R^1P^1\ot
q^1t_1p^1\\
\hspace*{-4cm}&=&\mu (Q^1)q^2t_2p^2S(Q^2)R^1P^1\ot
q^1t_1p^1R^2P^2\\
\hspace*{-4cm}{\rm (\ref{fu1})}&=&q^2t_2Q^1_2p^2S(Q^2)R^1P^1\ot
q^1t_1Q^1_1p^1R^2P^2\\
\hspace*{-4cm}{\rm (\ref{pqr})}&=&q^2t_2R^1P^1\ot
q^1t_1R^2P^2\\
\hspace*{-4cm}{\rm (\ref{qt3})}&=&q^2R^1t_1P^1\ot
q^1R^2t_2P^2\\
\hspace*{-4cm}{\rm (\ref{fu15})}&=&q^2R^1\b Q^1t_1P^1\ot
q^1R^2Q^2t_2P^2\\
\hspace*{-4cm}{\rm (\ref{fu14})}&=&q^2R^1\b S(q^1R^2)Q^1t_1P^1\ot
Q^2t_2P^2\\
\hspace*{-4cm}{\rm (\ref{fu18})}&=&S(q^1\b u\smi
(q^2))Q^1t_1P^1\ot Q^2t_2P^2\\
\hspace*{-4cm}{\rm (\ref{sqina}, \ref{qr},
\ref{q6})}&=&S(u)Q^1t_1P^1\ot Q^2t_2P^2,
\end{eqnarray*}
and this proves the first assertion.\\
2) Let $\l \in H^*$ be the element defined by (\ref{fu13}). By
Lemma \ref{lm6.2} we know that $\l $ is a non-zero left
cointegral on $H$. Consider now $r$ a non-zero right integral in
$H$ such that $\l (r)=1$, and take $r=\smi (t)$ for some non-zero 
left integral $t$ in $H$. Then, by Lemma \ref{lm6.1} we have
$$
\smi (\un{g}^{-1})=\l (q^1t_1p^1)q^2t_2p^2.
$$
Applying $id\ot \l $ to the equality (\ref{fu17}) we obtain
$$
\mu (Q^1)\smi (\un{g}^{-1})S(Q^2(R^2P^2\lh \mu ))R^1P^1=S(u),
$$
and since $\smi (\un{g})S(u)=S(uS^{-2}(\un{g}))=S(\un{g}u)$, it
follows that the above relation is equivalent to
\begin{equation}\label{fu19}
\mu (Q^1)S(Q^2(R^2P^2\lh \mu ))R^1P^1=S(u)S(\un{g}).
\end{equation}
On the other hand, if we denote by $r^1\ot r^2$ another copy of
$R$, we then have
\begin{eqnarray*}
&&\hspace*{-3cm}\mu (Q^1)S(Q^2(R^2P^2\lh \mu ))R^1P^1\\
\hspace*{-2cm}{\rm (\ref{qr})}&=&\mu
(X^1R^2_1P^2_1)S(X^2R^2_2P^2_2)\a X^3R^1P^1\\
\hspace*{-2cm}{\rm (\ref{qt2})}&=&\mu
(X^1R^2y^2P^2_1)S(r^2X^3y^3P^2_2)\a
r^1X^2R^1y^1P^1\\
\hspace*{-2cm}{\rm (\ref{exta}, \ref{sqina}, \ref{qr})}&=&\mu
(q^1R^2y^2P^2_1)S(S(q^2)y^3P^2_2)uR^1y^1P^1\\
\hspace*{-2cm}{\rm (\ref{pr}, \ref{qt3})}&=&\mu (q^1X^1_{(1,
1)}p^1_1R^2P^2S(X^3)f^1)S(S(q^2)X^1_2p^2S(X^2)f^2)uX^1_{(1,
2)}p^1_2R^1P^1\\
\hspace*{-2cm}{\rm (\ref{sqina}, \ref{qr1})}&=&\mu
(X^1q^1p^1_1R^2P^2S(X^3)f^1)S(S(q^2p^1_2)p^2S(X^2)f^2)uR^1P^1\\
\hspace*{-2cm}{\rm (\ref{pqr}, \ref{sqina})}&=&\mu
(X^1R^2P^2S(X^3)f^1)u\smi (S(X^2)f^2)R^1P^1.
\end{eqnarray*}
From de above computation and (\ref{fu19}) we obtain
\begin{equation}\label{fu20}
\mu (X^1R^2P^2S(X^3)f^1)\smi
(S(X^2)f^2)R^1P^1=u^{-1}S(u)S(\un{g}).
\end{equation}
But, as we have already seen, if $(H, R)$ is QT then
$\tilde{R}=R^{-1}_{21}=\ov{R}^2\ot \ov{R}^1$ is another
$R$-matrix for $H$. Repeating the above computations for $(H,
\tilde{R})$ instead of $(H, R)$, we find that
\begin{equation}\label{fu21}
\mu (X^1\ov{r}^1P^2S(X^3)f^1)\smi
(S(X^2)f^2)\ov{r}^2P^1=\tilde{u}^{-1}S(\tilde{u})S(\un{g}),
\end{equation}
where we denote by $\tilde{u}$ the element defined as in
(\ref{elmu}) for $(H, \tilde{R})$ instead of $(H, R)$, and 
where $\ov{r}^1\ot \ov{r}^2$ is another copy of $R^{-1}$. More 
precisely, we have that
\begin{equation}\label{fu22}
\tilde{u}=S(u^{-1}).
\end{equation}
Indeed, one can easily see that (\ref{fu18}) and (\ref{sqina}) 
imply
\begin{equation}\label{fu23}
\ov{r}^2\b S(\ov{r}^1)=\smi (\b )u^{-1}=u^{-1}S(\b ).
\end{equation}
Now, we compute
\begin{eqnarray*}
\tilde{u}&=&S(\ov{r}^1x^2\b S(x^3))\a \ov{r}^2x^1\\
\mbox{$(S(\b f^1)f^2=\a )$}&=&S(\b f^1\ov{r}^1x^2\b
S(x^3))f^2\ov{r}^2x^1\\
{\rm (\ref{ext}, \ref{qr})}&=&S(\ov{r}^2\b
S(\ov{r}^1)f^2p^2)f^1p^1\\
{\rm (\ref{fu23})}&=&S(\smi (f^1p^1)u^{-1}S(\b )f^2p^2)\\
\mbox{(\ref{sqina}, $S(\b f^1)f^2=\a $, \ref{qr}, \ref{q6})}&=&S(u^{-1}S(p^1)\a p^2)
=S(u^{-1}).
\end{eqnarray*}
Now, since $S^2(u)=u$ the relation (\ref{fu21}) becomes 
\begin{equation}\label{fu24}
\mu (X^1\ov{r}^1P^2S(X^3)f^1)\smi
(S(X^2)f^2)\ov{r}^2P^1=S(u)u^{-1}S(\un{g}).
\end{equation}
From (\ref{sqina}) it follows that $uS^{-1}(u)=S(u)u$, and since $S^2(u)=u$ 
we conclude that $uS(u)=S(u)u$, so $u^{-1}S(u)=S(u)u^{-1}$. Hence, by 
(\ref{fu20}) and (\ref{fu24}) we obtain 
$$
\mu (X^1R^2P^2S(X^3)f^1)\smi
(S(X^2)f^2)R^1P^1=\mu (X^1\ov{r}^1P^2S(X^3)f^1)\smi
(S(X^2)f^2)\ov{r}^2P^1.
$$
This comes out explicitly as 
$$
\mu (R^2P^2)R^1P^1=\mu (\ov{r}^1P^2)\ov{r}^2P^1 
$$
and implies 
$$
\mu (Q^1_1R^2P^2S(Q^2))Q^1_2R^1P^1=
\mu (Q^1_1\ov{r}^1P^2S(Q^2))Q^1_2\ov{r}^2P^1.
$$ 
From (\ref{qt3}) and (\ref{pqr}) we deduce that 
\begin{equation}\label{fu25}
\mu (R^2)R^1=\mu (\ov{r}^1)\ov{r}^2\Leftrightarrow 
\mu (R^2r^1)R^1r^2=1.
\end{equation}
Finally, the above relation allows us to compute 
\begin{eqnarray*}
{\cal Q}(\mu )&=&\mu (\tqra X^1R^2r^1p^1)
\tilde{q}^2_1X^2R^1r^2p^2S(\tilde{q}^2_2X^3)\\
{\rm (\ref{fu25}, \ref{qr})}&=&\mu (\tqra X^1x^1)\tilde{q}^2_1
X^2x^2\b S(\tilde{q}^2_2X^3x^3)\\
{\rm (\ref{q5}, \ref{ql})}&=&\mu (\a )\b ={\cal Q}(\mu (\a )\va ). 
\end{eqnarray*}  
If $(H, R)$ is factorizable then ${\cal Q}$ is bijective, so 
$\mu =\mu (\a )\va $. In particular, 
$1=\mu (1)=\mu (\a )\va (1)=\mu (\a)$. Hence $\mu =\va $, and this 
means that $H$ is unimodular.
\end{proof}

\begin{theorem}\label{te6.5}
Let $H$ be a finite dimensional quasi-Hopf algebra. Then the 
Drinfeld double $D(H)$ of $H$ is a unimodular quasi-Hopf algebra. 
\end{theorem}
\begin{proof}
It is an immediate consequence of Proposition \ref{pr2.3} 
and Theorem \ref{te6.4}. 
\end{proof}

\begin{center}
{\bf ACKNOWLEDGEMENTS}
\end{center}
The authors wish to express their sincere gratitude to Claudia 
Menini and Alessandro Ardizzoni for various comments during the 
development of this paper.\\
All the diagrams were made by using of "diagrams" software
program of Paul Taylor.


\end{document}